\documentstyle[12pt]{article}
\let\Bbb\relax
\newfont{\Bb }{msbm10 scaled 1000}
\newfont{\Bbb}{msbm10 scaled 1200}
\newfam\eufam%
\font\euzw=eufm10 scaled 1200%
\font\euac=eufm9%
\textfont\eufam=\euzw\scriptfont\eufam=\euac%
\def\fr{\fam\eufam\euzw}%
\newcommand{\Z}{{\Bbb Z}}
\newcommand{\R}{{\Bbb R}}
\newcommand{\C}{{\Bbb C}}

\newcommand{\Ha}{{\Bbb H}}
\newcommand{\J}{{\Bbb J}}
\newtheorem{thm}{Theorem}
\newtheorem{prop}{Proposition}
\newtheorem{cor}{Corollary}
\newtheorem{lemma}{Lemma}
\newtheorem{dof}{Definition}
\begin{document}
\begin{center}
{\LARGE\bf A New Construction of\\[0.5em]
Homogeneous Quaternionic Manifolds\\[0.5em]
and Related Geometric Structures}
\end{center}
\vskip 2.0 true cm 
\begin{center}
{\Large\bf  Vicente Cort\'es\footnote{Supported by SFB 256 
(Bonn University). e-mail: 
vicente@math.uni-bonn.de}}\end{center}
\noindent
{\small
Mathematisches Institut der Universit\"at Bonn,
Beringstr.\ 1,
53115 Bonn, Germany\\
}

\begin{abstract} 
Let $V = \mbox{\R}^{p,q}$ be the pseudo-Euclidean 
vector space of signature $(p,q)$, $p\ge 3$ and $W$ a module
over the even Clifford algebra $C\! \ell^0 (V)$. A homogeneous 
quaternionic manifold $(M,Q)$ is constructed for any 
${\fr spin}(V)$-equivariant linear map $\Pi : \wedge^2 W \rightarrow V$. 
If the  skew symmetric vector valued bilinear form  $\Pi$ is nondegenerate 
then $(M,Q)$ is endowed with a canonical pseudo-Riemannian metric $g$ such
that $(M,Q,g)$ is a homogeneous  quaternionic pseudo-K\"ahler manifold.  
If the metric $g$ is positive definite, i.e.\ a Riemannian metric, 
then the quaternionic K\"ahler manifold $(M,Q,g)$ is shown to admit a simply 
transitive solvable group of automorphisms. In this special case ($p=3$) 
we recover all the known homogeneous  
quaternionic K\"ahler manifolds of negative scalar curvature 
(Alekseevsky spaces) in a unified and direct way.
If $p>3$ then $M$ does not admit any transitive
action of a solvable Lie group and we obtain new families of  
quaternionic pseudo-K\"ahler manifolds. Then it is shown that for
$q = 0$ the noncompact quaternionic manifold $(M,Q)$ can be endowed with a 
Riemannian metric $h$ such that $(M,Q,h)$ is a homogeneous quaternionic 
Hermitian manifold, which  does not admit any  transitive solvable 
group of isometries if $p>3$. 

The twistor bundle $Z \rightarrow M$ 
and the canonical ${\rm SO}(3)$-principal bundle $S \rightarrow M$ 
associated to the quaternionic manifold $(M,Q)$ are shown to be 
homogeneous under the automorphism group of the base. More specifically, 
the twistor space is a homogeneous complex manifold carrying an
invariant holomorphic distribution $\cal D$ of complex codimension one, 
which is a complex contact structure if and only if $\Pi$ is nondegenerate. 
Moreover, an equivariant open holomorphic immersion $Z \rightarrow
\bar{Z}$ into a homogeneous complex manifold $\bar{Z}$ of complex algebraic
group is constructed. 

Finally, the construction is shown to have a natural mirror
in the category of supermanifolds. In fact, for any 
${\fr spin}(V)$-equivariant linear map $\Pi : \vee^2 W \rightarrow V$
a homogeneous quaternionic supermanifold $(M,Q)$ is constructed and,
moreover, a homogeneous quaternionic pseudo-K\"ahler supermanifold
$(M,Q,g)$ if the symmetric vector valued bilinear form $\Pi$ is nondegenerate.


   \end{abstract} 
\tableofcontents 
\section*{Introduction}
Let us start this introduction by recalling the notion of 
quaternionic manifold, see \cite{A-M2}.
A hypercomplex structure on a real vector space $E$ consists of 3
complex structures $(J_1,J_2,J_3)$ on $E$ satisfying $J_1J_2 = J_3$. It 
defines on $E$ the structure of (left-) vector space over the quaternions
$\mbox{\Ha} = \{ 1,i,j,k \}$ such that multiplication by $i$, $j$ and $k$
is given, respectively, by $J_1$, $J_2$ and $J_3$. The 3-dimensional
subspace $Q = {\rm span} \{ J_1,J_2,J_3\} \subset {\rm End}(E)$ is what
is called a quaternionic structure on $E$. A Euclidean scalar product
$\langle \cdot ,\cdot \rangle$ on $(E,Q)$ is called ($Q$-) Hermitian if
$Q$ consists of skew symmetric endomorphisms of $(E, 
\langle \cdot ,\cdot \rangle )$. Now let $M$ be a smooth manifold, 
$\dim M >4$. An almost quaternionic structure $Q$ on $M$ is a smooth
field $m \mapsto Q_m$ whose value at $m \in M$ is a quaternionic structure
on $T_mM$. $Q$ is called a quaternionic structure and $(M,Q)$
a {\bf quaternionic manifold} if there exists a torsionfree connection on $TM$
preserving the rank 3 subbundle $Q \subset {\rm End}(TM)$. 
Now let $g$ be a Riemannian metric on $M$, Hermitian with respect
to an (almost) quaternionic structure $Q$. Then $M$ with the 
structure $(Q,g)$ is called an (almost) quaternionic Hermitian
manifold. If, moreover, the Levi-Civita connection preserves $Q$ then
$(M,Q,g)$ is said to be a {\bf quaternionic K\"ahler manifold}. 
We remark that quaternionic K\"ahler manifolds represent one of the few
basic Riemannian geometries, as defined by Berger's list of possible 
Riemannian holonomy groups, see \cite{A1}, \cite{Bes}, 
\cite{Br1} and \cite{S2}. For the possible holonomy groups 
of, not necessarily Riemannian, 
torsionfree connections see \cite{Br2}, \cite{Schw} and references
therein.  

Next we review what is known about homogeneous quaternionic K\"ahler
manifolds. First of all,  quaternionic K\"ahler
manifolds $(M,Q,g)$ are Einstein manifolds, i.e.\ $Ric = cg$, 
see  e.g.\ \cite{A1}, \cite{Bes}, \cite{S2} and \cite{Kak}.
We discuss 3 cases depending on the sign of the constant
$c$ (which is the sign of the scalar curvature). 

$c=0$) Ricci-flat quaternionic K\"ahler
manifolds are better known as hyper-K\"ahler
manifolds, see \cite{Bes} Ch.\ 14. They are 
K\"ahler manifolds with respect to 3 complex structures 
$J_1$, $J_2$ and $J_3 = J_1 J_2$. It is a general fact
that any homogeneous Ricci-flat Riemannian manifold is
necessarily flat, see \cite{A-K}. In particular, all
homogeneous hyper-K\"ahler manifolds are flat. 
{}For  hyper-K\"ahler manifolds of small cohomogeneity
see \cite{Bi1}, \cite{Bi2}, \cite{B-G}, \cite{D-S1}, \cite{D-S2}, 
\cite{K-S2}, 
\cite{Sw1}, \cite{Sw2}
and \cite{C4}.

$c>0$) It follows from Myer's theorem that any complete Einstein
manifold of positive scalar curvature is necessarily compact. In particular,
any complete quaternionic K\"ahler manifold of positive scalar curvature
is compact. It was proven in \cite{L-S} that for every $n>1$ 
there is only a finite number of such manifolds of dimension $4n$ 
up to homothety, cf.\  
\cite{Bea}, \cite{G-S}, \cite{L1}, \cite{L4}, \cite{L5}, \cite{P-S}.  
The only known examples are, up to now, the Wolf spaces \cite{W1}. 
These are precisely the homogeneous quaternionic K\"ahler manifolds
of positive scalar curvature and are all symmetric of compact type,
cf.\ \cite{A2}. More generally, the Wolf spaces can be characterized 
as the compact quaternionic K\"ahler manifolds 
which admit an action of cohomogeneity $\le 1$ by a compact 
semisimple group of isometries and which are not scalar-flat, 
see \cite{D-S3}, cf.\ \cite{A-P}. 

$c<0$) Complete noncompact quaternionic K\"ahler manifolds
of negative scalar curvature exist in abundance. In fact, it was
proven in \cite{L3} that the moduli space of complete
quaternionic K\"ahler metrics on $\mbox{\R}^{4n}$, $n>1$, is infinite
dimensional. For explicit constructions of, in general not complete,
quaternionic K\"ahler manifolds see e.g.\ \cite{G1}, \cite{G-L}
and \cite{L2}.  

What about homogeneous examples? First of all, the noncompact duals
of the Wolf spaces are symmetric (and hence homogeneous) 
quaternionic K\"ahler manifolds
of negative Ricci and nonpositive sectional curvature. 
Moreover, like any Riemannian symmetric space of noncompact type, 
these manifolds admit smooth quotients by discrete cocompact
groups of isometries, see \cite{Bo}, \cite{Fi-S}. 
The first examples of quaternionic K\"ahler manifolds which are
not locally symmetric were found by D.V.\ Alekseevsky in 
\cite{A3}. An Alekseevsky space is a homogeneous quaternionic 
K\"ahler manifold of negative scalar curvature 
which admits a simply transitive splittable
solvable group of isometries. It follows from  
Iwasawa's decomposition theorem that the noncompact
duals of the Wolf spaces are precisely the symmetric
Alekseevsky spaces. Besides these there are 3 series
of nonsymmetric Alekseevsky spaces, see \cite{A3}, \cite{dW-VP2}  
and \cite{C2}.   In \cite{A3} it was conjectured that any
noncompact homogeneous quaternionic K\"ahler manifold admits a
transitive solvable group of isometries. This conjecture
is still open: up to now, the only known examples of homogeneous
quaternionic K\"ahler manifolds of negative scalar curvature
are the Alekseevsky spaces. However, by the construction presented in this
work we obtain many homogeneous quaternionic pseudo-K\"ahler manifolds (with 
indefinite metric) which do not admit any transitive action of a solvable
Lie group. Moreover, in \ref{notransSec} we 
construct a family of noncompact homogeneous quaternionic Hermitian
manifolds (with positive definite metric) with no  
transitive solvable group of isometries. 

It is natural to ask for examples of compact locally homogeneous 
quaternionic K\"ahler manifolds. The following negative result was proven
in \cite{A-C4}. Let $M$ be a compact quaternionic K\"ahler
manifold or, more generally, a quaternionic K\"ahler manifold of
finite volume. If the universal cover $\tilde{M}$ is a 
homogeneous quaternionic K\"ahler manifold then it is necessarily
symmetric. In particular, the only Alekseevsky spaces which admit
smooth quotients of finite volume by discrete groups of isometries
are the symmetric ones, this was as well proven in \cite{A-C1} by a simpler 
method. Additionally, the symmetric Alekseevsky spaces can 
be characterized by the property of having nonpositive 
curvature, see \cite{C2}. 

Given a simply transitive Lie group $L$ of isometries acting on 
a Riemannian manifold $M$, there exists an algorithm to compute
the full isometry group of $M$ \cite{A-W}, cf.\ \cite{Wo}.  However, this algorithm
involves the covariant derivatives (of all orders) 
of the curvature tensor and hence can only be applied effectively
in very simple situations.  If the simply transitive  
group $L$ is splittable solvable and unimodular then the
full isometry group is easily computed, see \cite{G-W}. 
Unfortunately, the splittable solvable groups of isometries
acting simply transitively on the Alekseevsky spaces are not 
unimodular. For that reason in \cite{A-C1} a new
algorithm was developed which completely avoids the 
curvature tensor and works also for nonunimodular
splittable solvable groups.  Using it the full isometry group of the 
nonsymmetric Alekseevsky spaces was determined  \cite{A-C1}. 
The Lie algebra of the full isometry group was previously described  
by the theoretical physicists de Wit, Vanderseypen and Van Proeyen 
by a different method, see \cite{dW-V-VP}.   

In the following, we shortly comment on the attention paid to
Alekseevsky's spaces in the physical literature. 
There is a concept of special geometry, which evolved in the 
theory of strings and supergravity, see e.g.\ \cite{Z}, \cite{B-W},
\cite{dW-VP1}, \cite{G-S-T}. More specifically, special
K\"ahler geometry is the geometry associated to $N=2$ supergravity
in $d=4$ space time dimensions coupled to vector multiplets
and was first described in \cite{dW-VP1}, cf.\ \cite{C-D-F}, \cite{St}, 
\cite{C3} and \cite{Fr}.  The (Kuranishi) moduli space
of a Calabi-Yau 3-fold bears this particular geometry. 
Moreover, there is a construction (the ``c-map''), 
related to Mirror Symmetry, which 
to any special K\"ahler manifold associates a quaternionic K\"ahler
manifold, see \cite{C-F-G} and \cite{F-S}; for the general framework
see \cite{Cr}.   Using the c-map, Cecotti \cite{Ce}
was the first to relate the classification problem for 
homogeneous special K\"ahler manifolds to Alekseevsky's classification 
\cite{A3}. In addition, he introduced Vinberg's theory of 
$T$-algebras \cite{V2} to describe the first nonsymmetric homogeneous special 
K\"ahler manifolds (the symmetric special K\"ahler manifolds were
described in terms of Jordan algebras, see \cite{C-VP} and \cite{G-S-T}).
Cecotti's classification of homogeneous special K\"ahler manifolds 
was extended in \cite{dW-VP2}, \cite{C3} and \cite{A-C5}. 
In \cite{C4}, the hyper-K\"ahlerian version of the c-map was used to
construct a  natural (pseudo-) hyper-K\"ahler structure on the bundle
of intermediate Jacobians over the moduli space of gauged Calabi-Yau
3-folds. 

In the last part of the introduction 
we describe the main components, results and the 
global structure of the present paper. The basic algebraic
data of our construction are a pseudo-Euclidean vector space $V$,
a module $W$ over the even Clifford algebra $C\! \ell^0 (V)$ 
and a ${\fr spin} (V)$-equivariant linear map 
$\Pi : \wedge^2 W \rightarrow V$.  
Any such $\Pi$ defines a 
$\mbox{\Z}_2$-graded Lie algebra ${\fr p} = {\fr p}(\Pi ) 
= {\fr p}_0 + {\fr p}_1$, where ${\fr p}_0 = {\rm Lie}\, {\rm Isom}(V)
= {\fr o}(V) + V$ and ${\fr p}_1 = W$. Here $V$ acts trivially on $W$
and ${\fr o}(V) \cong {\fr spin}(V)$ acts via the inclusion
${\fr spin}(V) \subset C\! \ell^0 (V)$ on the $C\! \ell^0 (V)$-module
$W$. The Lie bracket on ${\fr p}_1 \times {\fr p}_1$ is given by $\Pi$. 
The Lie algebras ${\fr p}(\Pi )$ were introduced in \cite{A-C2}, where
a basis for the vector space of ${\fr spin}(V)$-equivariant linear
maps $\Pi : \wedge^2 W \rightarrow V$ was explicitly constructed. 
The Lie algebra ${\fr p}(\Pi )$ is called an extended Poincar\'e algebra
of signature $(p,q)$ if $V \cong \mbox{\R}^{p,q}$ has signature
$(p,q)$. A mirror symmetric version of this construction
is obtained replacing $\Pi : \wedge^2 W \rightarrow V$ by a 
${\fr spin}(V)$-equivariant linear map $\Pi : \vee^2 W \rightarrow V$.
Here $\vee^2 W = {\rm Sym}^2W$ denotes the symmetric square of $W$. 
The corresponding algebraic structure ${\fr p}(\Pi ) 
= {\fr p}_0 + {\fr p}_1$ is now a super Lie algebra. It is called a 
superextended Poincar\'e algebra of signature $(p,q)$. For special signatures
$(p,q)$ of space time $V$ these super Lie algebras play an important role
in the physical literature since the early days of supersymmetry and 
supergravity, see e.g.\ \cite{Go-L} and \cite{O-S}; for more
recent contributions see e.g.\ \cite{F} and references therein. 
Notice that if $(p,q) = (1,3)$  then $V = \mbox{\R}^{1,3}$ is Minkowski
space and the even subalgebra ${\fr p}_0 = 
{\rm Lie}\, {\rm Isom}(\mbox{\R}^{1,3}) \subset {\fr p}(\Pi )$ is the
classical Poincar\'e algebra. The construction of all superextended
Poincar\'e algebras (of arbitrary signature) was carried out in 
\cite{A-C2}. In \cite{A-C-D-S} the twistor equation is interpreted
as the differential equation satisfied by infinitesimal automorphisms
of a geometric structure modelled on the linear Lie supergroup
associated to a superextended Poincar\'e algebra (for more on the 
twistor equation see e.g\ \cite{K-R} and references therein). In contrast
with the case of superextended Poincar\'e algebras, to our knowledge,
extended Poincar\'e algebras do not occur in the physical literature before the
publication of \cite{A-C2}: The first occurrence is \cite{D-L}. 

Any extended Poincar\'e algebra ${\fr p}(\Pi )$, as above, admits a derivation
$D$ with eigenspace 
decomposition ${\fr p}(\Pi ) = {\fr o}(V) + V + W$  and corresponding
eigenvalues $(0,1,1/2)$. We can extend ${\fr p}(\Pi)$ by $D$ obtaining 
a new  $\mbox{\Z}_2$-graded Lie algebra ${\fr g} = {\fr g}(\Pi) =
{\fr g}_0 + {\fr g}_1$, where ${\fr g}_0 = \mbox{\R}D + {\fr p}_0 =
\mbox{\R}D + {\fr o}(V) + V$ and
${\fr g}_1 = {\fr p}_1 = W$. The adjoint representation of the Lie 
algebra $\fr g$ is faithful and hence defines on $\fr g$ the structure
of linear Lie algebra. Let $G = G(\Pi ) \subset {\rm Aut}\, {\fr g}$ denote
the corresponding connected linear Lie group. It is the numerator
of the homogeneous quaternionic manifolds $M = M(\Pi ) = G/K$ we
are going to construct. To define the denominator $K$ let 
$E \subset V$ be a 3-dimensional Euclidean subspace and 
${\fr k} = {\fr k}(E) \subset {\fr o}(V) \subset {\fr g}$ the maximal 
subalgebra which preserves $E$, i.e.\ 
${\fr k} = {\fr o}(E) \oplus {\fr o}(E^{\perp})$. Now we define 
$K = K(E) \subset G$ to be the connected linear Lie group with
Lie algebra $\fr k$. Our main theorem is the following, see 
Thm.\ \ref{mainThm}:

\bigskip\noindent
{\bf Theorem A} \begin{it} 
Let $V$ be any pseudo-Euclidean 
vector space, $E \subset V$ a Euclidean $3$-dimensional subspace, $W$
any $C\! \ell^0 (V)$-module and $\Pi : \wedge^2W \rightarrow V$ any 
${\fr spin}(V)$-equivariant linear map. Let $M = G/K$ be the homogeneous
manifold associated to the Lie groups $G = G(\Pi )$ and $K = K(E) \subset
G$ constructed above. Then the following is true:  
\begin{enumerate}
\item[1)] There exists a $G$-invariant quaternionic structure $Q$ on $M$.
\item[2)] If $\Pi$ is nondegenerate (i.e.\ if $W\ni s\mapsto \Pi (s\wedge
\cdot ) \in W^{\ast} \otimes V$ is injective) then there exists a 
$G$-invariant pseudo-Riemannian metric $g$ on $M$ such that
$(M,Q,g)$ is a homogeneous quaternionic pseudo-K\"ahler manifold.  
\end{enumerate}  
\end{it}

\bigskip\noindent  
We remark that for $V = \mbox{\R}^{3,q}$ and $W$ arbitrary the map 
$\Pi$ can always
be chosen such that the metric $g$ in 2)
becomes positive definite.  In this special case we recover the 
3 series of Alekseevsky spaces
by a simple and unified construction, which completely avoids
the technicalities of constructing  complicated representations of 
K\"ahlerian Lie algebras, see \cite{A3} and \cite{C2}.  

If in all constructions $\Pi  : \wedge^2W \rightarrow V$ is replaced by a 
${\fr spin}(V)$-equivariant linear map  $\Pi  : \vee^2W \rightarrow V$
then we obtain the following analogue of Theorem A in the 
category of supermanifolds, see
Thm.\ \ref{supermainThm}:

\bigskip\noindent
{\bf Theorem B} \begin{it} 
Let $V$ be any pseudo-Euclidean 
vector space, $E \subset V$ a Euclidean $3$-dimensional subspace, $W$
any $C\! \ell^0 (V)$-module and $\Pi : \vee^2W \rightarrow V$ any 
${\fr spin}(V)$-equivariant linear map. Let $M = G/K$ be the homogeneous
supermanifold associated to the 
Lie supergroups $G = G(\Pi )$ and $K = K(E) \subset
G$ constructed in \ref{superSec}. Then the following holds: 
\begin{enumerate}
\item[1)] There exists a $G$-invariant quaternionic structure $Q$ on $M$.
\item[2)] If $\Pi$ is nondegenerate (i.e.\ if $W\ni s\mapsto \Pi (s\vee
\cdot ) \in W^{\ast} \otimes V$ is injective) then there exists a 
$G$-invariant pseudo-Riemannian metric $g$ on $M$ such that
$(M,Q,g)$ is a homogeneous quaternionic pseudo-K\"ahler supermanifold.  
\end{enumerate}  
(For the definition of quaternionic structure, pseudo-Riemannian metric etc.\  
on a supermanifold see the appendix.) 
\end{it}

\bigskip\noindent  
Furthermore, we remark  that replacing Euclidean 3-space 
$E \cong \mbox{R}^{3,0}$ by Lorentzian 3-space
$E \cong \mbox{\R}^{1,2}$ one obtains a para-quaternionic
version of our construction. 

Finally, we outline the structure of the paper:  
In section \ref{EPASec} we discuss extended Poincar\'e algebras. 
The basic definitions and formulas are given in \ref{BFSec}. 
To our fundamental map $\Pi : \wedge^2W \rightarrow V$
and to an oriented Euclidean subspace $E \subset V$, 
$\dim E \equiv 3\pmod{4}$, 
we associate a canonical symmetric bilinear form $b$ on 
$W$ and study its properties in \ref{CSBFSec}. 
Using the form $b$, in \ref{isomSec} we classify extended Poincar\'e
algebras of signature $(p,q)$, $p \equiv 3\pmod{4}$, up to isomorphism.

In section \ref{HQMSec} we construct the homogeneous quaternionic manifolds 
of Theorem A. The basic notions of quaternionic geometry
are recalled  in \ref{QSSec}. First, see \ref{homogSec}, we
describe the structure of the Lie group $G = G(\Pi )$ and the
coset spaces $M = G/K$, $K = K(E)$, where $E \subset V$ is any
pseudo-Euclidean subspace. The proof of Theorem A is given in 
\ref{MainSec}. A crucial observation is that $M = G/K$ contains
the locally symmetric quaternionic pseudo-K\"ahler submanifold
$M_0 = G_0/K$, where $G_0 \subset G$ is the connected
linear Lie group associated to ${\fr g}_0 \subset {\fr g} = {\fr g}_0 +
{\fr g}_1$. The first step consists in extending the $G_0$-invariant
quaternionic and pseudo-Riemannian structures on $M_0$ to 
$G$-invariant structures on $M$. Using the canonical symmetric
bilinear form $b$ on $W$ introduced in \ref{CSBFSec} the 
pseudo-Riemannian metric is 
extended, in the nondegenerate case, to a $G$-invariant 
pseudo-Riemannian metric $g$ on $M$.  
The quaternionic structure is always extended by the 
beautifully simple formula (\ref{extEqu}) to a $G$-invariant almost 
quaternionic structure $Q$ on $M$. To prove that $Q$ is a quaternionic
structure, we construct a $G$-invariant torsionfree connection $\nabla$
on $M$ which preserves $Q$. In the nondegenerate case $\nabla$ is simply
the Levi-Civita connection of the pseudo-Riemannian metric $g$. 
We use the description of invariant connections on homogeneous manifolds in 
terms of Nomizu maps, see \ref{ICHMSec}. 

Then we 
concentrate on  the Riemannian case, see \ref{RCSec}. The 
Riemannian manifolds $M(\Pi )$ are  classified up to isometry using results 
of \ref{isomSec}. We show that all these manifolds  admit
a non-Abelian simply transitive splittable solvable 
group of isometries and hence are Alekseevsky spaces, see Thm.\ \ref{PDThm}. 
Moreover, we explain how to obtain the 3 series of Alekseevsky spaces
by specifying $V$, $W$, and $\Pi$, see Thm.\ \ref{AThm}.

The quaternionic pseudo-K\"ahler manifolds $(M(\Pi ),Q,g)$  
of Theorem A do not admit 
any transitive action of a solvable group if $V = \mbox{\R}^{p,q}$ with
$p>3$, see Thm.\ \ref{NTThm}. Moreover, if $q = 0$, then we can replace
$g$ by a $G$-invariant Riemannian and $Q$-Hermitian metric $h$ such 
that  $(M(\Pi ),h,g)$ are homogeneous quaternionic Hermitian 
manifolds with no transitive
solvable group of isometries, see Thm.\ \ref{NTRThm}. 

In section \ref{bundleSec} we study 
the various bundles associated to the quaternionic
manifolds $M = M(\Pi) = G/K$: The twistor bundle $Z(M)$, the canonical 
${\rm SO}(3)$-principal bundle $S(M)$ and the  Swann bundle $U(M)$. 
We show that $G$ acts transitively on $Z(M)$ and $S(M)$ and with 
cohomogeneity one on $U(M)$. In particular, $Z = Z(M)$ is a homogeneous
complex manifold of the group $G = G(\Pi )$. We exhibit a $G$-invariant 
holomorphic tangent hyperplane distribution $\cal D$ on $Z$ and prove that
$\cal D$ defines a complex contact structure on $Z$ if and only if $\Pi$
is nondegenerate. Moreover, we construct an open $G$-equivariant holomorphic
immersion $Z \rightarrow \bar{Z}$ of $Z$ into a  
homogeneous complex manifold of the complex algebraic
group $G^{\mbox{\C}} \subset {\rm Aut}({\fr g}^{\mbox{\C}})$, see
Thm.\ \ref{EMBThm}. This immersion is a finite covering over an open
$G$-orbit. In the nondegenerate case $\bar{Z}$ is a 
homogeneous complex contact manifold of the group $G^{\mbox{\C}}$
and the immersion is a morphism of complex contact manifolds.  

In the final section \ref{superSec} we extend our construction to the category of 
supermanifolds, proving Theorem B. We have aimed at a straightforward
presentation, summarizing the needed supergeometric background in the  
appendix. 

\bigskip\noindent
{\bf Acknowledgements} 
It is a pleasure to thank my friend and coauthor Dmitry V.\ Alekseevsky
for many intensive discussions during our productive collaboration.
Also I am very grateful to Werner Ballmann and Ursula Hamenst\"adt for 
their encouragement and support. For reading the manuscript 
I thank Gregor 
Weingart.  Finally, I am indebted to many collegues for
hospitality and conversations, especially to Oliver Baues, Robert 
Bryant, Shiing-Shen Chern, Jost Eschenburg, Phillip A.\ Griffiths,
Ernst Heintze, Alan Huckleberry, Enrique Mac\'{\i}as, Robert Osserman,
Hans-Bert Rademacher, Gudlaugur Thorbergsson and Joseph A.\ Wolf. 

My research has received important support from the following institutions:
SFB 256 (Bonn University), Alexander von Humboldt Foundation and
Mathematical Sciences Research Institute (Berkeley). 
\section{Extended Poincar\'e algebras} 
\label{EPASec} 
\subsection{Basic facts} 
\label{BFSec} 
Let $V$ be a pseudo-Euclidean vector space with scalar product
$\langle \cdot , \cdot \rangle$. There exists an orthonormal
basis $(e_i)$, $i = 1, \ldots ,n = \dim V = p + q$,  of $V$
such that $\langle x,x \rangle = \sum_{i=1}^p (x^i)^2 -
\sum_{j=p+1}^n(x^j)^2$ for all $x = \sum_{i=1}^n x^ie_i \in V$.
Any such basis defines an isometry between $V$ and the standard
pseudo-Euclidean vector space $\mbox{\R}^{p,q}$ of signature
$(p,q)$. The isometry group of $V$ is the semidirect product
\[ {\rm Isom} (V) = {\rm O} (V) \mbox{\Bbb n} V\, .\] 

\begin{dof} The Lie group ${\rm P} (V):= 
{\rm Isom} (V)$ is called the {\bf Poincar\'e 
group} of $V$. Its Lie algebra ${\fr p} (V) = {\fr o}(V) + V$
is called the {\bf Poincar\'e algebra} of $V$. 
\end{dof}

Next we recall some basic facts concerning the Clifford algebra $C\! \ell (V) =
C\! \ell^0 (V) + C\! \ell^1 (V)$, see \cite{L-M}. Any unit
vector $x \in V$, $\langle x,x\rangle = \pm 1$, defines an 
invertible element $x\in C\! \ell (V)$. 
The group ${\rm Pin} (V) \subset C\! \ell (V)$ generated by all unit vectors
is called the {\bf pin group}. Its subgroup ${\rm Spin}(V) :=
{\rm Pin}(V)\cap C\! \ell^0 (V)$ consisting of even elements
is called the {\bf spin group}. The {\bf adjoint  representation}
\[ {\rm Ad} : {\rm Pin} (V) \longrightarrow {\rm O}(V)\, ,\] 
\[ {\rm Ad}(x)y = xyx^{-1}\in V\, , \quad x\in {\rm Pin}(V)\, , \quad y 
\in V\, ,\]
induces a two-fold covering of the special orthogonal group: 
\[ {\rm Spin} (V) \stackrel{2:1}{\longrightarrow} {\rm SO} (V) \, .\]
In fact, if $r_x \in {\rm O}(V)$ denotes the reflection in the
hyperplane $x^{\perp} \subset V$ orthogonal to a unit vector
$x\in V$ then the following formula 
\begin{equation}\label{AdEqu} {\rm Ad} (xy) = r_x \circ r_y\end{equation} 
holds for any two unit vectors $x,y \in V$. The groups ${\rm Pin} (V)$
and ${\rm Spin} (V)$ are Lie groups whose Lie algebra is the Lie 
subalgebra ${\fr spin} (V) \subset C\! \ell^0(V)$  generated by 
the commutators $[ x,y ] = xy - yx$ of all elements
$x,y \in V$. It is canonically isomorphic to the orthogonal
Lie algebra ${\fr o} (V)$ via the adjoint representation
\begin{equation}  \label{ad1Equ} {\rm ad} = {\rm Ad}_{\ast}: 
{\fr spin}(V) \stackrel{\sim}{\longrightarrow} {\fr o}(V)\, ,\quad
{\rm ad}(x)y = [ x,y ] \, , \quad x \in {\fr spin} (V)\, ,\quad y \in V\, . 
\end{equation} 
In fact, if we identify $\wedge^2 V = {\fr o} (V)$ by 
\begin{equation}\label{wedgeEqu} (x\wedge y) (z) : = \langle y,z\rangle x - 
\langle x,z\rangle y
\, , \quad x,y,z \in V\end{equation} 
then ${\rm ad}^{-1} : \wedge^2 V = {\fr o} (V) 
\stackrel{\sim}{\longrightarrow} {\fr spin}(V)$ is given by the 
following equation:
\begin{equation} \label{adEqu} {\rm ad}^{-1}(x\wedge y) = - \frac{1}{4}[ x,y ]\, ,\quad
x,y \in V\, .\end{equation} 
In particular, ${\rm ad} (xy) = -2x\wedge y$ if $x$ and $y$  are
orthogonal. 

Any $C\! \ell (V)$-module $W$ can be decomposed into irreducible
submodules. De\-pen\-ding on the signature $(p,q)$ of $V$, there 
exist one or two irreducible $C\! \ell (V)$-modules up to equivalence.
In case there are two, they are related by the unique automorphism of
$C\! \ell (V)$ which preserves $V$ and acts on $V$ as $-{\rm Id}$. 
The restriction of an irreducible $C\! \ell (V)$-module $S$ to
$C\! \ell^0 (V)$ (respectively, to ${\rm Spin} (V)$ and ${\fr spin}(V)$) is, 
up to equivalence, independent of $W$ and is called the 
{\bf spinor module} of $C\! \ell^0 (V)$ (respectively, of ${\rm Spin} (V)$ 
and ${\fr spin}(V)$). 
The spinor module $S$ is either irreducible or $S = S^+ \oplus S^-$ is 
the sum of two irreducible {\bf semispinor modules} $S^{\pm}$, 
which may be equivalent or not, depending on the signature of $V$. 
If $S^+$ and $S^-$ are not equivalent, they are related by an 
automorphism of $C\! \ell^0 (V)$ which preserves $V$ and acts
as an isometry on $V$.  In the following we will freely use standard
notations such as $C\! \ell_{p,q} =C\! \ell (\mbox{\R}^{p,q})$, 
${\rm Spin}(p,q) = {\rm Spin}(\mbox{\R}^{p,q})$, ${\rm Spin}(p)
= {\rm Spin}(p,0)$ etc., cf.\ \cite{L-M}.

Now let $W$ be a module of the  even Clifford algebra 
$C\! \ell^0 (V)$ and $\Pi : \wedge^2 W \rightarrow V$ 
a ${\fr spin} (V)$-equivariant linear map. 
Given these data we extend the Lie bracket on  ${\fr p}(V)$ to a Lie 
bracket $[ \cdot , \cdot ]$ 
on the vector space ${\fr p} (V) + W$ by the following
requirements: 
\begin{enumerate}
\item[1)] $W$ is a ${\fr p} (V)$-submodule with trivial
action of $V$ and action of ${\fr o}(V)$ defined by 
${\rm ad}^{-1}: {\fr o} (V) \rightarrow {\fr spin}(V) \subset
C\! \ell^0 (V)$, see equation (\ref{ad1Equ}),
\item[2)] $[s,t] = \Pi (s\wedge t)$ for all $s,t \in W$. 
\end{enumerate}
The reader may observe that the Jacobi identity follows from 1) and 2). 
The resulting Lie algebra will be denoted by ${\fr p} (\Pi )$. 
Note that ${\fr p} = {\fr p} (\Pi )$ has a $\mbox{\Z}_2$-grading
\[ {\fr p} = {\fr p}_0 + {\fr p}_1\, ,\quad 
{\fr p}_0 = {\fr p} (V)\, , \quad {\fr p}_1 = W\, ,\]   
compatible with the Lie bracket, i.e.\ $[{\fr p}_a,{\fr p}_b] 
\subset {\fr p}_{a+b}$, $a,b\in \mbox{\Z}_2 = \mbox{\Z}/2\mbox{\Z}$. 
In other words, ${\fr p} (\Pi )$ is a $\mbox{\Z}_2$-graded Lie
algebra. 

\begin{dof}
\label{EPADef} Any $\mbox{\Z}_2$-graded Lie algebra ${\fr p} (\Pi ) = 
{\fr p} (V) + W$ as above is called an {\bf extended Poincar\'e 
algebra} (of {\bf signature} $(p,q)$ if $V \cong \mbox{\R}^{p,q}$). 
${\fr p} (\Pi )$ is called {\bf nondegenerate} if $\Pi$  is nondegenerate, 
i.e.\ if the map $W \ni s\mapsto \Pi (s\wedge \cdot ) \in W^{\ast} \otimes
V$ is injective. 
\end{dof}

The structure of extended Poincar\'e algebra on the vector space
${\fr p} (V) + W$ is completely determined by the ${\fr o}(V)$-equivariant
map $\Pi : \wedge^2 W \rightarrow V$ (${\fr o}(V)$ acts on $W$ via
${\rm ad}^{-1} :{\fr o} (V) \rightarrow {\fr spin} (V) 
\subset C\! \ell^0 (V)$). The set of all ${\fr o}(V)$-equivariant
linear maps $\wedge^2 W \rightarrow V$ is naturally a vector space.
In fact, it is the subspace $(\wedge^2 W^{\ast} \otimes V)^{{\fr o}(V)}$
of ${\fr o}(V)$-invariant elements of the vector space $\wedge^2 W^{\ast}
\otimes V$ of all linear maps $\wedge^2 W \rightarrow V$. 
In the classification \cite{A-C2} an explicit basis for the vector space 
$(\wedge^2 W^{\ast} \otimes V)^{{\fr o}(V)}$ of extended Poincar\'e
algebra structures on ${\fr p}(V) + W$ is constructed  
for all possible signatures $(p,q)$ of $V$ and any $C\! \ell^0 (V)$-module
$W$.  

\subsection{The canonical symmetric bilinear form $b$}
\label{CSBFSec}
Let $V = \mbox{\R}^{p,q}$ be the standard
pseudo-Euclidean vector space with scalar product $\langle \cdot , \cdot
\rangle$ of signature $(p,q)$. 
{}From now on we fix a decomposition $p = p' + p''$ and assume
that $p' \equiv 3\pmod{4}$, see Remark 2 below. 

\noindent 
{\bf Remark 1:} Notice that $p$ and $q$ are on equal footing, since 
any extended Poincar\'e 
algebra of signature $(q,p)$ is isomorphic to an extended Poincar\'e 
algebra of signature $(p,q)$. 
In fact, the canonical antiisometry which maps the standard orthonormal
basis of $\mbox{\R}^{q,p}$ to that of $\mbox{\R}^{p,q}$ induces an
isomorphism of the corresponding Poincar\'e algebras which is trivially
extended to an isomorphism of extended Poincar\'e algebras. 

We denote by $(e_i) = (e_1,\ldots ,e_{p'})$ the first $p'$ basis vectors
of the standard basis of $V$ and by $(e_i') = (e_1',\ldots , e_{p''+q}')$
the remaining ones. The two complementary orthogonal subspaces  
of  $V$ spanned by these bases are denoted by $E =  \mbox{\R}^{p'} = 
\mbox{\R}^{p',0}$ and $E' = E^{\perp} = \mbox{\R}^{p'',q}$ respectively. The vector 
spaces $V$, $E$ and $E'$ are 
oriented by their standard orthonormal bases. E.g,
the orientation of Euclidean $p'$-space $E$ defined by the basis 
$(e_i)$ is $e_1^{\ast}\wedge \cdots \wedge 
e_{p'}^{\ast} \in \wedge^{p'} E^{\ast}$. Here $(e_i^{\ast})$ denotes
the basis of $E^{\ast}$ dual to $(e_i)$. Now let 
${\fr p} (\Pi ) = {\fr p} (V) + W$ be an extended Poincar\'e algebra
of signature $(p,q)$ and $(\tilde{e}_i)$ any orhonormal
basis of $E$. Then we define a $\mbox{\R}$-bilinear form 
$b_{\Pi , (\tilde{e}_i)}$
on the $C\! \ell^0(V)$-module $W$ by: 
\begin{equation} b_{\Pi , (\tilde{e}_i)} (s,t) = \langle \tilde{e}_1, 
[\tilde{e}_2\ldots \tilde{e}_{p'}s,t] \rangle = \langle \tilde{e}_1, \Pi 
(\tilde{e}_2 \ldots \tilde{e}_{p'}s\wedge t)\rangle \, ,\quad s,t \in W \, .
\label{bEqu} 
\end{equation} 
We put $b = b(\Pi ) := b_{\Pi , (e_i)}$ for the standard basis
$(e_i)$ of $E$. 

\noindent 
{\bf Remark 2:} Equation (\ref{bEqu}) defines a skew symmetric bilinear
form on $W$ if $p'  \equiv 1\pmod{4}$. For even $p'$ the above formula does
not make sense, unless one assumes that $W$ is a   $C\! \ell (V)$-module
rather than a $C\! \ell^0(V)$-module. Here we are only interested in the case 
$p'  \equiv 3\pmod{4}$. Moreover, later on, for the construction
of homogeneous quaternionic manifolds we will put $p' =3$. 

\begin{thm} \label{bThm} The bilinear form $b$ has the following properties:
\begin{enumerate}
\item[1)] $b_{\Pi , (\tilde{e}_i)} = \pm b$ if $\tilde{e}_1
\wedge \cdots \wedge \tilde{e}_{p'} = \pm e_1 \wedge \cdots \wedge 
e_{p'}$. In particular, $b_{\Pi , (\tilde{e}_i)} =  b$  for any  
positively oriented orthonormal basis $(\tilde{e}_i)$ of $E$. 
\item[2)] $b$ is symmetric.
\item[3)] $b$ is invariant under the maximal connected subgroup $K(p',p'') 
= {\rm  Spin}(p') \cdot {\rm Spin}_0(p'',q) \subset 
{\rm Spin}(p,q)$ which preserves the orthogonal decomposition $V = E + E'$ 
(and is not ${\rm Spin}_0(p,q)$-invariant, 
unless $p''+q=0$). 
\item[4)] Under the identification 
${\fr o} (V) = \wedge^2 V = \wedge^2 E + \wedge^2 E' + E\wedge E'$, 
see equation (\ref{wedgeEqu}),  
the subspace $E\wedge E'$ acts on $W$ by $b$-symmetric endomorphisms
and the subalgebra $\wedge^2 E \oplus \wedge^2 E' \cong {\fr o}(p') 
\oplus {\fr o}(p'',q)$ acts on $W$ by $b$-skew symmetric endomorphisms. 
\end{enumerate}
\end{thm}

\noindent
{\bf Proof:} Obviously 4) $\Rightarrow$ 3). We show first that 
3) $\Rightarrow$ 1). If $(\tilde{e}_i)$ is a positively
oriented orthonormal basis then there exists $\varphi \in {\rm Spin}(p')$
such that ${\rm Ad}(\varphi ) e_i = \tilde{e}_i$, 
$i = 1,\ldots,p'$. Now $\Pi = [\cdot ,\cdot ]| \wedge^2W : 
\wedge^2W \rightarrow V$ is ${\fr spin}(V)$-equivariant
and hence equivariant under the connected group ${\rm Spin}_0(V)
\subset {\rm Spin}(V)$. In particular, $\Pi$ is 
${\rm Spin}(p')$-equivariant. Under the condition 3), this implies
that
\begin{eqnarray*} b (s,t) &=& b(\varphi^{-1}s,\varphi^{-1}t) =\langle 
e_1,[e_2\ldots e_{p'}\varphi^{-1}s,\varphi^{-1}t]\rangle \\
&=& \langle e_1,[\varphi^{-1}{\rm Ad}_{\varphi}e_2 \ldots {\rm Ad}_{\varphi}
e_{p'}s, \varphi^{-1}t]\rangle \\
&=&  \langle e_1,{\rm Ad} ({\varphi}^{-1}) 
[\tilde{e}_2\ldots \tilde{e}_{p'}s,t]
\rangle = \langle \tilde{e}_1, 
[\tilde{e}_2\ldots \tilde{e}_{p'}s,t]\rangle \\
&=& 
b_{\Pi , (\tilde{e}_i)} (s,t)\,  ,\quad s,t \in W\, .
\end{eqnarray*} 
Here we have used the notation ${\rm Ad}_{\varphi} = 
{\rm Ad}(\varphi )$. 
The case of negatively oriented orthonormal basis $(\tilde{e}_i)$
follows now from the Clifford relation $\tilde{e}_i\tilde{e}_j
 = - \tilde{e}_j\tilde{e}_i$, $i \neq j$. 

Next we prove 2) using first the ${\fr spin}(V)$-equivariance 
of $\Pi$ then equation (\ref{adEqu}) and eventually $p' \equiv 3\pmod{4}$:
\begin{eqnarray*} 
b(t,s) &=& \langle e_1,[e_2 \ldots e_{p'}t,s]\rangle\\
&=& -\langle e_1 ,[e_4\ldots e_{p'}t,e_2e_3s]\rangle + \langle e_1 , 
{\rm ad} (e_2e_3)
[e_4 \ldots e_{p'}t,s]\rangle\\
&=& -\langle e_1 ,[e_4 \ldots e_{p'}t,e_2e_3s]\rangle = \cdots =
-\langle e_1 ,[t,e_2 \ldots e_{p'}s] \rangle \\
&=&  \langle e_1 ,[e_2 \ldots e_{p'}s,t] \rangle  = b(s,t)\, . 
\end{eqnarray*}
Finally we prove 4). By equation (\ref{adEqu}) we have to check that
$e_ie_j, e_k'e_l'\in {\fr spin}(V)$ ($i \neq
j$ and $k\neq l$) act by $b$-skew symmetric endomorphisms and $e_i e_k'$ by
a $b$-symmetric endomorphism on $W$.   This is done in 
the next computation, in which we use again the equivariance of $\Pi$
and equations (\ref{AdEqu}) and (\ref{adEqu}) to express the 
adjoint representations $\rm Ad$ and $\rm ad$ respectively:  
\begin{eqnarray*}
b(e_1e_2s,t) &=& \langle e_1,[e_2\ldots e_{p'}e_1e_2s,t]\rangle =
\langle e_1,[e_1e_2e_2 \ldots e_{p'}s, e_1e_2e_1e_2t]\rangle \\
 &=& \langle e_1,{\rm Ad}(e_1e_2)[e_2\ldots e_{p'}s,e_1e_2t]\rangle =
-\langle e_1,[e_2\ldots e_{p'}s,e_1e_2t]\rangle\\
&=& -b(s,e_1e_2t)\, ,\\
b(e_2e_3s,t) &=& \langle e_1,[e_2\ldots e_{p'}e_2e_3s,t]\rangle = 
\langle e_1,[e_2e_3e_2 \ldots e_{p'}s,t]
\rangle \\
&=& \langle e_1,{\rm ad}(e_2e_3)[e_2 \ldots e_{p'}s,t]\rangle  
-\langle e_1,[e_2\ldots e_{p'}s,e_2e_3t]\rangle \\
&=& -b(s,e_2e_3t)\, ,\\
b(e_k'e_l's,t) &=& \langle e_1,[e_2\ldots e_{p'}e_k'e_l's,t]\rangle = 
\langle e_1,[e_k'e_l'e_2\ldots e_{p'}s,t]\rangle \\
&=& -\langle e_1,[e_2\ldots e_{p'}s,e_k'e_l't]\rangle + 
\langle e_1,{\rm ad}(e_k'e_l')[e_2\ldots e_{p'}s,t]\rangle \\
&=& -b(s,e_k'e_l't) + 0 = -b(s,e_k'e_l't)\, .
\end{eqnarray*}
This already proves 3) and hence 1); in particular we have: 
\[ b_{\Pi,(e_1,e_2,\ldots , e_{p'})} =  b_{\Pi,(e_2,\ldots , e_{p'},e_1)}
\, .\]
Due to this symmetry, it is sufficient to check that $e_2e_k'$ acts as 
$b$-symmetric endomorphism on $W$:
\begin{eqnarray*}
b(e_2e_k's,t) &=& \langle e_1,[e_2\ldots e_{p'}e_2e_k's,t]\rangle =
\langle e_1,[e_3\ldots e_{p'}e_k's,t]\rangle \\
&\stackrel{(\ast )}{=}& -\langle e_1,[s,e_3\ldots e_{p'}e_k't]\rangle = 
-\langle e_1,[s,e_2\ldots e_{p'}e_2e_k't]\rangle \\
&=& \langle e_1,[e_2\ldots e_{p'}e_2e_k't,s]\rangle = b(e_2e_k't,s)\\
&=& b(s,e_2e_k't)\, .
\end{eqnarray*}
At ($\ast$) we have used $(p'-1)/2$ times the ${\fr spin}(V)$-equivariance
of $\Pi$ and the fact that $(p'-1)/2$ is odd if $p' \equiv 3\pmod{4}$. 
$\Box$

\begin{dof} The bilinear form $b = b(\Pi ) = b_{\Pi ,(e_1,\ldots , e_{p'})}$ 
defined above is called
the {\bf canonical symmetric bilinear form} on $W$ associated 
to the ${\fr spin}(V)$-equivariant map $\Pi :\wedge^2W 
\rightarrow V = \mbox{\R}^{p,q}$ and the decomposition $p = p' + p''$.  
\end{dof}

\begin{prop} \label{kerProp} The kernels of the linear 
maps $\Pi : W \rightarrow 
W^{\ast}\otimes V$ and $b = b(\Pi ) : W \rightarrow W^{\ast}$ coincide:
${\rm ker}\, \Pi = {\rm ker}\, b$. 
\end{prop}

\noindent {\bf Proof:} It follows from Thm.~\ref{bThm}, 4)
that $W_0 := {\rm ker}\, b \subset W$ is ${\fr o}(V)$-invariant. This
implies that $\Pi (W_0 \wedge W) \subset V$ is an ${\fr o}(V)$-submodule.
The definition of $W_0$ implies that $\Pi (W_0 \wedge W) \subset E' =
E^{\perp}$ and hence by Schur's lemma  $\Pi (W_0 \wedge W) = 0$.
This proves that ${\rm ker}\, b \subset {\rm ker}\, \Pi$. On the other hand, 
we have the obvious inclusion ${\rm ker}\, \Pi \subset 
{\rm ker}\, (b\circ e_2 \ldots e_{p'}) = {\rm ker}\, b$. Here the 
equation follows
from the ${\rm Spin}(p')$-invariance of $b$, see Thm.~\ref{bThm}, 3). $\Box$

\begin{cor} \label{bpiCor} ${\fr p} (\Pi )$ is nondegenerate (see Def.\ 
\ref{EPADef}) if and only if 
$b (\Pi )$ is nondegenerate.
\end{cor} 

\begin{thm}\label{decompThm} Let ${\fr p} (\Pi ) = {\fr p} (V) + W$ 
be any extended 
Poincar\'e algebra of signature $(p,q)$, $p = p' + p''$, 
$p' \equiv 3\pmod{4}$, and $b$ the canonical symmetric bilinear
form associated to these data.  
Then there exists a $b$-orthogonal decomposition $W = \oplus_{i=0}^{l+m}W_i$ 
into $C\! \ell^0 (V)$-submodules with the following properties
\begin{enumerate}
\item[1)] $[W_0,W] = 0$, $[W_i,W_j] = 0$ if $i\neq j$ and 
$[W_i,W_i] = V$ for all $i = 1,2,\ldots , l$. 
\item[2)] $W_0 = {\rm ker}\, b$ and $W_i$ is $b$-nondegenerate for all
$i \ge 1$.  
\item[3)] For $i = 1,\ldots ,l$ the $C\! \ell^0 (V)$-submodule 
$W_i$ is irreducible and for $j = l+1, \ldots ,l+m$ the 
$C\! \ell^0 (V)$-submodule 
$W_j = X_j \oplus X_j'$ is the direct sum of two irreducible 
$b$-isotropic $C\! \ell^0 (V)$-submodules.
\item[4)] The restriction of $b$ to a bilinear form on any irreducible
$C\! \ell^0 (V)$-submodule of $X = \oplus_{j=l+1}^{l+m}W_j$ vanishes. 
\end{enumerate}
\end{thm}

\noindent {\bf Proof:}  By Prop.\ \ref{kerProp}, $W_0 := {\rm ker}\, b =
{\rm ker}\, \Pi$ satisfies $[W_0,W] =0$. As kernel of the 
${\fr o}(V)$-equivariant map $\Pi$ the subspace $W_0$ is 
${\fr o}(V)$-invariant and 
hence a $C\! \ell^0 (V)$-submodule. We denote by $W'$ a complementary
$C\! \ell^0 (V)$-submodule. Every such submodule is $b$-nondegenerate.
Let $W_i \subset W'$ be any irreducible
$C\! \ell^0 (V)$-submodule. By Thm.~\ref{bThm}, 4), 
${\rm ker}(b|W_i \times W_i$) is ${\fr o}(V)$-invariant and
hence a $C\! \ell^0 (V)$-submodule. Now by Schur's lemma we conclude
that either $b|W_i \times W_i = 0$ or $b$ is nondegenerate on $W_i$.
In particular, we can decompose $W' = \oplus_{i=1}^{l} W_i \oplus
X$ as direct $b$-orthogonal sum of $b$-nondegenerate $C\! \ell^0
(V)$-submodules such that $W_i$ is irreducible and the restriction
of $b$ to a bilinear form on any irreducible $C\! \ell^0 (V)$-submodule
of $X$ vanishes. Let $Y, Z \subset X$ be two such submodules, $Y\neq Z$.
The bilinear form $b$ induces a linear map $Y \rightarrow Z^{\ast}$. 
By Thm.~\ref{bThm}, 4), the kernel of this map is ${\fr o}(V)$-invariant and
hence a $C\! \ell^0 (V)$-submodule. Now Schur's lemma implies that 
either the kernel is $Y$ and hence the restriction of $b$ to a bilinear
form on $Y \oplus Z$ vanishes or the kernel is trivial and $b$ is 
nondegenerate on $Y \oplus Z$. In the second case $X$ splits as direct 
$b$-orthogonal sum: $X = (Y\oplus Z) \oplus (Y \oplus Z)^{\perp}$.
This shows that $X = \oplus_{i = l+1}^{l+m}W_i$ is the direct orthogonal
sum of $b$-nondegenerate $C\! \ell^0 (V)$-submodules $W_i$ such that 
$W_i = X_i \oplus X_i'$ is the direct sum of two $b$-isotropic 
irreducible submodules. 
This proves 2), 3) and 4). Now  1) is established 
applying Schur's lemma to the  
${\fr o} (V)$-equivariant map $\Pi$. In fact, $b(W_i,W_j) = 0$ 
(respectively, $b(W_i,W_i) \neq 0$) implies $\Pi (W_i \wedge W_j) 
\subset E'$ (respectively, $\Pi (\wedge^2 W_i) \neq 0$)
and thus $[W_i,W_j] = \Pi (W_i \wedge W_j) = 0$ 
(respectively, $[W_i,W_j] = \Pi (\wedge^2 W_i) = V$). $\Box$ 

Next we will construct the subgroup $\hat{K}(p',p'') \subset {\rm Spin}(p,q)$
which consists of all elements preserving the orthogonal decomposition
$V = E + E'$ and the canonical symmetric bilinear form $b$ on $W$.
Its identity component is the group $K(p',p'') \subset {\rm Spin}_0(p,q)$ 
introduced above. We will see that if $p''= 0$ then 
$\hat{K}(p',p'') = \hat{K}(p,0)$ is a maximal compact subgroup
of ${\rm Spin}(p,q)$, which together with the element 
$1 \in C\! \ell^0_{p,q}$ generates the even Clifford algebra 
$C\! \ell^0_{p,q}$. This property will be very useful in the next section. 

We denote by $x \mapsto x'$ the linear map $\mbox{\R}^q = 
\mbox{\R}^{q,0} \rightarrow \mbox{\R}^{0,q}$ which maps the standard
orthonormal basis $(e_1,\ldots ,e_q)$ of $\mbox{\R}^q$ to the standard
orthonormal basis $(e_1',\ldots ,e_q')$ of $\mbox{\R}^{0,q}$. It is
an antiisometry: $\langle x',x' \rangle = - \langle x,x \rangle$. Let
$\omega_{p'} = e_1 \ldots e_{p'}$ be the volume element of $C\! \ell_{p'} = 
C\! \ell_{p',0}$, $(e_1,\ldots ,e_{p'})$ the standard orthonormal
basis of $\mbox{\R}^{p'} = \mbox{\R}^{p',0} \subset \mbox{\R}^{p,q}$. 
Note that, since $p'$ is odd, the volume element $\omega_{p'}$ commutes
with $\mbox{\R}^{p'}$ and anticommutes with $\mbox{\R}^{p'',q} \supset 
\mbox{\R}^{0,q}$. Moreover, it 
satisfies $\omega_{p'}^2 = 1$, due to $p'\equiv 3\pmod{4}$. 

\begin{lemma} \label{iotaLemma} The map 
\[ \mbox{\R}^q \ni x \mapsto \omega_{p'} x' \in C\! \ell_{p,q}^0\] 
extends to an embedding $\iota : C\! \ell_q \hookrightarrow C\! \ell_{p,q}^0$
of algebras, which restricts to an embedding 
$\iota |{\rm Pin}(q) : {\rm Pin}(q) \hookrightarrow 
{\rm Spin}(p,q)$ of groups. 
\end{lemma}

\noindent {\bf Proof:} It follows from $(\omega_{p'}x')^2 = -
\omega_{p'}^2{x'}^2 = -{x'}^2 = \langle x',x' \rangle = -\langle x,x
\rangle$ that the map $x \mapsto \omega_{p'} x'$ extends to a homomorphism
$\iota$ of Clifford algebras. Recall that $\C\! \ell_q$ is either
simple or the sum of two simple ideals. In the first case, 
we can immediately conclude that ${\rm ker}\, \iota$ is trivial
and hence $\iota$ an embedding. In the second case, the  
two simple ideals of $\C\! \ell_q$ are $\C\! \ell_q^{\pm} := (1\pm \omega_q)
\C\! \ell_q$, where $\omega_q =e_1\ldots e_q$ is the volume element
of $\C\! \ell_q$. Now it is sufficient to check that $\iota (1\pm \omega_q)
\neq 0$. Using the fact that $q$ is odd if $C\! \ell_q$ is not simple we 
compute 
\[ \iota (\omega_q) = \pm \omega_{p'}^qe_1' e_2'\ldots e_q' = 
\pm e_1e_2\ldots e_{p'}e_1'e_2' \ldots e_q'\, .\]
This shows that $\iota (1\pm \omega_q ) = 1 \pm \iota (\omega_q) \neq 0$. 
$\Box$ 

We denote by $\hat{K}(p',p'') \subset {\rm Spin}(p,q)$ the subgroup
generated by the subgroups ${\rm Spin}(p')\cdot {\rm Spin}_0(p'',q)
\subset {\rm Spin}(p,q)$ and $\iota ({\rm Pin}(q)) \subset {\rm Spin}(p,q)$.
\begin{thm}\label{KHATThm} 
The group $\hat{K}(p',p'')$ has the following properties:
\begin{enumerate}
\item[1)] $\hat{K}(p',p'') \subset {\rm Spin}(p,q)$ consists of all elements
preserving the orthogonal decomposition
$V = E + E'$ and the canonical symmetric bilinear form $b$ on $W$. 
Its identity component is the group $K(p',p'') =
{\rm Spin}(p')\cdot {\rm Spin}_0(p'',q) = \hat{K}(p',p'')\cap 
{\rm Spin}_0(p,q)$. 
\item[2)] The homogeneous space ${\rm Spin}(p,q)/\hat{K}(p',p'')$ is connected.
\item[3)]  $\hat{K}(p',p'')$ is compact if and only if $q=0$ or $p'' =0$
(and hence $p=p'$). In the latter case $\hat{K}(p',p'') = \hat{K}(p,0) =
{\rm Spin}(p)\cdot \iota ( {\rm Pin}(q)) \subset 
{\rm Spin}(p,q)$ is a maximal compact subgroup and 
$\hat{K}(p,0) \cong  ({\rm Spin}(p) \times {\rm Pin}(q))/\{ \pm 1\}$. Finally,
in this case, the even Clifford algebra $C\! \ell^0_{p,q}$ is 
generated by $1$ and 
$\hat{K}(p,0)$. 
\end{enumerate}
\end{thm}

\noindent {\bf Proof:} The first part of 1) can be checked using 
Thm.\ \ref{bThm}, 4) and implies the second part of 1). To prove 2)
it is sufficient to observe that $\iota ({\rm Pin}(q)) \cong {\rm Pin}(q)$ 
has nontrivial intersection with all connected components of 
${\rm Spin}(p,q)$ (due to our assumption $p\ge 3$ there are 
two such components if $q\neq 0$). The first part of 3) now 
follows simply from the fact that ${\rm Spin}_0(p'',q)$ is
compact if and only if $p'' =0$ or $q=0$. The compact group
$\hat{K}(p,0) \subset {\rm Spin}(p,q)$ is maximal compact, 
because it has the same number of connected components as ${\rm Spin}(p,q)$, 
and from ${\rm Spin}(p) \cap \iota ({\rm Pin}(q)) = \{ \pm 1\}$ we 
obtain the isomorphism $\hat{K}(p,0) = {\rm Spin}(p) \cdot \iota ({\rm Pin}(q))
\cong ({\rm Spin}(p) \times {\rm Pin}(q))/\{ \pm 1\}$.  Finally,
to prove the last statement, one easily checks that $\hat{K}(p,0)$ contains all
quadratic monomials $xy$ in unit vectors $x,y \in \mbox{\R}^{p,0}
\cup \mbox{\R}^{0,q}$. $\Box$

\begin{cor} \label{injCor} The correspondence $\Pi \mapsto b(\Pi )$ 
defines an injective
linear map $(\wedge^2W^{\ast} \otimes V)^{{\fr o}(V)} \rightarrow 
(\vee^2W^{\ast})^{\hat{K}(p',p'')}$. 
\end{cor} 

\noindent {\bf Proof:} The existence of the map follows from   
Thm.~\ref{KHATThm}. We prove the injectivity. From $b(\Pi ) = 0$  
it follows that $\Pi ( \wedge^2 W) \subset E'$ and 
hence by Schur's lemma $\Pi = 0$. $\Box$

\subsection{The set of isomorphism classes of extended Poincar\'e
algebras} \label{isomSec} 
Starting from the decomposition proven in Thm.~\ref{decompThm} we will 
derive the classification of extended Poincar\'e
algebras of signature $(p,q)$, $p \equiv 3\pmod{4}$, up to isomorphism. 
It will turn out that the space of isomorphism classes is naturally
parametrized by a finite number of integers. 
We fix the decomposition $p = p' + p''$, $p' = p$, $p'' = 0$, and 
for any extended Poincar\'e algebra ${\fr p}(\Pi )$ of signature
$(p,q)$ as above we consider the canonical symmetric bilinear form $b = 
b_{\Pi ,(e_1, \ldots e_p)}$. 

\begin{thm}\label{decompIIThm} Let ${\fr p} (\Pi ) = {\fr p} (V) + W$ 
be any extended 
Poincar\'e algebra of signature $(p,q)$, $p \equiv 3\pmod{4}$. 
Then there exists a $b$-orthogonal decomposition $W = \oplus_{i=0}^lW_i$ 
into $C\! \ell^0 (V)$-submodules with the following properties
\begin{enumerate}
\item[1)] $[W_0,W] = 0$, $[W_i,W_j] = 0$ if $i\neq j$ and 
$[W_i,W_i] = V$ for all $i = 1,2,\ldots , l$. 
\item[2)] $W_0 = {\rm ker}\, b$ and $W_i$ is $b$-nondegenerate for all
$i \ge 1$.  
\item[3)] $W_i$, $i\ge 1$, is an irreducible
$C\! \ell^0 (V)$-submodule on which $b$ is (positive or negative)
definite. 

\end{enumerate}
\end{thm}

\noindent {\bf Proof:} Let $W = \oplus_{i=1}^{l+m}W_i$ be a decomposition
as in Thm.~\ref{decompThm}. It only remains to prove
that $b$ is definite on $W_i$ for $i =1, \ldots ,l$ and that
$X = \oplus_{i=l+1}^{l+m} W_i = 0$. This follows from 
Lemma~\ref{defLemma} and Lemma~\ref{ineqLemma} below. $\Box$

\begin{lemma} \label{KirredLemma} The restriction of an 
irreducible $C\! \ell_{p,q}^0$-module
$\Sigma$ to a module of the ma\-xi\-mal compact subgroup 
$\hat{K} = \hat{K}(p,0)   = {\rm Spin}(p) \cdot
\iota ({\rm Pin}(q)) \subset {\rm Spin}(p,q)$ is irreducible. Here 
$\iota :{\rm Pin}(p) \hookrightarrow {\rm Spin}(p,q)$ is the embedding of
Lemma \ref{iotaLemma}. Moreover, $\Sigma$ is irreducible as module
of the connected group $K = K(p,0) = {\rm Spin}(p)\cdot \iota ({\rm Spin}(q)) =
{\rm Spin}(p)\cdot {\rm Spin}(q)$ if and only if $n = p+q \equiv 2,4,5$
or $6\pmod{8}$. If $n \equiv 0,1,3$ or $7\pmod{8}$ then $\Sigma$ is the
sum of two irreducible $K$-submodules. 
\end{lemma} 

\noindent {\bf Proof:} Recall that $C\! \ell_{p,q} = 
C\! \ell_p \hat{\otimes} C\! \ell_{0,q}$ is (identified with)
the $\mbox{\Z}_2$-graded tensor product of the Clifford algebras
$C\! \ell_p = C\! \ell_{p,0}$ and $C\! \ell_{0,q}$. It is easily checked, 
using the classification of Clifford algebras and their modules,
see \cite{L-M}, that any irreducible $C\! \ell_{p,q}^0$-module $\Sigma$
is irreducible as module of the subalgebra $C\! \ell^0_p \otimes
C\! \ell^0_{0,q} \subset C\! \ell^0_{p,q} = C\! \ell^0_p \otimes
C\! \ell^0_{0,q} + C\! \ell^1_p \otimes
C\! \ell^1_{0,q}$ if $n \equiv 2,4,5$ or $6\pmod{8}$ and is the sum of
two irreducible submodules if $n \equiv 0,1,3$ or $7\pmod{8}$. 
Now Lemma \ref{KirredLemma} follows from the fact that 
$C\! \ell^0_p \otimes C\! \ell^0_{0,q}$ (respectively, $C\! \ell^0_{p,q}$) 
is the subalgebra of $C\! \ell_{p,q}$ generated by $1$ and $K = {\rm Spin}(p)
\cdot {\rm Spin}(0,q)$ (respectively, by $1$ and $\hat{K} = {\rm Spin}(p)
\cdot \iota ({\rm Pin}(q))$). $\Box$

\begin{lemma} \label{defLemma} A $C\! \ell_{p,q}^0$-module 
$W$ is irreducible if and only if
it is irreducible as module of the maximal compact subgroup $\hat{K} \subset 
{\rm Spin}(p,q)$. In this case $(\vee^2W^{\ast})^{\hat{K}}$ is one-dimensional
and is spanned by a positive definite scalar product on $W$. 
Let $W$ be an irreducible  $C\! \ell_{p,q}^0$-module and  
$\Pi : \wedge^2W \rightarrow V$ be any ${\fr o}(V)$-equivariant 
linear map. Then either $\Pi = 0$ or $b (\Pi )$ is a definite 
$\hat{K}$-invariant symmetric bilinear form. 
\end{lemma} 

\noindent {\bf Proof:} The first statement follows from 
Lemma~\ref{KirredLemma}. Since $\hat{K}$ is compact there exists
a positive definite $\hat{K}$-invariant symmetric bilinear form
on $W$. From the irreducibility of $W$ we conclude by Schur's lemma that 
$(\vee^2 W^{\ast})^{\hat{K}}$ is spanned by this form. 
Now the last statement is an immediate consequence of Cor.~\ref{injCor}. 
 $\Box$ 

\begin{lemma} \label{ineqLemma} Let $W$ be a $C\! \ell^0(V)$-module and 
$\Pi : \wedge^2W \rightarrow V$ an ${\fr o}(V)$-equivariant linear map
such that $b = b(\Pi )$ is nondegenerate. Suppose that $W = \Sigma 
\oplus \Sigma'$ is the direct sum of two irreducible submodules
$\Sigma$ and $\Sigma'$.
Then there exists a $b$-orthogonal decomposition $W = \Sigma_1
\oplus \Sigma_2$ into two $b$-nondegenerate (and hence $b$-definite
by Lemma \ref{defLemma}) irreducible
submodules $\Sigma_1$ and $\Sigma_2$. 
\end{lemma}

\noindent {\bf Proof:} It is sufficient to show that $W$ contains
a $b$-nondegenerate irreducible $C\! \ell^0(V)$-submodule
$\Sigma_1$. Then $\Sigma_2 := \Sigma_1^{\perp}$ is a $b$-nondegenerate 
$\hat{K}$-submodule. It is also a $C\! \ell^0(V)$-submodule, because the 
algebra $C\! \ell^0(V)$ is generated by $1$ and $\hat{K}$, and it 
is irreducible 
since the $C\! \ell^0(V)$-module $W$ is the direct sum of only two 
irreducible submodules. If $W$ does not contain any $b$-nondegenerate
irreducible $C\! \ell^0(V)$-submodule
$\Sigma_1$ then, by Schur's lemma, the restriction of $b$ to a bilinear
form on any irreducible $C\! \ell^0(V)$-submodule vanishes. 
In the following we derive a contradiction from this assumption. 
Since the bilinear form $b$ is nondegenerate it defines a nondegenerate 
pairing between the $b$-isotropic subspaces $\Sigma$ and $\Sigma'$.
Due to the $\hat{K}$-invariance of $b$ (Thm.~\ref{KHATThm}) 
$b: \Sigma' \stackrel{\sim}{\rightarrow} \Sigma^{\ast}$ is a
$\hat{K}$-equivariant isomorphism. On the other hand $\Sigma^{\ast} \cong 
\Sigma$ as irreducible modules of the compact group $\hat{K}$. 
This shows that $\Sigma$ and $\Sigma'$ are equivalent as $\hat{K}$-modules and
thus as $C\! \ell^0(V)$-modules, because the 
algebra $C\! \ell^0(V)$ is generated by $1$ and $\hat{K}$. 
Hence there exists a $C\! \ell^0(V)$-equivariant isomorphism
$\varphi :\Sigma \stackrel{\sim}{\rightarrow} \Sigma'$. We 
define two $\hat{K}$-invariant bilinear forms $\beta_{\pm}$ on 
$\Sigma$ by:
\[ \beta_{\pm} (s,t) := b(\varphi (s),t) \pm b(\varphi (t),s)
\, ,\quad s,t\in \Sigma \, .\]
$\beta_+$ is symmetric and $\beta_-$ is skew symmetric. If 
$\beta_+ \neq 0$ then it is a definite $\hat{K}$-invariant scalar
product on $\Sigma$, since $\Sigma$ is an irreducible module
of the compact group $\hat{K}$. So for $s\in \Sigma -\{ 0\}$ we 
obtain 
\[ 0 \neq \beta_+(s,s) = b(\varphi (s),s) + b(s,\varphi (s)) =
b(s+\varphi (s),s+\varphi (s))\, .\] 
This implies that the irreducible $C\! \ell^0(V)$-submodule 
$\Sigma_{\varphi} := \{ s + \varphi (s)| s\in \Sigma \} \subset W$
is $b$-definite, which contradicts our assumption. We conclude
that $\beta_+ = 0$ and hence $\beta_- = 2b\circ \varphi$ is a 
$\hat{K}$-invariant 
symplectic
form on $\Sigma$. Let $\beta$ be a $\hat{K}$-invariant positive definite
scalar product on $\Sigma$ (such scalar products  exist since
$\hat{K}$ is compact). We define a $\hat{K}$-equivariant isomorphism
$\chi :\Sigma \rightarrow \Sigma$ by the equation
\[ \beta (s,t) = \beta_- (\chi (s),t)\, ,\quad s,t\in \Sigma \,.\]
Then $\psi := \varphi \circ \chi : \Sigma \rightarrow \Sigma'$ is 
a $\hat{K}$-equivariant and hence $C\! \ell^0(V)$-equivariant isomorphism.
Using $\beta_+ = 0$, we compute for $s \in \Sigma -\{ 0\}$:
\begin{eqnarray*} b(s + \psi (s), s + \psi (s)) &=& 
b(\psi (s),s) + b(s,\psi (s)) = b(\psi (s),s) - b(\varphi (s), \chi (s))\\
&=& \beta_-(\chi (s),s) = \beta (s,s) \neq 0\, .
\end{eqnarray*} 
As above, this implies that $W$ contains a $b$-definite irreducible
$C\! \ell^0(V)$-submodule $\Sigma_{\psi} \cong \Sigma$ contradicting
our assumption. 
$\Box$ 

\begin{thm}\cite{A-C2} \label{dimThm}
Let $V = \mbox{\R}^{p,q}$, $p \equiv 3\pmod{4}$. If $W$ is an irreducible
$C\! \ell^0(V)$-module then 
\[ \dim (\wedge^2 W^{\ast} \otimes V)^{{\fr o}(V)} = 1\, .\] 
\end{thm}

\noindent {\bf Proof:} Cor.~\ref{injCor} and Lemma.~\ref{defLemma} show that
$\dim (\wedge^2 W^{\ast} \otimes V)^{{\fr o}(V)} \le 1$. On the other
hand, there exists a nontrivial ${\fr o}(V)$-equivariant
linear map $\wedge^2W \rightarrow V$, see \cite{A-C2}, and hence
$\dim (\wedge^2 W^{\ast} \otimes V)^{{\fr o}(V)} \ge 1$.  This proves
that $\dim (\wedge^2 W^{\ast} \otimes V)^{{\fr o}(V)} = 1$.  $\Box$ 

Next, in order to parametrize the isomorphism classes of extended Poincar\'e 
algebras ${\fr p} (\Pi ) = {\fr p}(V) + W$ we associate a certain number 
of nonnegative integers to ${\fr p} (\Pi )$. Let us first consider
the case when there is only one irreducible $C\! \ell^0(V)$-module
$\Sigma$ up to equivalence. Then $W$ is necessarily isotypical
and, due to Thm.~\ref{decompIIThm}, there exists a $b$-orthogonal
decomposition $W = W_0 \oplus \oplus_{i=1}^l W_i$ with $W_0  = 
{\rm ker}\, b = {\rm ker}\, \Pi \cong l_0 \Sigma$ and irreducible
$C\! \ell^0(V)$-submodules $W_i \cong \Sigma$ for $i = 1, \ldots ,l$
on which $b$ is definite. We denote by $l_+$ (repectively $l_-$) the
number of summands $W_i$ on which $b$ is positive (respectively, 
negative) definite. Note that the triple $(l_0,l_+,l_-)$ does not depend
on the choice of decomposition. 

\begin{thm} \label{moduliThmI} Let $p \equiv 3\pmod{4}$ and 
$q\not\equiv 3\pmod{4}$. Then 
the isomorphism class $[ {\fr p} (\Pi )]$ of an extended Poincar\'e
algebra  ${\fr p} (\Pi )$ of signature $(p,q)$ is completely
determined by the triple $(l_0,l_+,l_-)$ introduced above. 
We put ${\fr p}(p,q,l_0,l_+,l_-) := [ {\fr p} (\Pi )]$. Then 
${\fr p}(p,q,l_0,l_+,l_-) = {\fr p}(p,q,l_0',l_+',l_-')$ if and
only if $l_0 = l_0'$ and $\{ l_+,l_-\} = \{ l_+',l_-'\}$. 
\end{thm}

\noindent {\bf Proof:} 
If $p \equiv 3\pmod{4}$ then there is only one irreducible
$C\! \ell^0_{p,q}$-module $\Sigma$ up to equivalence if and only if
$q \not\equiv 3\pmod{4}$. Let ${\fr p}(\Pi ) = {\fr p}(V) + W$ and 
${\fr p}(\Pi' ) =  {\fr p}(V) + W'$ be two extended Poincar\'e
algebras of signature $(p,q)$ with the same integers $l_0 = l_0(\Pi ) 
= l_0(\Pi' )$, $l_+ = l_+(\Pi ) =  l_+(\Pi' )$ and 
$l_- = l_-(\Pi ) =  l_-(\Pi' )$. Then the modules $W$ and $W'$ are
equivalent and we can assume that $W = W' = W_0 \oplus \oplus_{i=1}^{l=l_++l_-}
W_i$ is a decomposition as above. In particular, it is $b(\Pi )$- and 
$b(\Pi')$-orthogonal, $b(\Pi )$ and $b(\Pi')$ are both positive definite or 
both negative definite on $W_i$ for $i\ge 1$, $\Pi (W_0\wedge W) =
\Pi' (W_0\wedge W) = 0$, $\Pi (W_i\wedge W_j) = \Pi' (W_i\wedge W_j) = 0$ if 
$i \neq j$ and $\Pi (\wedge^2W_i) = \Pi' (\wedge^2W_i) = V$ if $i\ge 1$. 
So the maps $\Pi$ and $\Pi'$ are completely determined by their restrictions
$\Pi_i := \Pi|\wedge^2W_i \neq 0$ and $\Pi_i' := \Pi'|\wedge^2W_i
\neq 0$ ($i\ge 1$) respectively. By Thm.~\ref{dimThm} $\Pi_i' = \lambda_i\Pi_i$
($i\ge 1$) for some constant $\lambda_i \in \mbox{\R}^{\ast}$. Now
$b = b(\Pi )$ and $b' = b(\Pi')$ are both positive definite or both negative
definite on $W_i$ and hence $\lambda_i = \mu_i^2 > 0$. Now we can
define an isomorphism $\varphi :{\fr p}(\Pi ) \rightarrow
{\fr p}(\Pi')$ by $\varphi|{\fr p}(V)+W_0 = {\rm Id}$ and 
$\varphi|W_i = \mu_i {\rm Id}$. This shows that the integers 
$(l_0,l_+,l_-)$ determine the extended Poincar\'e algebra
${\fr p}(\Pi )$ of signature $(p,q)$ up to isomorphism. The  
$\mbox{\Z}_2$-graded Lie algebras ${\fr p}(\Pi )$ and 
${\fr p}(-\Pi )$ are isomorphic via $\alpha : 
{\fr p}(\Pi ) \rightarrow {\fr p}(-\Pi )$ defined by: 
$\alpha|{\fr o}(V)+W = {\rm Id}$ and $\alpha|V = -{\rm Id}$. 
This proves ${\fr p}(p,q,l_0,l_+,l_-)  = {\fr p}(p,q,l_0,l_-,l_+)$.
It remains to show that ${\fr p}(p,q,l_0,l_+,l_-) = 
{\fr p}(p,q,l_0',l_+',l_-')$ implies $l_0 = l_0'$ and 
$\{ l_+,l_-\} = \{ l_+',l_-'\}$. Let ${\fr p}(\Pi )
= {\fr p}(V) + W \in {\fr p} (p,q,l_0,l_+,l_-)$ and 
${\fr p}(\Pi' ) = {\fr p}(V) + W' \in {\fr p} (p,q,l_0',l_+',l_-')$ 
be representative extended Poincar\'e algebras and $W = W_0 \oplus
\oplus_{i=1}^l W_i$ ($l = l_+ + l_-$) a decomposition as above. 
We assume that there exists an isomorphism $\varphi : 
{\fr p} (\Pi ) \rightarrow {\fr p} (\Pi')$ of $\mbox{\Z}_2$-graded Lie 
algebras, i.e.\ $\varphi {\fr p}(V) =  {\fr p}(V)$ and 
$\varphi W = W$. The automorphism $\varphi|{\fr p}(V)$ preserves the 
radical $V$ and maps the Levi subalgebra ${\fr o}(V)$ to an 
other Levi subalgebra of ${\fr p}(V)$. Now by Malcev's theorem
any two Levi subalgebras are conjugated by an inner automorphism,
see \cite{O-V}.  So, using an inner automorphism of ${\fr p}(\Pi )$,
we can assume that $\varphi {\fr o}(V) = {\fr o}(V)$. 
The subalgebra $\varphi ({\fr o}(p) \oplus {\fr o}(q)) 
\subset {\fr o}(V)$ is  maximal compact (i.e.\ the Lie algebra
of a maximal compact subgroup of ${\rm O}(V) = {\rm O}(p,q)$)  and hence 
conjugated by an 
inner automorphism 
to the maximal compact subalgebra ${\fr o}(p) \oplus {\fr o}(q) 
\subset {\fr o}(V)$. So, again, we can assume that 
$\varphi$ preserves ${\fr o}(p) \oplus {\fr o}(q)$.  
Moreover, since $p\neq q$ any automorphism of 
${\fr o}(p) \oplus {\fr o}(q)$ is inner and we can assume that 
$\varphi |{\fr o}(p) \oplus {\fr o}(q) = {\rm Id}$. From the fact that
$\varphi$ is an automorphism of ${\fr p}(V)$ we obtain that $\phi :=
\varphi |V \in {\rm GL}(V)$ normalizes ${\fr o}(V)$ and 
${\fr o}(p) \oplus {\fr o}(q)$. This implies 
$\phi \in {\rm O}(p) \times {\rm O}(q)$, in particular,
$\varphi \mbox{\R}^p = \mbox{\R}^p$ and  $\varphi \mbox{\R}^{0,q} 
= \mbox{\R}^{0,q}$. Using an inner automorphism of ${\fr p}(\Pi )$
we can further assume that $\varphi |\mbox{\R}^p = \epsilon (\varphi )
{\rm Id}$, $\epsilon (\varphi ) \in \{ +1, -1\}$, and hence
$\varphi |{\fr o}(p) = {\rm Ad}_{\phi}|{\fr o}(p) = {\rm Id}$. This 
implies that $\varphi |W$ is ${\rm Spin}(p)$-equivariant. 
Now we can compute
\begin{eqnarray*}  b'(\varphi s,\varphi t)  &=& 
\langle e_1,[e_2 \ldots e_p\varphi s, \varphi t]\rangle =
\langle e_1,[\varphi e_2 \ldots e_ps,\varphi t]\rangle \\
&=& \langle e_1,\varphi [e_2 \ldots e_p s, t]\rangle =
\epsilon (\varphi ) \langle e_1,[e_2 \ldots e_p s, t]\rangle\\
&=& \epsilon (\varphi )b(s,t)\, ,\quad s,t\in W\, .
\end{eqnarray*}
Put $W_i' := \varphi W_i$. Then, since $\varphi^{\ast}b' = 
\epsilon (\varphi )b$,  we obtain a 
$b'$-orthogonal decomposition $W' = W_0' \oplus
\oplus _{i=1}^{l'}W_i'$ ($l' = l_+' +l_-'$) as above;
$W_0' = {\rm ker}\, \Pi' \cong l_0'\Sigma$, $W_i' \cong \Sigma$
($i\ge 1$) etc.\ This shows that $l_0' = l_0$, $l_{\pm}' = l_{\pm}$ if 
$\epsilon (\varphi ) = +1$ and $l_{\pm}' = l_{\mp}$ if 
$\epsilon (\varphi ) = -1$. $\Box$ 

Now we discuss the complementary case $p\equiv q\equiv 3\pmod{4}$. In this
case, the spinor module $S$ of $C\! \ell^0(V)$ is the sum $S^+ \oplus S^-$
of  two irreducible inequivalent semispinor modules $S^+$ and $S^-$ and
any irreducible $C\! \ell^0(V)$-module is equivalent to $S^+$ or $S^-$. 
As ${\rm Spin}_0(V)$-modules, $S^+$ and $S^-$ are dual: 
$S^- \cong (S^+)^{\ast}$. Let ${\fr p}(\Pi ) = {\fr p}(V) + W$ be an 
extended Poincar\'e algebra of signature $(p,q)$ as above. Thanks to
Thm.~\ref{decompIIThm} there exists a $b$-orthogonal decomposition 
$W = W_0 \oplus \oplus_{i=1}^l W_i$ as above with the following 
$b$-orthogonal refinements: $W_0 = W_0^+ \oplus W_0^-$, $W_0^{\pm} 
\cong l_0^{\pm}S^{\pm}$, and $\oplus_{i=1}^l W_i = 
\oplus_{i=1}^{l^+} W_i^+ \oplus  \oplus_{i=1}^{l^-} W_i^-$, 
$W_i^{\pm} \cong S^{\pm}$.  We denote by $l^{\epsilon}_+$ 
(respectively, $l^{\epsilon}_-$) the number of submodules 
$W_i^{\epsilon}$, $i = 1, \ldots , l^{\epsilon}$, on which $b$ is positive
(respectively, negative) definite, $\epsilon \in \{ +,-\}$. So to 
${\fr p}(\Pi )$ we have associated the nonnegative integers
$(l_0^+,l_+^+,l_-^+,l_0^-,l_+^-,l_-^-)$.

\begin{thm}\label{moduliThmII} Let $p \equiv q\equiv 3\pmod{4}$. 
Then the isomorphism class
$[{\fr p}(\Pi )]$ of an extended Poincar\'e algebra ${\fr p}(\Pi )$ 
of signature $(p,q)$ is completely determined by the tuple
$(l_0^+,l_+^+,l_-^+,l_0^-,l_+^-,l_-^-)$ introduced above. 
We put ${\fr p}(p,q,l_0^+,l_+^+,l_-^+,l_0^-,l_+^-,l_-^-) := 
[{\fr p}(\Pi )]$. Then  ${\fr p}(p,q,l_0^+,l_+^+,l_-^+,l_0^-,l_+^-,l_-^-) 
=  {\fr p}(p,q,\tilde{l}_0^+,\tilde{l}_+^+,\tilde{l}_-^+,\tilde{l}_0^-,
\tilde{l}_+^-,\tilde{l}_-^-)$ if and only if 
$(\tilde{l}_0^+,\tilde{l}_+^+,\tilde{l}_-^+,\tilde{l}_0^-,
\tilde{l}_+^-,\tilde{l}_-^-) \in  \Gamma 
(l_0^+,l_+^+,l_-^+,l_0^-,l_+^-,l_-^-)$, where $\Gamma \cong 
\mbox{\Z}_2 \times \mbox{\Z}_2$ is the group ge\-ne\-ra\-ted by the 
following two involutions:
\[ (l_0^+,l_+^+,l_-^+,l_0^-,l_+^-,l_-^-) \mapsto
(l_0^+,l_-^+,l_+^+,l_0^-,l_-^-,l_+^-)\] 
and 
\[ (l_0^+,l_+^+,l_-^+,l_0^-,l_+^-,l_-^-) \mapsto
(l_0^-,l_+^-,l_-^-,l_0^+,l_+^+,l_-^+)\, .\] 
\end{thm}

\noindent {\bf Proof:} 
The proof uses again Thm.~\ref{decompIIThm} and Thm.~\ref{dimThm} and is
similar to that of Thm.~\ref{moduliThmI}. Therefore, we explain only the reason
for the appearance of the second involution. In terms of the standard
basis $(e_1,\ldots ,e_p,e_1',\ldots ,e_q')$ of $V = \mbox{\R}^{p,q}$
we define an isometry $\phi \in {\rm SO}(p) \times {\rm O}(q)$ by:  
$\phi e_1 := e_1$, $\phi e_i := -e_i$ ($i\ge 2$) and 
$\phi e_j' = -e_j'$ ($j\ge 1$). Then ${\rm Ad}_{\phi} \in {\rm Aut} 
\,{\fr p}(V)$ induces an (outer) automorphism of ${\fr o}(V)$
interchanging the two semispinor modules. Let 
$({\fr p}(V) + W,[\cdot ,\cdot ]) \in {\fr p}(p,q,l_0^+,l_+^+,l_-^+,
l_0^-,l_+^-,l_-^-)$ be an extended Poincar\'e algebra of 
signature $(p,q)$. Then we define a new extended Poincar\'e algebra 
$({\fr p}(V) + W,[\cdot ,\cdot ]') \in {\fr p}(p,q,l_0^-,l_+^-,l_-^-,
l_0^+,l_+^+,l_-^+)$ by: 
\[ [\cdot ,\cdot ]' :=  [\cdot ,\cdot ] \quad \mbox{on} \quad
\wedge^2{\fr p}(V) \oplus \wedge^2W \oplus V\wedge W\]
and
\[ [A,s]' := [{\rm Ad}_{\phi}(A),s] \quad \mbox{for} \quad
A\in {\fr o}(V) \, , s\in W\, .\] 
The two $\mbox{\Z}_2$-graded Lie 
algebras are isomorphic via $\varphi : ({\fr p}(V) + W,[\cdot ,\cdot ]) 
\stackrel{\sim}{\rightarrow}  ({\fr p}(V) + W,[\cdot ,\cdot ]')$ defined
by: $\varphi |{\fr p}(V) = {\rm Ad}_{\phi}$ and $\varphi |W = {\rm Id}$. $\Box$
 
\section{The homogeneous quaternionic manifold $(M,Q)$ associated to
an extended Poincar\'e algebra}
\label{HQMSec} 
\subsection{Homogeneous manifolds associated to extended Poin\-ca\-r\'e
algebras} \label{homogSec}
Any extended Poincar\'e algebra ${\fr p} = {\fr p}(\Pi) =
{\fr p}(V) + W$ has an even derivation $D$ with eigenspace 
decomposition ${\fr p} = {\fr o}(V) + V + W$  and corresponding
eigenvalues $(0,1,1/2)$. Therefore, the $\mbox{\Z}_2$-graded Lie 
algebra ${\fr p} = {\fr p}_0 + {\fr p}_1 = {\fr p}(V) + W$ is canonically
extended to a $\mbox{\Z}_2$-graded Lie 
algebra ${\fr g} = {\fr g}(\Pi ) = \mbox{\R} D + {\fr p} =
{\fr g}_0 + {\fr g}_1$, where ${\fr g}_0 = \mbox{\R} D + {\fr p}_0 =
\mbox{\R} D + {\fr p}(V) =: {\fr g}(V)$ and ${\fr g}_1 =
{\fr p}_1 = W$. The next proposition describes the basic structure
of the Lie algebra ${\fr g} = {\fr g}(\Pi ) = {\fr g}(V) + W$. 
To avoid trivial exceptions, in the following we assume that 
$n = \dim V >2$ and hence that ${\fr o}(V)$ is semisimple  
(for the construction  of homogeneous quaternionic manifolds
we will put $V = \mbox{\R}^{p,q}$ and $p\ge 3$). 

\begin{prop} \label{LeviProp} The Lie algebra ${\fr g} = {\fr g}(\Pi )$
has the Levi decomposition 
\begin{equation} \label{LeviEqu} 
{\fr g} = {\fr o}(V) + {\fr r}
\end{equation}
into the radical ${\fr r} = \mbox{\R}D + V +W$ and the complementary
(maximal semisimple) Levi subalgebra ${\fr o}(V)$. The nilradical
${\fr n} = [{\fr r},{\fr r}] = V + W$ is two-step  nilpotent if
$\Pi \neq 0$ and Abelian otherwise. 
\end{prop} 

{}For any Lie algebra 
$\fr l$ we denote by ${\rm der} ({\fr l})$ the Lie algebra
of its derivations.   
\begin{prop} \label{faithProp} The adjoint representation
${\fr g} \rightarrow  {\rm der} ({\fr r})$ of ${\fr g} = 
{\fr o}(V) + {\fr r}$ on its ideal ${\fr r} = \mbox{\R}D + V +W$
is faithful.
\end{prop}

\noindent {\bf Proof:} Let $x\in {\fr g}$. We show that $[x, {\fr r}] = 0$ 
implies $x = 0$. First $[x,D] = 0$ implies $x \in {\fr co}(V) = 
\mbox{\R}D + {\fr o}(V)$. Then $[x,V] = 0$ implies $x=0$, because
the conformal Lie algebra ${\fr co}(V)$ acts faithfully on $V$. $\Box$ 

By Prop.~\ref{faithProp} we can consider ${\fr g}(V)$ and ${\fr g} =
{\fr g}(\Pi ) = {\fr g}(V) + W$ as linear Lie algebras via the embedding
${\fr g}(V) \subset {\fr g} \hookrightarrow {\rm der} ({\fr r}) \subset 
{\fr gl} ({\fr r})$. The corresponding connected Lie groups of 
${\rm Aut} ({\fr r}) \subset {\rm GL}({\fr r})$ will be denoted
by $G(V)$ and $G = G(\Pi )$ respectively: ${\rm Lie}\, G(V) = 
{\fr g}(V) \subset {\rm Lie}\, G = {\fr g}$. Now let 
$V = \mbox{\R}^{p,q}$ and fix a decomposition $p = p' + p''$. 
The subalgebra of ${\fr o}(V)$ preserving the corresponding orthogonal
splitting $V = E + E' = \mbox{\R}^{p',0} +  \mbox{\R}^{p'',q}$ is 
${\fr k} = {\fr k}(p',p'') = {\fr o}(p') \oplus {\fr o}(p'',q) 
\subset {\fr o}(p,q) = {\fr o}(V)$. We consider ${\fr k} 
\subset {\fr o}(V) \subset {\fr g}$ as a subalgebra of the
linear Lie algebra ${\fr g} \hookrightarrow {\rm der} ({\fr r})$ and denote
by $K = K(p',p'') \subset G(V) \subset G \subset {\rm Aut} ({\fr r})$
the corresponding connected linear Lie group. $K$ is a closed Lie subgroup, 
see Cor.\ \ref{closedCor} below.  We are interested
in the homogeneous spaces
\[ M(V) := G(V)/K \subset M=M(\Pi ) := G/K = G(\Pi )/K\, .\]

\begin{prop}
\label{algProp} The Lie subalgebras ${\fr k}$, ${\fr n}$, ${\fr r}$, 
${\fr o}(V)$, ${\fr p}(V)$, $\fr p$, ${\fr g}(V)$, ${\fr g} \subset
{\rm der} ({\fr r})$ are algebraic subalgebras of the (real) algebraic
Lie algebra ${\rm der}({\fr r})$. 
\end{prop}

\noindent {\bf Proof:} We use the following sufficient conditions
for algebraicity, see \cite{O-V} Ch.\ 3 \S 3 8$^{\rm o}$:
\begin{enumerate}
\item[a)] A linear Lie algebra coinciding with its derived algebra
is algebraic.
\item[b)] The radical of an algebraic linear Lie algebra is algebraic.
\item[c)] A linear Lie algebra generated by algebraic subalgebras is
algebraic. 
\end{enumerate}
The subalgebras ${\fr o}(V) \subset {\fr p}(V) \subset {\fr p} \subset 
{\rm der} ({\fr r})$ are algebraic by a). 
By c), to prove that ${\fr g} = \mbox{\R}D + {\fr p}$
and ${\fr g}(V)$ are algebraic it is now sufficient to show that 
$\mbox{\R}D \subset {\rm der}({\fr r})$ is algebraic. $D$ 
preserves ${\fr n} = V+W$ and acts trivially on the complement
$\mbox{\R}D + {\fr o}(V)$. The Lie algebra $\mbox{\R}D 
\hookrightarrow {\rm der}({\fr n})$ is the Lie algebra
of the algebraic group
\[ \{ \lambda {\rm Id}_V \oplus \mu {\rm Id}_W|\lambda = \mu^2 \neq 0\}
\subset {\rm GL}(V\oplus W)\, .\]  
This shows that $\mbox{\R}D$, ${\fr g}(V)$, ${\fr g} \subset
{\rm der}({\fr r})$ are algebraic.  Now ${\fr n}$, the radical of 
$\fr p$, and $\fr r$, the radical of $\fr g$, are algebraic by b). 
Finally, ${\fr k} \subset {\fr o}(V)$
is the subalgebra which preserves the orthogonal splitting 
$V = E + E'$ and is
hence algebraic. $\Box$ 

\begin{cor}
\label{closedCor} The connected linear Lie groups $K$, $S \cong 
{\rm Spin}_0(V)$, $R = \exp {\fr r}$, $G(V)$, $G \subset {\rm Aut}({\fr r})$
with Lie algebras $\fr k$, ${\fr s} := {\fr o}(V)$, $\fr r$, ${\fr g}(V)$,
${\fr g} \subset {\rm der}({\fr r})$  are closed. 
\end{cor}

\noindent {\bf Proof:} This follows from Prop.\ \ref{algProp} and the fact
that the identity component of a real algebraic linear group is a closed
linear group. $\Box$ 
\begin{prop} The Lie group $G(\Pi )$ has the following Levi decomposition:
\begin{equation} G(\Pi ) = S \mbox{\Bbb n} R \, .
\label{LeviIIEqu} \end{equation}
\end{prop} 

\noindent {\bf Proof:} This follows from the corresponding Levi
decomposition (\ref{LeviEqu}) of Lie algebras since $S \cap R = \{ e\}$. 
$\Box$ 
 
Next we show that $M(V)$ can be naturally endowed with a $G(V)$-invariant 
structure of pseudo-Riemannian locally symmetric space. With 
this in mind, 
we consider the pseudo-Euclidean vector space 
$\tilde{V} = \mbox{\R}^{p+1,q+1}$ with
scalar product $\langle \cdot ,\cdot \rangle$ and orthogonal decomposition
$\tilde{V} = \mbox{\R}e_0 + V + \mbox{\R}e_0'$, $\langle e_0,e_0\rangle =
- \langle e_0',e_0'\rangle = 1$. Recall that ${\fr o}(\tilde{V})$ is 
identified with $\wedge^2\tilde{V}$ via the pseudo-Euclidean 
scalar product  $\langle \cdot ,\cdot \rangle$, see (\ref{wedgeEqu}). 

\begin{prop}\label{embProp} The subspace
\[  \mbox{\R}e_0\wedge e_0' + \wedge^2V + (e_0-e_0')\wedge V \subset
\wedge^2\tilde{V} = {\fr o}(\tilde{V})\] 
is a subalgebra isomorphic
to ${\fr g}(V)$. 
\end{prop}

\noindent {\bf Proof:} The canonical embedding ${\fr o}(V) = 
\wedge^2V \hookrightarrow \wedge^2\tilde{V} = {\fr o}(\tilde{V})$
is extended to an embedding ${\fr g}(V) \hookrightarrow 
{\fr o}(\tilde{V})$ via 
\[ D\mapsto e_0\wedge e_0'\, ,\quad V\ni v \mapsto (e_0-e_0')\wedge v\, .
\Box \]

It is easy to see that the embedding ${\fr g}(V) \hookrightarrow
{\fr o}(\tilde{V})$ lifts to a homomorphism of connected
Lie groups $G(V) \rightarrow {\rm SO}_0(\tilde{V})$ with
finite kernel. In particular, we have a natural isometric
action of $G(V)$ on the pseudo-Riemannian symmetric space
$\tilde{M}(V) := {\rm SO}_0(p+1,q+1)/({\rm SO}(p'+1) \times
{\rm SO}_0(p'',q+1))$. We  denote by $[e] = e\tilde{K}  \in \tilde{M}(V)
:= {\rm SO}_0(\tilde{V})/\tilde{K}$, $\tilde{K} := 
{\rm SO}(p'+1) \times {\rm SO}_0(p'',q+1)$, the canonical base point
and by $G(V)[e]$ its $G(V)$-orbit.  
\begin{prop}\label{stProp} The action of $G(V)$ on $\tilde{M}(V)$ induces a 
$G(V)$-equivariant open embedding 
\[ M(V) = G(V)/K \stackrel{\sim}{\rightarrow} G(V)[e] \subset 
\tilde{M}(V)\, .\]
$G(V)$ acts transitively on $\tilde{M}(V)$ if and only if 
$\tilde{M}(V)$ is Riemannian, i.e.\ if and only if $p'' = 0$. 
In that case $M(V) \cong \tilde{M}(V)$ is the noncompact dual of 
the Grassmannian ${\rm SO}(p+q+2)/({\rm SO}(p+1)\times {\rm SO}(q+1))$
and admits a simply transitive splittable solvable subgroup
$I(S) \mbox{\Bbb n} \exp (\mbox{\R}D +V) \subset G(V)$. Here
$I(S)$ denotes the (solvable) Iwasawa subgroup of $S \cong {\rm Spin}_0(V)$.
\end{prop}

\noindent {\bf Proof:} The stabilizer $G(V)_{[e]}$ of $[e]$ in $G(V)$ coincides
with $K$ and hence $G(V)[e] \cong G(V)/G(V)_{[e]} = G(V)/K = M(V)$. Now 
a simple dimension count shows that the orbit $G(V)[e] \subset 
\tilde{M}(V)$ is open. If $\tilde{M}(V)$ is Riemannian then it is a
Riemannian symmetric space of noncompact type and, by the Iwasawa
decomposition theorem (see \cite{H}) the Iwasawa subgroup
$I({\rm SO}_0(p+1,q+1)) \subset {\rm SO}_0(p+1,q+1)$ is a splittable
solvable subgroup which acts simply transitively on $\tilde{M}(V)$.
Now let $G(V) = S \mbox{\Bbb n} R(V)$ be the Levi decomposition
associated to the Levi decomposition ${\fr g}(V) = {\fr o}(V) +
(\mbox{\R}D + V)$, $S \cong {\rm Spin}_0(V)$, 
$R(V) = \exp (\mbox{\R}D + V)$  (cf.\ 
Prop.~\ref{LeviProp}). We denote by $I(S) \subset S$ 
the Iwasawa subgroup of $S$. Then $I(S) \mbox{\Bbb n} R(V) \subset G(V)$
is mapped isomorphically onto $I({\rm SO}_0(p+1,q+1)) \subset 
{\rm SO}_0(\tilde{V})$ by the homomorphism $G(V) \rightarrow
{\rm SO}_0(\tilde{V})$ introduced above. This shows that 
$I(S) \mbox{\Bbb n} R(V)$ and hence $G(V)$ acts transitively on
$\tilde{M}(V)$ in the Riemannian case. 

If  $\tilde{M}(V)$ is not Riemannian then the homogeneous
spaces $\tilde{M}(V)$ and $M(V)$ are not homotopy equivalent
and hence $G(V)$ does not act transitively on $\tilde{M}(V)$. 
In fact, $\tilde{M}(V)$ (respectively, $M(V)$) has the 
homotopy type of the Grassmanian ${\rm SO}(p+1)/(
{\rm SO}(p'+1)\times {\rm SO}(p''))$ (respectively, ${\rm SO}(p)/
({\rm SO}(p')\times {\rm SO}(p''))$). Now it is sufficient
to observe that the Stiefel manifolds ${\rm SO}(p+1)/
{\rm SO}(p'+1)$ and ${\rm SO}(p)/
{\rm SO}(p')$ are homotopy equivalent only if $p=p'$ and hence
$p''=0$, see \cite{O3}. $\Box$

Before we can formulate the main result of the present paper in section
\ref{MainSec} we need to recall the notions of quaternionic manifold and 
of (pseudo-) quaternionic K\"ahler manifold, cf.\ \cite{A-M2}. 
The reader familiar with these concepts should skip the next section. 
\subsection{Quaternionic structures}
\label{QSSec}
It is instructive to introduce the basic concepts of quaternionic
geometry as analogues of the more familiar concepts of complex 
geometry. 
\begin{dof} Let $E$ be a (finite dimensional) real vector space.
A {\bf complex structure} on $E$ is an endomorphism $J \in {\rm End}(E)$ 
such that $J^2 = -{\rm Id}$. A {\bf hypercomplex structure} on $E$
is a triple $(J_{\alpha}) =(J_1,J_2,J_3)$ of complex structures on $E$
satisfying $J_1J_2 = J_3$. A {\bf quaternionic 
structure} on $E$ is the three-dimensional subspace $Q \subset {\rm End}(E)$ 
spanned by a hypercomplex structure $(J_{\alpha})$: $Q = {\rm span}
\{ J_1,J_2,J_3\}$. In that case, we say that the hypercomplex
structure $(J_{\alpha})$ is {\bf subordinate} to the quaternionic
structure $Q$. 
\end{dof}
Note, first, that $Q \subset {\fr gl}(E)$ is a Lie subalgebra isomorphic
to ${\fr sp}(1) \cong {\rm Im} \mbox{\Bbb H} = {\rm span} \{ i,j,k \}$ 
the Lie algebra of the group
$Sp(1) = S^3 \subset \mbox{\Bbb H} = {\rm span} \{ 1, i,j,k \}$ of unit 
quaternions and, second,  
that the real associative subalgebra of ${\rm End}(E)$ generated
by ${\rm Id}$ and a quaternionic structure $Q$ on $E$ is isomorphic
to the algebra of quaternions $\mbox{\Bbb H}$.
In both cases, the choice of such isomorphism is equivalent
to the choice of a hypercomplex structure $(J_{\alpha}) = (J_1,J_2,J_3)$ 
subordinate to $Q$. In fact, given $(J_{\alpha})$ we can define an 
isomorphism of associative algebras by $({\rm Id},J_1,J_2,J_3) \mapsto 
(1,i,j,k)$ this induces also an isomorphism of Lie algebras
$Q \stackrel{\sim}{\rightarrow} {\fr sp}(1)$. 

\begin{dof}\label{acsDef}
Let $M$ be a (smooth) manifold. An {\bf almost complex structure} $J$
(respectively, {\bf almost hypercomplex structure} $(J_1,J_2,J_3)$,
{\bf almost quaternionic structure} $Q$) on $M$ is a (smooth) 
field $M \ni m \mapsto J_m \in {\rm End}(T_mM)$ of complex structures 
(respectively,  $m \mapsto (J_1,J_2,J_3)_m$ of hypercomplex structures,
$m \mapsto Q_m$ of quaternionic structures). The pair $(M,J)$ (respectively, 
$(M,(J_{\alpha}))$, $(M,Q))$ is called an {\bf almost complex manifold} 
(respectively, {\bf almost hypercomplex manifold}, {\bf almost 
quaternionic manifold}).   
\end{dof}

We recall that a connection on a manifold $M$ (i.e.\ a covariant derivative
$\nabla$ on its tangent bundle $TM$) induces a covariant derivative $\nabla$
on the full tensor algebra over $TM$ and, in particular, on ${\rm End}(TM)
\cong TM \otimes T^{\ast}M$. We will say that $\nabla$ preserves a subbundle
$B\subset {\rm End}(TM)$ if it maps (smooth) sections of $B$ into
sections of $T^{\ast}M \otimes B$. 

\begin{dof}
\label{connDef} A connection $\nabla$ on an almost complex manifold 
$(M,J)$ (respectively, almost hypercomplex manifold $(M,(J_{\alpha}))$, 
almost quaternionic manifold $(M,Q))$ is called {\bf almost complex} 
(respectively, {\bf almost hypercomplex}, {\bf almost quaternionic})
{\bf connection} if $\nabla J =0$ (respectively, if $\nabla J_1 = 
\nabla J_2 = \nabla J_3 = 0$, if $\nabla$ preserves the rank 3 subbundle
$Q  \subset {\rm End}(TM))$. A {\bf complex} (respectively, 
{\bf hypercomplex}, {\bf quaternionic}) {\bf connection}  is a
torsionfree  almost complex (respectively, almost hypercomplex, 
almost quaternionic) connection. 
\end{dof}
Note that the equation $\nabla J = 0$ is equivalent to the condition
that $\nabla$ preserves the rank 1 subbundle of  ${\rm End}(TM)$ spanned
by $J$, i.e.\ $\nabla J = \theta \otimes J$ for some 1-form $\theta$ on $M$. 
Therefore, the notion of 
almost quaternionic connection (as well as that of almost hypercomplex 
connection) is a direct quaternionic analogue of the 
notion of almost complex connection. 
\begin{dof}
\label{integrDef} Let $M$ be a manifold. An almost complex structure 
(respectively, almost hypercomplex structure, 
almost quaternionic structure) on $M$ is called {\bf 1-in\-te\-gra\-ble}  
if there exists a complex (respectively, hypercomplex, quaternionic)
connection on $M$. A {\bf complex structure} (respectively,
{\bf hypercomplex structure}, {\bf quaternionic structure}) on $M$ 
is a 1-integrable almost complex structure (respectively, 
almost hypercomplex structure, almost quaternionic structure) on $M$. 
\end{dof}

\noindent 
{\bf Remark 3:} It is well known, see \cite{N-N} and \cite{K-NII}, 
that an almost complex manifold $(M,J)$ is integrable, i.e.\ 
admits an atlas with holomorphic transition maps, if and only if it is
1-integrable. This justifies the following definition of complex
manifold.   

\begin{dof} \label{quatmfDef}
A {\bf complex manifold} (respectively,  {\bf hypercomplex manifold})
is a manifold $M$ together with a complex structure $J$ (respectively,  
hypercomplex structure $(J_{\alpha}))$ on $M$. A {\bf quaternionic manifold}
of dimension $d>4$ is a manifold $M$ of dimension $d$ together with
a quaternionic structure $Q$ on $M$. Finally, a {\bf quaternionic manifold}
of dimension $d = 4$ is a 4-dimensional manifold $M$ together with an 
almost  
quaternionic structure $Q$ which annihilates the Weyl tensor of 
the conformal structure defined by $Q$, see
Remark 4 below. 
\end{dof}
{}For examples of hypercomplex and quaternionic manifolds (without
metric condition) see \cite{J1}, \cite{J2} and \cite{B-D}.  

\noindent 
{\bf Remark 4:} Notice that an almost quaternionic structure on a 
4-manifold induces an oriented conformal structure. This follows
from the fact that the normalizer in ${\rm GL}(4,\mbox{\R})$ of the 
standard quaternionic structure(s) on $\mbox{\R}^4 = \mbox{\Ha}$ is the
special conformal group ${\rm CO}_0(4)$. The definition of 
quaternionic 4-manifold $(M,Q)$ implies that this conformal structure
is half-flat. More precisely, if the orientation of $M$ is 
chosen such that $Q$ is locally generated by  positively oriented almost
complex structures, then the self-dual half of the  Weyl tensor vanishes. 
The special treatment
of the 4-dimensional case in Def.~\ref{quatmfDef} and also in 
Def.~\ref{KaehlerDef} below has the advantage that with these 
definitions many important properties of quaternionic 
(and also of  quaternionic K\"ahler) manifolds of dimension $>4$ 
remain true in dimension 4, cf.\ Remark 5. In particular, all future statements
about quaternionic manifolds and quaternionic K\"ahler manifolds in the
present paper, such as the integrability
of the canonical almost complex structure on the twistor space (see section
\ref{bundleSec}), are valid also in the 4-dimensional case. 

Next we discuss (almost) complex, hypercomplex and quaternionic
structures on a pseudo-Riemannian manifold $(M,g)$. 

\begin{dof}
A pseudo-Riemannian metric $g$ on an almost complex manifold
$(M,J)$ (respectively, almost hypercomplex manifold $(M,(J_{\alpha}))$,
almost quaternionic manifold $(M,Q))$ is called {\bf Hermitian} if
$J$ is skew symmetric (respectively, the $J_{\alpha}$ are skew symmetric,
$Q$ consists of skew symmetric endomorphisms). 
\end{dof}
Note that, due to $J^2 = -{\rm Id}$, an almost complex structure
$J$ on a pseudo-Riemannian manifold $(M,g)$ is skew symmetric
if and only if $J$ is orthogonal, i.e.\ if and only if $g(JX,JY) =
g(X,Y)$ for all vector fields $X,Y$ on $M$. Similarly, an almost
quaternionic structure $Q$ on a pseudo-Riemannian manifold $(M,g)$
consists of skew symmetric endomorphisms if and only if 
$Z := \{ A \in Q| A^2 = -{\rm Id}\}$ consists of orthogonal
endomorphisms (here the equation $A^2 = -{\rm Id}$ is on 
$T_{\pi A}M$, $\pi :Q\rightarrow M$ the bundle projection). 
\begin{dof} \label{HermDef}
A {\bf complex} (pseudo-) {\bf Hermitian manifold} $(M,J,g)$ 
(respectively, {\bf hypercomplex} (pseudo-) {\bf Hermitian manifold} 
$(M,(J_{\alpha}),g)$, {\bf quaternionic}  (pseu\-do-) {\bf Hermitian manifold} 
$(M,Q,g)$, {\bf almost complex} (pseudo-) {\bf Hermitian manifold} 
$(M,J,g)$ etc.) is a complex manifold $(M,J)$ 
(respectively,  hypercomplex manifold
$(M,(J_{\alpha}))$, quaternionic manifold $(M,Q)$, almost complex
manifold $(M,J)$ etc.) with a Hermitian (pseudo-) Riemannian metric $g$. 
\end{dof}
Next we define the hypercomplex and quaternionic analogues of 
(pseudo-) K\"ahler manifolds. 
\begin{dof}
\label{KaehlerDef} A  (pseudo-) {\bf K\"ahler manifold} (respectively,
(pseudo-) {\bf hyper-K\"ah\-ler manifold},   {\bf quaternionic} (pseudo-) 
{\bf K\"ahler manifold} of dimension $d>4$) is an almost complex 
(pseudo-) Hermitian manifold $(M,J,g)$ (respectively,
almost hypercomplex (pseudo-) Hermitian manifold $(M,(J_{\alpha}),g)$,
almost quaternionic  (pseudo-)\linebreak[3] Hermitian manifold $(M,Q,g)$ of 
dimension $d>4$) with the property that the Levi-Civita connection
$\nabla^g$ of the (pseudo-) Riemannian metric $g$ is complex
(respectively, hypercomplex, quaternionic). An almost  
quaternionic Hermitian 4-manifold $(M,Q,g)$ is called {\bf quaternionic 
K\"ahler manifold} if $Q$ annihilates the curvature tensor $R$ of 
$\nabla^g$. 
\end{dof}

See \cite{Bes}, \cite{H-K-L-R}, \cite{C-F-G}, \cite{Sw1}, \cite{D-S1},
\cite{D-S2}, \cite{Bi1}, \cite{Bi2}, \cite{K-S2}  and \cite{C4}  
for examples of hyper-K\"ahler manifolds and \cite{W1}, \cite{A3}, 
\cite{Bes}, \cite{G1}, \cite{G-L}, \cite{L2}, \cite{L3}, 
\cite{F-S}, \cite{dW-VP2},  \cite{A-G}, 
\cite{K-S1}, \cite{A-P}, \cite{D-S3} and 
\cite{C2}  for examples of quaternionic K\"ahler manifolds. 

\noindent
{\bf Remark 5:} As explained in Remark 4, an almost quaternionic structure 
$Q$ on a 4-manifold $M$ defines a conformal structure. It is clear that a 
pseudo-Riemannian metric $g$ on $M$ defines the same conformal structure as 
$Q$ if and only if it is $Q$-Hermitian and that any such metric is definite.
Moreover, the Levi-Civita connection $\nabla^g$ of an almost quaternionic 
Hermitian 4-manifold $(M,Q,g)$ automatically preserves $Q$.  
In fact, its holonomy group at $m\in M$ is a subgroup of 
${\rm SO}(T_mM,g_m)$ because $M$ is oriented (by $Q$) and
the latter group normalizes $Q$ because $g$ is $Q$-Hermitian. 
Now let $(M,Q,g)$ be a quaternionic K\"ahler 4-manifold. Since $Q$
annihilates the curvature tensor $R$ and the metric $g$ it must also
annihilate the the 
Ricci tensor $Ric$ and the Weyl tensor of $(M,g)$. This shows first that  
$Ric$ is $Q$-Hermitian and hence proportional to $g$: $Ric = cg$. 
In other words, $(M,g)$ is an Einstein manifold. Second, $(M,Q)$
is a quaternionic manifold because $Q$ annihilates the Weyl tensor
of the conformal structure defined by $g$, which coincides with
the conformal structure defined by $Q$. For any quaternionic
pseudo-K\"ahler manifold $(M,g)$ (of arbitrary dimension) 
it is known, see e.g.\ 
\cite{A-M2}, that $(M,g)$ is Einstein and that $Q$ annihilates $R$.

As next, we introduce the appropriate notions in order to discuss 
transitive group actions on manifolds with the special geometric structures
defined above.

\subsection{Invariant connections on homogeneous manifolds and 1-integrability
of homogeneous almost quaternionic manifolds} 
\label{ICHMSec}
\begin{dof}
The {\bf automorphism group} of an almost complex manifold $(M,J)$ 
(respectively, almost hypercomplex manifold $(M,(J_{\alpha}))$,
almost quaternionic  manifold $(M,Q)$, almost complex Hermitian
manifold $(M,J,g)$, etc.) is the group of diffeomorphisms of M
which preserves $J$ (respectively, $(J_{\alpha})$, $Q$, $(J,g)$, etc.). 
An   almost complex manifold (respectively, almost hypercomplex manifold,
almost quaternionic  manifold, almost complex Hermitian
manifold, etc.) is called {\bf homogeneous} if it has a transitive 
automorphism group. 
\end{dof}
In the next section, we are going to construct an 
almost quaternionic structure $Q$ on certain homogeneous 
manifolds $M = G/K$ ($G$ is a Lie group and $K$ a closed subgroup). 
The almost quaternionic structure $Q$ will be $G$-invariant
by construction and hence $(M,Q)$ will be a homogeneous
almost quaternionic manifold. Similarly, we will construct
homogeneous almost quaternionic (pseudo-) Hermitian manifolds
$(M = G/K,Q,g)$. In order to prove that $(M = G/K,Q,g)$ is 
a quaternionic (pseudo-) K\"ahler manifold, or simply to
establish the 1-integrability of the almost quaternionic
structure $Q$ it is useful to have an appropriate description
of the affine space of $G$-invariant connections on the
homogeneous manifold $M = G/K$. This is provided by the notion
of Nomizu map, which we now recall, see \cite{A-V-L}, \cite{V1}. 
Let $\nabla$ be a connection on a manifold $M$. For any vector
field $X$ on $M$ one defines the operator 
\begin{equation} \label{Nomizu0Equ} L_X := {\cal L}_X -\nabla_X \, ,
\end{equation}
where ${\cal L}_X$ is the Lie derivative (i.e.\ ${\cal L}_XY = [X,Y]$ for
any vector field $Y$on $M$). $L_X$ is a $C^{\infty}(M)$-linear map
on the $C^{\infty}(M)$-module  $\Gamma (TM)$ of vector fields on $M$,
so it can be identified with a section $L_X \in \Gamma ({\rm End}(TM))$.
In particular, $L_X|_m \in {\rm End}(T_mM)$ for all $m\in M$. 

Now let $M = G/K$ be a homogeneous space and suppose that $\nabla$ is 
$G$-invariant. The action of $G$ on $M$ defines an antihomomorphism of
Lie algebras $\alpha :{\fr g} = {\rm Lie} \, G \rightarrow \Gamma (TM)$
from the Lie algebra $\fr g$ of {\em left}-invariant vector fields
on $G$ to the Lie algebra $\Gamma (TM)$ of vector fields on $M$. 
This antihomomorphism maps an element $x\in {\fr g}$  to the 
{\bf fundamental vector 
field}\label{fundvfDef} $\alpha (x):= X \in\Gamma (TM)$ defined by 
$X(m) := \frac{d}{dt}|_{t=0} \exp (tx)m$. The statement that $\alpha$
is an antihomomorphism means that $[\alpha (x),\alpha (y)] = -\alpha
([x,y])$ for all $x,y \in {\fr g}$. Note that $\alpha$ becomes a 
homomorphism if we replace $\fr g$ by the Lie algebra of 
{\em right}-invariant vector fields on $G$. 
Without restriction of generality, we can assume that the action
is almost effective, i.e.\ ${\fr g} \rightarrow \Gamma (TM)$ is 
injective. Then we can identify $\fr g$ with its faithful image in 
$\Gamma (TM)$.  The isotropy subalgebra ${\fr k} = {\rm Lie} \, K$ is mapped
(anti)isomorphically onto a subalgebra of vector fields vanishing at the base
point $[e] = eK \in M =G/K$. In this situation we define the {\bf Nomizu
map} $L = L(\nabla ): {\fr g} \rightarrow {\rm End}(T_{[e]}M)$,    
$x\mapsto L_x$, by the equation
\[ L_x := L_X|_{[e]}\, ,\] 
where again $X$ is the fundamental vector field on $M$ associated to $x \in {\fr g}$. 
The ope\-ra\-tors $L_x \in {\rm End}(T_{[e]}M)$ will be called 
{\bf Nomizu operators}. They have the following properties:
\begin{equation} \label{NomizuIEqu} 
L_x = d\rho (x) \quad \mbox{for all} \quad x \in {\fr k}  
\end{equation}
and 
\begin{equation} \label{NomizuIIEqu} 
L_{{\rm Ad}_kx} = \rho (k) L_x \rho (k)^{-1} \quad \mbox{for all} 
\quad x \in {\fr g}\, , \quad k\in K\, , 
\end{equation}
where $\rho :K \rightarrow {\rm GL}(T_{[e]}M)$ is the isotropy 
representation. The first equation follows directly from equation 
(\ref{Nomizu0Equ}) since $(\nabla_XY)_m = 0$ if $X(m) = 0$. The 
second equation expresses the $G$-invariance of $\nabla$. 
Conversely, any linear map $L : {\fr g} \rightarrow {\rm End}(T_{[e]}M)$ 
satisfying (\ref{NomizuIEqu}) and (\ref{NomizuIIEqu})
is the Nomizu map of  a uniquely defined $G$-invariant connection
$\nabla = \nabla (L)$ on $M$. Its torsion tensor $T$ and curvature tensor
$R$  are expressed at $[e]$ by: 
\[ T(\pi x, \pi y) = -(L_x \pi y - L_y \pi x + \pi [x,y]) \]
and 
\[ R(\pi x, \pi y) = [L_x,L_y]  + L_{[x,y]} \, ,\quad x,y\in {\fr g}\]
where $\pi  : {\fr g} \rightarrow T_{[e]}M$ is the  canonical projection
$x \mapsto \pi x =  X([e]) = \frac{d}{dt}|_{t=0} \exp (tx)K$. 

\noindent {\bf Remark 6:} The difference between our formulas 
for torsion and curvature and 
those of \cite{A-V-L} is due to the fact that in  \cite{A-V-L}
everything is expressed in terms of right-invariant vector 
fields on the Lie group $G$ whereas we use left-invariant vector fields. 
To obtain the corresponding expressions 
in terms of right-invariant vector fields from the expressions in terms
of left-invariant vector fields and vice versa it is sufficient to 
replace $[x,y]$ by $-[x,y]$ in all formulas ($x,y \in {\fr g}$). 
The same remark applies to formula (\ref{KoszulEqu}) below, which
expresses the Levi-Civita connection on a pseudo-Riemannian homogeneous
manifold. 

Suppose now that we are given a $G$-invariant geometric structure
$S$ on $M$ (e.g.\ a $G$-invariant almost quaternionic structure $Q$)
defined by a corresponding $K$-invariant geometric structure
$S_{[e]}$ on $T_{[e]}M$. Then a $G$-invariant connection $\nabla$
preserves $S$ if and only if the corresponding Nomizu operators
$L_x$, $x\in {\fr g}$, preserve $S_{[e]}$. So to construct a 
$G$-invariant connection preserving $S$ it is sufficient to
find a Nomizu map $L: {\fr g} \rightarrow {\rm End} (T_{[e]}M)$
such that $L_x$ preserves $S_{[e]}$ for all $x\in {\fr g}$. 
We observe that, due to the $K$-invariance of $S_{[e]}$, the Nomizu
operators $L_x$ preserve $S_{[e]}$ already for $x\in {\fr k}$. 
The above considerations can be specialized as follows: 
\begin{prop}
Let $Q$ be a $G$-invariant almost quaternionic structure on a
homogeneous manifold $M = G/K$. There is a natural one-to-one
correspondence between $G$-invariant almost quaternionic connections on 
$(M,Q)$ and Nomizu maps $L: {\fr g} \rightarrow {\rm End} (T_{[e]}M)$,
whose image normalizes $Q_{[e]}$, i.e.\ whose Nomizu operators $L_x$,
$x\in {\fr g}$, belong to the normalizer ${\rm n}(Q) \cong {\fr sp}(1)
\oplus {\fr gl}(d,\mbox{\Ha})$ ($d =\dim M/4$) of the 
quaternionic structure $Q_{[e]}$ in the Lie algebra ${\fr gl}(T_{[e]}M)$. 
\end{prop}
\begin{cor}
\label{quatCor} Let $(M = G/K,Q)$ be a homogeneous almost quaternionic
manifold and $L: {\fr g} \rightarrow {\rm End} (T_{[e]}M)$ a Nomizu
map such that 
\begin{enumerate}
\item[(1)] $L_x\pi y -L_y\pi x = -\pi [x,y]$ for all $x,y \in {\fr g}$
(i.e.\ $T=0$) and 
\item[(2)] $L_x$ normalizes $Q_{[e]}\subset {\rm End} (T_{[e]}M)$.
\end{enumerate}
Then the connection $\nabla (L)$ associated to the Nomizu map $L$ 
is a $G$-invariant quaternionic connection on
$(M,Q)$ and hence $Q$ is 1-in\-te\-gra\-ble. 
\end{cor}
{}For future use we give the well known formula for the Nomizu map
$L^g$ associated to the Levi-Civita connection $\nabla^g$ of a
$G$-invariant pseudo-Riemannian metric $g$ on a homogeneous
space $M = G/K$. Let $\langle \cdot ,\cdot \rangle = g_{[e]}$
be the $K$-invariant scalar product on $T_{[e]}M$ induced by $g$.
Then $L_x^g \in {\rm End} (T_{[e]}M)$, $x\in {\fr g}$, is given by the 
following Koszul type formula:
\begin{equation}\label{KoszulEqu}
-2\langle L_x^g\pi y,\pi z\rangle = 
\langle \pi [x,y],\pi z\rangle - \langle \pi x, \pi [y,z]\rangle 
- \langle \pi y, \pi [x,z]\rangle 
\, ,\quad x,y,z \in {\fr g}\, .
\end{equation}
\begin{cor}
\label{qKCor} Let $(M = G/K,Q,g)$ be a homogeneous almost quaternionic
(pseudo-) Hermitian manifold and assume that $L_x^g$ normalizes
$Q_{[e]}$ for all $x\in {\fr g}$. Then the Levi-Civita connection
$\nabla^g = \nabla (L^g)$ is a $G$ invariant quaternionic connection
on $(M,Q,g)$ and hence $(M,Q,g)$ is a  quaternionic (pseudo-)
K\"ahler manifold if $\dim M >4$. 
\end{cor}
\subsection{The main theorem} \label{MainSec}
Let ${\fr p}(\Pi ) = {\fr p}(V) + W$ be an 
extended Poincar\'e algebra of signature $(p,q)$ and $p\ge 3$. We fix the 
decomposition $p = p' + p''$, where $p'=3$. For notational convenience
we put $r := p'' = p-3$. 
Then we consider the linear groups $K = K(p',p'') = K(3,r) 
\subset G(V) \subset G = G(\Pi ) \subset {\rm Aut} ({\fr r})$ 
introduced in section \ref{homogSec}. As before $\fr r$ denotes the
radical of ${\fr g} = {\rm Lie} \, G$.  
\begin{thm}\label{mainThm}
\begin{enumerate}
\item[1)] There exists a $G$-invariant quaternionic structure $Q$ on
$M = M(\Pi ) = G(\Pi )/K$.
\item[2)] If $\Pi$ is nondegenerate (see Def.~\ref{EPADef}) then
there exists a $G$-invariant pseudo-Riemannian metric $g$ on $M$
such that $(M,Q,g)$ is a  quaternionic pseudo-K\"ahler manifold.
\end{enumerate}
\end{thm}

\noindent
{\bf Proof:} The main idea of the proof, which we will carry out
in detail, is first to observe that the submanifold $M(V) = G(V)/K
\subset G/K = M$ has a natural $G(V)$-invariant structure of
quaternionic  pseudo-K\"ahler manifold and  then to study the problem of
extending this structure to the manifold $M$. As shown in section
\ref{homogSec}, we can embed $M(V)$ as open $G(V)$-orbit
into the pseudo-Riemannian symmetric space $\tilde{M}(V) = 
{\rm SO}_0(\tilde{V})/\tilde{K}$, $\tilde{V} = 
\mbox{\R}^{p+1,q+1}$, $\tilde{K} = {\rm SO}(p'+1)\times {\rm SO}_0(p'',q+1)
={\rm SO}(4) \times {\rm SO}_0(r,q+1)$. Now we claim that 
$\tilde{M}(V)$ carries an ${\rm SO}_0(\tilde{V})$-invariant
(and hence $G(V)$-invariant) almost quaternionic structure $Q$. 
It is sufficient to specify the corresponding $\tilde{K}$-invariant
quaternionic structure $Q_{[e]} \subset {\rm End}(T_{[e]}\tilde{M}(V))$
at the canonical base point $[e] = e\tilde{K} \in \tilde{M}(V)$.
We first define a hypercomplex structure $(J_{\alpha})$ on the
isotropy module $T_{[e]}\tilde{M}(V) \cong \mbox{\R}^4\otimes 
\mbox{\R}^{r,q+1}$: 
\begin{equation}
\label{JalphaEqu} J_1 := \mu_l(i) \otimes {\rm Id}\, ,\quad
J_2 := \mu_l(j) \otimes {\rm Id}\, ,\quad
J_3 := \mu_l(k) \otimes {\rm Id}\, .
\end{equation}
Here $\mu_l(x) \in{\rm End}(\mbox{\R}^4)$ stands for left-multiplication
by the quaternion $x \in \mbox{\Ha} = {\rm span} \{ 1,i,j,k\}$:
$\mu_l(x)y = xy$ for all $x,y \in  \mbox{\Ha} =  \mbox{\R}^4$. 
Right-multiplication by $x$ will be denoted by $\mu_r(x)$. 
Now it is easy to check that the quaternionic structure 
$Q_{[e]} := {\rm span}\{ J_1,J_2,J_3\}$ is normalized under the
isotropy representation of $\tilde{K}$ and hence extends to an  
${\rm SO}_0(\tilde{V})$-invariant almost quaternionic structure $Q$ on
$\tilde{M}(V)$. Note that the complex structures $J_{\alpha}$ are
not invariant under the isotropy representation and hence do not
extend to ${\rm SO}_0(\tilde{V})$-invariant almost complex structures
on $\tilde{M}(V)$. Next we observe that the tensor product of the 
ca\-no\-ni\-cal 
pseudo-Euclidean scalar products on $\mbox{\R}^4$ and $\mbox{\R}^{r,q+1}$
defines a pseudo-Euclidean scalar product $-g_{[e]}$ on $T_{[e]}\tilde{M}(V)
\cong \mbox{\R}^4 \otimes \mbox{\R}^{r,q+1}$. Notice that $g_{[e]}$ is 
{\em positive} definite if $r = 0$ and indefinite otherwise. It is invariant
under the isotropy representation and hence extends to a ($Q$-Hermitian)
${\rm SO}_0(\tilde{V})$-invariant pseudo-Riemannian metric on 
$\tilde{M}(V)$, which is unique, up to scaling, and defines on
$\tilde{M}(V)$ the well known structure of pseudo-Riemannian symmetric
space. In fact, it is the symmetric space associated to the  
symmetric pair $({\fr o}(\tilde{V}),\tilde{{\fr k}} := 
{\rm Lie}\, \tilde{K})$. 
Let $\tilde{{\fr m}} \subset {\fr o}(\tilde{V})$ be the 
$\tilde{{\fr k}}$-invariant complement to the isotropy algebra
$\tilde{{\fr k}}$. Then $[\tilde{{\fr m}},\tilde{{\fr m}}] \subset
\tilde{{\fr k}}$ and hence the Nomizu operators $L_x^g \in 
{\rm End} (T_{[e]}\tilde{M}(V))$ associated to the Levi-Civita 
connection $\nabla^g$ of $g$ vanish for $x \in \tilde{{\fr m}}$, see 
(\ref{KoszulEqu}). On the other hand if $x\in \tilde{{\fr k}}$ then
the Nomizu operator $L^g_x$ on $T_{[e]}\tilde{M}(V) \cong \tilde{{\fr m}}$
coincides with the image of $x$ under the isotropy representation:
$L_x^g = {\rm ad}_x|\tilde{{\fr m}}$.  This shows that $L_x^g$
normalizes the quaternionic structure $Q_{[e]}$ for all $x\in 
{\fr o}(\tilde{V})$ and by Cor.~\ref{qKCor} we conclude that
$(\tilde{M}(V),Q,g)$ is a homogeneous quaternionic  pseudo-K\"ahler 
manifold if $\dim \tilde{M}(V) >4$. The manifold
$\tilde{M}(V)$ is 4-dimensional only if $r=q=0$ and in this case
$\tilde{M}(V) = {\rm SO}_0(4,1)/{\rm SO}(4)$ reduces to 
real hyperbolic 4-space $H_{\mbox{\R}}^4$, i.e.\  to 
the quaternionic hyperbolic line   $H_{\mbox{\Ha}}^1 ={\rm Sp}(1,1)/
{\rm Sp}(1)\cdot {\rm Sp}(1)$, which is a standard example of 
(conformally half-flat Einstein) quaternionic K\"ahler 4-manifold. 
Of course, the pair $(Q,g)$ defined on the manifold
$\tilde{M}(V)$ (for all $r$ and $q$) restricts
to a $G(V)$-invariant  quaternionic pseudo-K\"ahler structure
on the open $G(V)$-orbit $M(V) \hookrightarrow \tilde{M}(V)$. 
So we have proven that $(M(V) =G(V)/K,Q,g)$ is a homogeneous
quaternionic pseudo-K\"ahler manifold. 

Our strategy is now to extend the geometric structures $Q$ and 
$g$ from $M(V)$ to $G$-invariant structures on $M = G/K \supset
G(V)/K = M(V)$. 
First we will extend $Q$ to a $G$-invariant almost quaternionic 
structure on $M$. Using the infinitesimal action of ${\fr g} =
{\fr g}(\Pi ) = {\fr g}(V) + W$ on $M$ we can identify
$T_{[e]}M = ({\fr g}(V)/{\fr k}) \oplus W = T_{[e]}M(V) \oplus W$.
Note that the isotropy representation of $K$ on $T_{[e]}M$ preserves
this decomposition and acts on $T_{[e]}M(V) = T_{eK}M(V) \cong 
T_{e\tilde{K}}\tilde{M}(V)$ as restriction of the isotropy representation
of $\tilde{K}$ on $T_{e\tilde{K}}\tilde{M}(V)$ to the subgroup
$K \subset \tilde{K}$. We extend the hypercomplex structure
$(J_{\alpha})$ defined above on $T_{e\tilde{K}}\tilde{M}(V)
\cong T_{eK}M(V)$ to a hypercomplex structure on $T_{eK}M =
T_{eK}M(V) \oplus W$ as follows:
\begin{equation}\label{extEqu}  
J_{\alpha}s := e_{\beta}e_{\gamma}s\, ,\quad s\in W\, ,
\end{equation}
where $(e_{\alpha},e_{\beta},e_{\gamma})$ is a cyclic permutation of 
the standard orthonormal basis of $E = \mbox{R}^{3,0} \subset  V = E + E'$
(the product $e_{\beta}e_{\gamma}$ is in the Clifford algebra). 
Now we can extend the quaternionic structure $Q_{[e]}$ on 
$T_{[e]}M(V)$ to a quaternionic structure on $T_{[e]}M$ such that
$Q_{[e]} = {\rm span}\{ J_1,J_2,J_3\}$. To prove that 
$Q_{[e]} \subset {\rm End}(T_{[e]}M)$ extends to a $G$-invariant
almost quaternionic structure on $M$ it is sufficient to
check that the isotropy algebra $\fr k$ normalizes $Q_{[e]}$. 
It is obvious that the subalgebra ${\fr o}(E') = {\fr o}(r,q) 
\subset {\fr k} = {\fr o}(E) \oplus {\fr o}(E')$
centralizes $Q_{[e]}$. We have to show that ${\fr o}(E) = {\fr o}(3)$
normalizes $Q_{[e]}$. In terms of the standard basis $(e_1,e_2,e_3) 
= (i,j,k)$ of $E = \mbox{\R}^3 = {\rm Im} \mbox{\Ha}$ a basis
of ${\fr o}(3) = \wedge^2E$ is given by $(e_2\wedge e_3, e_3\wedge
e_1, e_1\wedge e_2)$. For any cyclic permutation $(\alpha , \beta ,
\gamma )$ of $(1,2,3)$ we have the following easy formulas, cf.\ 
(\ref{adEqu}):
\begin{equation} -2 d\rho (e_{\beta}\wedge e_{\gamma}) = (\mu_l(e_{\alpha}) -
\mu_r(e_{\alpha})) \otimes {\rm Id} \label{dRHOEqu}\end{equation}
on $T_{eK}M(V) \cong T_{e\tilde{K}}\tilde{M}(V) \cong \mbox{\R}^4 \otimes 
\mbox{\R}^{r,q+1}$ and 
\begin{equation} \label{JaEqu} -2 d\rho (e_{\beta}\wedge e_{\gamma}) = 
e_{\beta}e_{\gamma} = J_{\alpha}
\end{equation}  
on $W$. As before $d\rho : {\fr k} \rightarrow {\fr gl}(T_{[e]}M)$ denotes
the isotropy representation. 
{}From the equations (\ref{JalphaEqu})-(\ref{JaEqu}) we immediately 
obtain that
\[ [d\rho (e_{\beta}\wedge e_{\gamma}),J_{\alpha}] = 0 \quad
\mbox{and} \quad  [d\rho (e_{\beta}\wedge e_{\gamma}),J_{\beta}] = 
-J_{\gamma}\, .\]
This shows that $Q_{[e]} \subset {\rm End}(T_{[e]}M)$ is invariant
under $\fr k$ and hence extends to a $G$-invariant almost quaternionic 
structure $Q$ on $M$. The 1-integrability of $Q$ will be proven 
in the sequel.

Let us first treat the case where $\Pi$ is nondegenerate. First of all, 
we extend the $G(V)$-invariant $Q$-Hermitian pseudo-Riemannian
metric $g$ on $M(V)$ to a $G$-invariant $Q$-Hermitian pseudo-Riemannian
metric on $M$. It is sufficient to extend the $K$-invariant pseudo-Euclidean 
scalar product $g_{[e]}$ on $T_{[e]}M(V)$ to a 
$K$-invariant and $J_{\alpha}$-invariant ($\alpha = 1,2,3$) scalar product on
$T_{[e]} = T_{[e]}M(V) \oplus W$. We do this in such a way that the
above decomposition is orthogonal for the extended scalar product $g_{[e]}$ on 
$T_{[e]}M$ and define
\[ g_{[e]}|W\times W := -b\, ,\]
where $b = b_{\Pi,(e_1,e_2,e_3)}$ is the canonical symmetric
bilinear form associated to $\Pi$ and to the decomposition $p = p' + p''
=3+r$. By Prop.~\ref{kerProp} the nondegeneracy of $\Pi$ implies
the nondegeneracy of $b$. This shows that the symmetric bilinear
form $g_{[e]}$ on $T_{[e]}M$ defined above is indeed a 
pseudo-Euclidean scalar product. This scalar product is invariant under
the isotropy group $K = K(p',p'') = K(3,r)$, by virtue of 
Thm.~\ref{bThm}, 3), and hence extends to a $G$-invariant
pseudo-Riemannian metric $g$ on $M$. Moreover, the metric
$g$ is $Q$-Hermitian. In fact, it is sufficient to observe that $b$
is $J_{\alpha}$-invariant. This property follows from the 
$K$-invariance of $b$, since $J_{\alpha}|W \in d\rho ({\fr k})|W$, see
(\ref{JaEqu}). Summarizing, we obtain: $(M,Q,g)$ is a homogeneous
pseudo-Hermitian almost quaternionic  manifold. The next step is to
compute the Nomizu map $L^g: {\fr g} \rightarrow {\rm End} (T_{[e]}M)$
associated to the Levi-Civita connection $\nabla^g$ of $g$. 
It is convenient to identify $T_{[e]}M$ with a $\fr k$-invariant
complement $\fr m$ to $\fr k$ in $\fr g$. Such complement is easily
described in terms of the $K$-invariant orthogonal decomposition
$V = E + E' = \mbox{\R}^{3,0} + \mbox{\R}^{r,q}$. Indeed, using the 
canonical identification ${\fr o}(V) = \wedge^2 V$ we have the 
following $\fr k$-invariant decompositions:
\[ {\fr o}(V) = {\fr k} + E\wedge E'\, ,\]
\[ {\fr g}(V) = {\fr o}(V) + \mbox{\R}D + V = {\fr k} + {\fr m}(V)\, ,
\quad {\fr m}(V) = E\wedge E' + \mbox{\R}D + V\, ,\]
\[ {\fr g} = {\fr g}(V) +W = {\fr k} + {\fr m}\, ,\quad 
{\fr m} = {\fr m}(\Pi ) := {\fr m}(V) + W\, .\] 
In the following, we will make the identifications 
$T_{[e]}M(V) = {\fr g}(V)/{\fr k} = {\fr m}(V) \subset
T_{[e]}M = {\fr g}/{\fr k} = {\fr m}$\label{mPage}. Let 
$(e_a')$, $a = 1, \ldots , r+q$, be the standard orthonormal basis of 
$E' = \mbox{\R}^{r,q}$. Then we have the following 
orthonormal basis of $({\fr m}(V),g_{[e]})$:
\[ (e_{\alpha}\wedge e_a',D,e_{\alpha},e_a')\, ,\]
where $\alpha = 1,2,3$ and $a = 1, \ldots , r+q$. For notational convenience
we denote this basis by \label{eiotaiPage} 
$(e_{\iota i})$, $\iota = 0, \ldots 3$, 
$i= 0,\ldots ,r+q$, where 
\begin{eqnarray*} 
e_{00} &:=& e_0' := e_0 := D\\
e_{0a} &:=& e_a' \quad (a = 1,\ldots , r+q)\\
e_{\alpha 0} &:=& e_{\alpha} \quad (\alpha =1,2,3)\\
e_{\alpha a} &:=& e_{\alpha}\wedge e_a'\, .
\end{eqnarray*} 
In terms of this basis the hypercomplex structure $J_{\alpha}$ is expressed
on ${\fr m}(V)$ simply by:
\[ J_{\alpha}e_{0i} = J_{\alpha}e_i' = e_{\alpha i}\quad \mbox{and} 
\quad J_{\alpha}e_{\beta i} = e_{\gamma i}\, , 
\] 
where $i = 0,\ldots ,r+q$ and $({\alpha},{\beta},{\gamma})$ is a 
cyclic permutation of $(1,2,3)$. The scalar pro\-duct $g_{[e]}$
is completely determined on  ${\fr m}(V)$ by the condition that 
$(e_{\iota i})$ is orthonormal and that 
\[ \epsilon_i := g(e_{\iota i},e_{\iota i}) = \left\{ 
\begin{array}{r@{\quad \mbox{if} \quad}l}
-1 & i = 1,\ldots , r\\
+1 & i = 0, r+1,\ldots ,r+q\, .
\end{array}
\right.
\]
The scalar product $g_{[e]}$ on $\fr m$ induces an identification
$x\wedge y \mapsto x\wedge_gy$ of the exterior square 
$\wedge^2{\fr m}$ with the vector space
of $g_{[e]}$-skew symmetric endomorphisms on $\fr m$, 
where $x\wedge_gy(z) := g(y,z)x -g(x,z)y$, $x,y,z\in {\fr m}$. 
We denote by ${\rm n}(Q)$ (respectively, ${\rm z}(Q)$) the normalizer
(respectively, centralizer) of the quaternionic structure $Q_{[e]}$
in the Lie algebra ${\fr gl}({\fr m})$. 
\begin{lemma} \label{LCLemma}
The Nomizu map $L^g = L(\nabla^g)$ associated to the Levi-Civita connection 
$\nabla^g$ of the homogeneous pseudo-Riemannian manifold $(M = G/K,g)$
is given by the following formulas:
\begin{eqnarray*}
L^g_{e_0} &=&  L^g_{e_{\alpha a}} = 0\, ,\\
L^g_{e_{\alpha}} &=& \frac{1}{2}J_{\alpha} + \bar{L}^g_{e_{\alpha}}\, ,
\end{eqnarray*}
where 
\[ \bar{L}^g_{e_{\alpha}} = \frac{1}{2}\sum_{i=0}^{r+q}
\epsilon_i e_{\alpha  i} \wedge_ge_i' - \frac{1}{2}\sum_{i=0}^{r+q}
\epsilon_i e_{\gamma i}  \wedge_ge_{\beta i} \in {\rm z}(Q)\, ,\]
$L^g_{e_a'} \in {\rm z}(Q)$ is given by:
\begin{eqnarray*}
L^g_{e_a'} |{\fr m}(V) &=& \sum_{\iota =0}^3 e_{\iota a} \wedge_g e_{\iota}
\, ,\\
L^g_{e_a'} |W &=& \frac{1}{2} e_1e_2e_3e_a'\, .
\end{eqnarray*}
{}For all $s\in W$ the Nomizu operator $L_s^g\in {\rm z} (Q)$ maps the
subspace ${\fr m}(V) \subset {\fr m} = {\fr m}(V) + W$ into $W$ and
$W$ into ${\fr m}(V)$. The restriction $L_s^g|W$ ($s\in W$) is completely 
determined by $L^g_s|{\fr m}(V)$ (and vice versa) according to the
relation
\[ g(L_s^gt,x) = -g(t,L_s^gx)\, ,\quad s,t\in W\, ,\quad x\in {\fr m}(V)\, .
\]
Finally,  $L^g_s|{\fr m}(V)$ ($s\in W$) is completely 
determined by its values on the quaternionic basis $(e_i')$, 
$i = 0,\ldots , r+q$, which are as follows:
\begin{eqnarray*}
L_s^ge_0' &=& \frac{1}{2}s\, ,\\
L_s^ge_a' &=& \frac{1}{2}e_1e_2e_3e_a's\, .
\end{eqnarray*}
(It is understood that $L_x^g = d\rho (x) = {\rm ad}_x|{\fr m}$ for all
$x\in {\fr k}$, cf.\ equation (\ref{NomizuIEqu}).)
In the above formulas $a = 1, \ldots , r+q$, $\alpha = 1,2,3$ and 
$({\alpha},{\beta},{\gamma})$ is a 
cyclic permutation of $(1,2,3)$. 
\end{lemma}

\noindent
{\bf Proof:} This follows from equation (\ref{KoszulEqu}) by a straightforward
computation. $\Box$ 

\begin{cor}
\label{LCCor} The Levi-Civita connection $\nabla^g$ of the homogeneous
almost quaternionic pseudo-Hermitian manifold $(M,Q,g)$ preserves
$Q$ and hence $(M,Q,g)$ is a quaternionic  pseudo-K\"ahler manifold.
\end{cor}

\noindent
{\bf Proof:} From the $G$-invariance of $Q$, we know already that 
$L_x^g \in {\rm n}(Q)$ for all $x\in {\fr k}$ and the formulas of 
Lemma \ref{LCLemma} show that $L_x^g \in {\rm n}(Q)$ for all $x\in {\fr m}$.
Now the corollary follows from Cor.~\ref{qKCor}. $\Box$ 

By Cor.~\ref{LCCor} we have already established part 2) of 
Thm.~\ref{mainThm}. Part 1) is a consequence of 2) provided that
$\Pi$ is nondegenerate. It remains to discuss the case of degenerate
$\Pi$. 

By Thm.~\ref{decompThm}, we have a direct decomposition $W = W_0 + W'$ of  
$C\! \ell^0(V)$-modules, where $W_0 = {\rm ker}\, \Pi$. We put 
$\Pi':= \Pi|\wedge^2W'$ and denote by $L': {\fr g}(\Pi') \rightarrow
{\rm End}({\fr m}(\Pi'))$ the Nomizu map associated to the Levi-Civita
connection of the  quaternionic  pseudo-K\"ahler manifold $(M':= M(\Pi'),
Q,g)$. Note that the quaternionic structure $Q$ of $M'$ coincides with the
almost quaternionic structure induced by the obvious $G(\Pi')$-equivariant
embedding $M' = G(\Pi')/K \subset G(\Pi)/K =M$. By the next lemma, 
we can extend the map $L'$ to a torsionfree Nomizu map 
$L:{\fr g}(\Pi) \rightarrow {\rm End}({\fr m}(\Pi))$, whose
image normalizes $Q_{[e]}$. This proves the 1-integrability of $Q$
(by Cor.~\ref{quatCor}), completing the proof of
Thm.~\ref{mainThm}. $\Box$  

By Cor.\ \ref{LCCor} we can decompose $L_x' = \sum_{\alpha=1}^3
\omega_{\alpha}'(x)J_{\alpha} + \bar{L}_x'$, where 
the $\omega_{\alpha}'$ are 1-forms on ${\fr m}(\Pi')$ and 
$\bar{L}_x' \in {\rm z}(Q)$ belongs to the centralizer of the 
quaternionic structure on ${\fr m}(\Pi')$. 

\begin{lemma} \label{NomizuLemma} 
The Nomizu map $L':{\fr g}(\Pi') \rightarrow {\rm End}({\fr m}(\Pi'))$ 
associated to the Levi-Civita connection of $(M(\Pi') = G(\Pi')/K,g)$
can be extended to the Nomizu map $L:{\fr g}(\Pi) \rightarrow 
{\rm End}({\fr m}(\Pi))$ of a  $G(\Pi)$-invariant 
quaternionic connection $\nabla$ on the homogeneous almost 
quaternionic manifold
$(M(\Pi) = G(\Pi)/K,Q)$. The extension is defined as follows:
\[ L_x:= \sum_{\alpha=1}^3 \omega_{\alpha}(x)J_{\alpha} + \bar{L}_x\, ,
\quad x\in {\fr m}(\Pi)\, ,\]
where $\bar{L}_x \in{\rm z}(Q)$, with centralizer  taken 
in ${\fr gl}({\fr m}(\Pi))$, is defined below and 
the 1-forms $\omega_{\alpha}$ on ${\fr m}(\Pi) := {\fr m}(\Pi') + W_0$
satisfy $\omega_{\alpha}|{\fr m}(\Pi') := \omega_{\alpha}'$ and 
$\omega_{\alpha}|W_0 := 0$. The operators $\bar{L}_x$ are given by:
\[ \bar{L}_x|{\fr m}(\Pi') := \bar{L}_x' \quad \mbox{if} \quad
x\in {\fr m}(\Pi')\, ,\]
\[ \bar{L}_{e_0}|W_0 := \bar{L}_{e_{\alpha a}}|W_0 
:= \bar{L}_{e_{\alpha}}|W_0 := 0\, ,
\]
\[ \bar{L}_{e_a'}|W_0 := \frac{1}{2} e_1e_2e_3e_a'\, ,\]
\[ \bar{L}_s|W_0 := 0 \quad \mbox{if} \quad s\in W\, ,\]
\[ \bar{L}_sx := L_xs -[s,x] \quad \mbox{if} \quad s\in W_0\, ,\quad
x\in {\fr m}(\Pi')\, .\]
(It is understood that $L_x = {\rm ad}_x|{\fr m}(\Pi)$ for all $x\in{\fr k}$.
In the above formulas, as usual, $a = 1, \ldots , r+q$ and $\alpha = 1,2,3$.)
\end{lemma}

\noindent
{\bf Proof:} To check that $L: {\fr g}(\Pi) \rightarrow 
{\rm End}({\fr m}(\Pi))$ is a Nomizu map it is sufficient to check 
(\ref{NomizuIIEqu});  (\ref{NomizuIEqu}) is satisfied by definition.
The equation (\ref{NomizuIIEqu}) expresses the $K$-invariance
of $L\in {\fr g}(\Pi)^{\ast}\otimes {\rm End}({\fr m}(\Pi))$ which
is equivalent to the invariance of $L$ under the Lie algebra
$\fr k$, since $K$ is a connected Lie group. So (\ref{NomizuIIEqu})
is equivalent to the following equation:
\begin{equation}
\label{NomizuIII} L_{[x,y]}z = [{\rm ad}_x,L_y]z\quad x\in {\fr k}\, ,
\quad y,z\in {\fr m}(\Pi)\, .
\end{equation}
We first check this equation. Let always $x\in {\fr k}$. Due to $L_{W}
W_0 = 0$ and $[{\fr k},W_0] \subset W_0$ we have 
\[ L_{[x,y]}z  = [{\rm ad}_x,L_y]z = 0\]
if $y\in W$ and $z\in W_0$. Also (\ref{NomizuIII}) is satisfied if $y,z\in 
{\fr m}(\Pi')$ because $L_x|{\fr m}(\Pi') = L_x'$ is the 
Nomizu operator associated to the Nomizu map $L'$.
Now let $y\in W_0$ and $z\in {\fr m}(\Pi')$. Then we compute:
\[ L_{[x,y]}z - [{\rm ad}_x,L_y]z = (L_z[x,y] - 
[[x,y],z]) - [x,L_yz] + L_y[x,z] \]
\[ = L_z[x,y] - [[x,y],z] -[x,L_zy] +[x,[y,z]]
+ L_{[x,z]}y - [y,[x,z]]\]
\[ = L_{[x,z]}y -[x,L_zy] + L_z[x,y] = L_{[x,z]}y -[{\rm ad}_x,L_z]y\,.\]
This shows that it is sufficient to check  (\ref{NomizuIII}) for 
$y\in {\fr m}(\Pi')$ and $z\in W_0$; the case $y\in W_0$ and 
$z\in  {\fr m}(\Pi')$ then follows from the above computation. 
In the following let $x\in {\fr k}$ and $z \in W_0$. We check
(\ref{NomizuIII}) for all $y \in {\fr m}(\Pi')$.  From $L_{e_0} = 0$
and $[{\fr k},e_0] = 0$ it follows that 
\[ L_{[x,e_0]}z - [x,L_{e_0}z] +  L_{e_0}[x,z] = 0\, .\]
Next we check  (\ref{NomizuIII}) for $y = e_{\alpha}$: 
\[ L_{[x,e_{\alpha}]}z - [x,L_{e_{\alpha}}z] +  L_{e_{\alpha}}[x,z] = 
L_{[x,e_{\alpha}]}z - \frac{1}{2} [x,J_{\alpha}z] +\frac{1}{2}J_{\alpha}
[x,z]\, .\]
It is clear that the first summand and the sum of the second
and third summands vanish if $x\in {\fr o}(E') \subset 
{\fr k} = {\fr o}(E) \oplus {\fr o}(E')$. For 
$x = e_{\alpha}\wedge e_{\beta} \in {\fr o}(E) = {\fr o}(3)$ 
($(\alpha , \beta , \gamma )$ cyclic) we compute:
\[ L_{[e_{\alpha}\wedge e_{\beta},e_{\alpha}]}z 
- \frac{1}{2} [e_{\alpha}\wedge e_{\beta},J_{\alpha}z] +\frac{1}{2}J_{\alpha}
[e_{\alpha}\wedge e_{\beta},z] \]
\[ = -L_{e_{\beta}}z + \frac{1}{4}e_{\alpha}e_{\beta}e_{\beta}e_{\gamma}z
- \frac{1}{4}e_{\beta}e_{\gamma}e_{\alpha}e_{\beta}z\]
\[ = -\frac{1}{2} e_{\gamma}e_{\alpha}z + \frac{1}{2} e_{\gamma}e_{\alpha}z = 0
\]
and for $x = e_{\beta} \wedge e_{\gamma}$ we obtain:
 \[ L_{[e_{\beta}\wedge e_{\gamma},e_{\alpha}]}z 
-  \frac{1}{2} [e_{\beta}\wedge e_{\gamma},J_{\alpha}z]
+ \frac{1}{2}J_{\alpha}[e_{\beta}\wedge e_{\gamma},z] \]
\[ = 0 +\frac{1}{4} e_{\beta} e_{\gamma}e_{\beta} e_{\gamma}z -
\frac{1}{4} e_{\beta} e_{\gamma}e_{\beta} e_{\gamma}z = 0\, .\]
Next we check  (\ref{NomizuIII}) for $y = e_a'$: 
\[ L_{[x,e_a']}z -[x,L_{e_a'}z] + L_{e_a'}[x,z] =\]
\[ \frac{1}{2}e_1e_2e_3[x,e_a']z - \frac{1}{2}[x,e_1e_2e_3e_a'z] 
+\frac{1}{2}e_1e_2e_3e_a'[x,z]\, .\]
It is easy to see that this is zero if $[x,e_a'] = 0$. So we can put 
$x = e_b'\wedge e_a'$ obtaining:
\[ \frac{1}{2}\langle e_a',e_a'\rangle e_1e_2e_3e_b'z 
+\frac{1}{4}e_b'e_a'e_1e_2e_3e_a'z
- \frac{1}{4}e_1e_2e_3e_a'e_b'e_a'z = 0 \, .\] 
Finally, for $y \in E\wedge E'$ we immediately obtain:
\[ L_{[x,y]}z -[x,L_yz] + L_y[x,z] = L_{[x,y]} z = 0\, ,\]
since $L_y = 0$ for $y \in E\wedge E'$ and $[{\fr k}, E\wedge E']
\subset E\wedge E'$. So we have proven that $L$ is a Nomizu map. It
is easily checked that $L_xy -L_yx = -\pi [x,y]$ for all
$x,y\in {\fr m}(\Pi )$, where $\pi : {\fr g}(\Pi ) \rightarrow
{\fr m}(\Pi) \cong T_{[e]}M(\Pi)$ is the projection along $\fr k$. 
This shows that $\nabla = \nabla (L)$ has zero torsion. 
It only remains to check that $L_x \in {\rm n}(Q)$ for all
$x\in {\fr g}(\Pi )$. This is easy to see for $x\in 
{\fr g}(\Pi' )$ from the definition of the map $L$ as extension
of the Nomizu map $L'$. We present the 
calculation only for  $L_s$, $s\in W_0$:
\[ L_se_0 = L_{e_0}s - [s,e_0] = 0 +\frac{1}{2}s = \frac{1}{2}s\, ,\]
\[ L_sJ_{\alpha}e_0 = L_se_{\alpha} = L_{e_{\alpha}}s - 
[s,e_{\alpha}] = \frac{1}{2}J_{\alpha}s - 0 = J_{\alpha} L_se_0\, ,\]
\[ L_se_a' = L_{e_a'}s - [s,e_a'] = 
\frac{1}{2}e_1e_2e_3e_a's - 0 = \frac{1}{2}e_1e_2e_3e_a's \, ,\]
\[ L_sJ_{\alpha}e_a' = L_s e_{\alpha a} = L_{e_{\alpha a}} s -
[s,e_{\alpha a}] = 0 - \frac{1}{2}e_{\alpha}e_a's\]
\[ = e_{\beta}e_{\gamma}(\frac{1}{2}e_{\alpha}e_{\beta}e_{\gamma}
e_a's) = J_{\alpha}L_se_a'\, ,\]
$(\alpha , \beta ,\gamma )$ cyclic, 
\[ L_st = L_ts - [s,t]  = 0\, ,\]
\[ L_sJ_{\alpha}t = L_{J_{\alpha}t}s - [s, J_{\alpha}t] = 0 = J_{\alpha}L_ts\]
if $t\in W$. We have used that $L_WW_0 = 0$ and $[W_0,W] = 0$. 
This shows that $L_s\in {\rm z}(Q)$
for all $s\in W_0$ finishing the proof of the lemma. $\Box$

\subsection{The Riemannian case} 
\label{RCSec} 
\begin{prop}
\label{RiemProp} Let ${\fr p}(\Pi ) = {\fr p}(V) + W$ be a nondegenerate
extended Poincar\'e algebra of signature $(p,q)$, $p\ge 3$, and 
$(M(\Pi ),Q,g)$ the corresponding homogeneous  quaternionic pseudo-K\"ahler 
manifold, see Thm.\ \ref{mainThm}. Then the pseudo-Riemannian
metric $g$ is positive definite, and hence a Riemannian metric, if 
and only if $-b$ is positive definite and $p=3$. In all other cases
$g$ is indefinite. Here $b = b_{\Pi , (e_1,e_2,e_3)}$ is the
canonical symmetric bilinear form associated to $\Pi$ and to the decomposition
$p = 3 + r$. 
\end{prop}

\noindent
{\bf Proof:} By construction, the restriction of $g$ to the submanifold
$M(V) \subset M(\Pi )$ is a (positive definite) Riemannian metric if
and only if $p=3$ and is indefinite otherwise. Now the proposition
follows from the fact that $-b$ is precisely the restriction of the
scalar product $g_{[e]}$ on $T_{[e]}M(\Pi ) \cong T_{[e]}M(V) \oplus W$
to the subspace $W$. $\Box$ 

Next we will use the classification of extended Poincar\'e algebras of
signature $(3,q)$ up to isomorphism (see \ref{isomSec}) to derive 
the classification of the quaternionic K\"ahler manifolds 
$(M(\Pi ),Q,g)$ up to isometry. We recall that ${\fr p}(3,q,0,0,l) =
{\fr p}(p=3,q,l_0=0,l_+=0,l_-=l)$ (respectively, ${\fr p}(3,q,0,0,l^+,0,0,l^-)
= {\fr p}(p=3,q,l_0^+=0,l_+^+=l^+,l_0^-=0,l_+^-=0,l_-^-=l^-)$) is the
set of isomorphism classes of extended Poincar\'e algebras for which $-b$ is
positive definite if $q\not\equiv 3\pmod{4}$ (respectively, if 
$q\equiv 3\pmod{4}$). We denote by $M(q,l)$ (respectively, $M(q,l^+,l^-)$)
the homogeneous quaternionic K\"ahler manifold $(M(\Pi ),Q,g)$ 
associated to $\Pi \in {\fr p}(3,q,0,0,l)$  (respectively, 
$\Pi \in {\fr p}(3,q,0,0,l^+,0,0,l^-)$). 

\begin{thm}\label{PDThm}
Every homogeneous quaternionic  pseudo-K\"ahler manifold of the form 
$(M(\Pi ),Q,g)$ for which $g$ is positive definite is isometric to 
one of the homogeneous quaternionic K\"ahler manifolds $M(q,l)$ 
($q\not\equiv 3\pmod{4})$ or $M(q,l^+,l^-) \cong M(q,l^-,l^+)$
($q\equiv 3\pmod{4})$. In particular, there are only  countably
many such Riemannian manifolds up to isometry. 
\end{thm}

\noindent
{\bf Proof:} This is a direct consequence of Prop.\ \ref{RiemProp}, 
Thm.\ \ref{moduliThmI} and Thm.\ \ref{moduliThmII}. $\Box$ 

Any real vector space $E$ admitting a quaternionic structure has 
$\dim E \equiv 0\pmod{4}$. Therefore, we can define its 
{\bf quaternionic  dimension} $\dim_{\mbox{\Ha}} E := \dim E/4$. 
Similarly, the quaternionic dimension of a quaternionic manifold 
$(M,Q)$ is $\dim_{\mbox{\Ha}} M := \dim M/4$. We denote by 
$N(q)$ the quaternionic dimension of an irreducible 
$C\! \ell^0_{3,q}$-module.  

\begin{prop}
\begin{enumerate}
\item[1)]  $N(0) = N(1)
= N(2) = N(3) = 1$, $N(4) = 2$, $N(5) = 4$, $N(6) = N(7) = 8$ and  
$N(q + 8) = 16N(q)$ for all $q\ge 0$. In particular,  $N(q)$ coincides 
with the dimension of an irreducible
$\mbox{\Z}_2$-graded $C\! \ell_{q-3}$-module if $q\ge 3$. 
\item[2)] The quaternionic dimension of the  homogeneous quaternionic 
K\"ahler manifolds $M(q,l)$ and $M(q,l^+,l^-)$ is given by:
\[ \dim_{\mbox{\Ha}} M(q,l) = q+1+lN(q)\]
and 
\[ \dim_{\mbox{\Ha}} M(q,l^+,l^-) = q+1+(l^++l^-)N(q)\, .\]
\end{enumerate}
\end{prop}

\noindent
{\bf Proof:} The first part follows from the classification of Clifford 
algebras. The second part follows from $\dim M = \dim M(V) + \dim W$. $\Box$

The next theorem identifies the spaces $M(q,l)$ and $M(q,l^+,l^-)$ with
Alekseevsky's quaternionic K\"ahler manifolds. We recall that an
{\bf Alekseevsky space} is a quaternionic K\"ahler manifold which admits
a simply transitive non Abelian splittable sol\-va\-ble group of isometries, 
see \cite{A3} and \cite{C2}. Due to Iwasawa's decomposition theorem 
any symmetric quaternionic K\"ahler manifold of noncompact type
is an Alekseevsky space. These are precisely the noncompact duals
of the Wolf spaces. We recall that a {\bf Wolf space} is a symmetric 
quaternionic K\"ahler manifold of compact type  and that such manifolds
are in 1-1-correspondence with the complex simple Lie algebras \cite{W1}. 
The nonsymmetric Alekseevsky spaces are grouped into 3 series: 
$\cal V$-spaces, $\cal W$-spaces and $\cal T$-spaces, see \cite{A3}, 
\cite{dW-VP2} and \cite{C2}. These 3 series 
contain also all symmetric  Alekseevsky spaces of rank $>2$ and no 
symmetric spaces of smaller rank. 
By definition an Alekseevsky space can be presented as metric 
Lie group, i.e.\ as homogeneous Riemannian manifold of the form
$(L,g)$, where $L$ is a Lie group and $g$ a left-invariant
Riemannian metric on $L$. 
\begin{thm}
\label{AThm} Let $(M = M(\Pi ) = G(\Pi )/K,Q,g)$ be a homogeneous
quaternionic K\"ahler manifold as in Prop.\ \ref{RiemProp}, 
\[ G(\Pi ) = S  \mbox{\Bbb n} R \cong {\rm Spin}_0(V) \mbox{\Bbb n} R\] 
the Levi decomposition (\ref{LeviIIEqu}) and $I(S)$ the 
Iwasawa subgroup of $S$. Then $L := L(\Pi ) :=  I(S) 
 \mbox{\Bbb n} R \subset G(\Pi )$ is a (non Abelian) splittable
solvable Lie subgroup which acts simply transitively on $M$. 
In particular, $(M,Q,g)$ is an Alekseevsky space. More precisely,
we have the following identifications with the 
$\cal V$-spaces, $\cal W$-spaces, $\cal T$-spaces and symmetric
Alekseevsky spaces: 
\begin{enumerate}
\item[1)] $M(q,l) = {\cal V}(l,q-3)$ and 
$M(q,l^+,l^-) = {\cal V}(l^+,l^-,q-3)$ if $q\ge 4$,
\item[2)] $M(3,l^+,l^-) = {\cal W}(l^+,l^-)$, 
\item[3)] $M(2,l) = {\cal T}(l)$,
\item[4)] $M(1,l) = {\rm SU}(l+2,2)/{\rm S}({\rm U}(l+2)\times {\rm U}(2))$,
\item[5)] $M(0,l) = {\rm Sp}(l+1,1)/{\rm Sp}(l+1) {\rm Sp}(1) = 
\mbox{\Ha} H^{l+1}$ (quaternionic hyperbolic $(l+1)$-space). 
\end{enumerate}
The above Riemannian manifolds $M(q,l)$ and $M(q,l^+,l^-) \cong M(q,l^-,l^+)$
are pairwise nonisometric and exhaust all Alekseevsky spaces with two
symmetric exceptions: $\mbox{\C}H^2 = {\rm SU}(1,2)/{\rm U}(2) =:  
M(1,-1)$ (complex hyperbolic plane) and ${\rm G}_2^{(2)}/{\rm SO}(4)$. 
(Note that these two symmetric spaces have rank $\le 2$ and so do not 
belong to any of the 3 series $\cal V$, $\cal W$ and $\cal T$ of
Alekseevsky spaces.) 
\end{thm}

\noindent
{\bf Proof:} The fact that $L$ acts simply transitively on $M$ follows
from Prop.\ \ref{stProp}. This shows that $(M,Q,g)$ is an Alekseevsky space.
The Riemannian metric $g$ induces a left-invariant metric $g_L$ on the 
Lie group $L$. To establish the identifications given in the theorem
it is sufficient to check that $(L,g_L)$ is isomorphic (as metric
Lie group) to one of the metric
Lie groups which occur in the classification of Alekseevsky spaces,
see \cite{A3} and  \cite{C2}. 
(The quaternionic structure can be reconstructed
from the holonomy of the Levi-Civita connection, 
up to an automorphism of the full isometry group which preserves the 
isotropy group.) Finally, to prove that $M(q,l)$ and 
$M(q,l^+,l^-) \cong M(q,l^-,l^+)$ are pairwise nonisometric it 
is, by \cite{A4}, sufficient to check that the corresponding metric Lie groups
(which occur in the classification of Alekseevsky spaces) are pairwise 
nonisomorphic. This was done in \cite{C2}. $\Box$

\subsection{A class of noncompact homogeneous  quaternionic Hermitian
manifolds with no transitive solvable group of isometries} 
\label{notransSec}
There is a widely known conjecture by D.V.\ Alekseevsky which says that any
noncompact homogeneous quaternionic K\"ahler manifold admits a
transitive solvable group of isometries \cite{A3}. The next theorem
shows that this conjecture becomes false if we replace ``K\"ahler''
by ``Hermitian''. 

\begin{thm}\label{NTRThm} 
Let ${\fr p}(\Pi )$ be any extended Poincar\'e algebra of signature
$(p,q) = (3+r,0)$, $r\ge 0$, and $(M = G/K,Q)$ the corresponding 
homogeneous quaternionic manifold, see Thm.\ \ref{mainThm}, 1). 
Then there exists a $G$-invariant $Q$-Hermitian Riemannian metric
$h$ on $M$. Moreover, the noncompact homogeneous quaternionic Hermitian
manifold $(M,Q,h)$ does not admit any transitive solvable Lie group
of isometries if $r>0$. 
\end{thm}

\noindent
{\bf Proof:} Since $K = K(3,r) \cong {\rm Spin}(3)\cdot {\rm Spin}(r)$ is
compact, one can easily construct a $K$-invariant $Q_{[e]}$-Hermitian
Euclidean scalar product $h_{[e]}$ on $T_{[e]}M$ by the standard
averaging procedure and extend it to a $G$-invariant $Q$-Hermitian
Riemannian metric $h$ on $M$. More explicitly, we can construct
such a scalar product $h_{[e]}$ on $T_{[e]}M \cong T_{[e]}M(V) \oplus W$
as orthogonal sum of $K$-invariant Euclidean scalar pro\-ducts on
$T_{[e]}M(V)$ and $W$ as follows. Using the open embedding
$M(V) \hookrightarrow \tilde{M}(V)$ we can identify 
$T_{eK}M(V) \cong T_{e\tilde{K}}\tilde{M}(V) \cong 
\mbox{\R}^4\otimes \mbox{\R}^{r,1} \cong \mbox{\R}^4 \otimes \mbox{\R}^{r+1}$
and choose the standard ${\rm O}(4) \times {\rm O}(r+1)$-invariant
Euclidean scalar product on $\mbox{\R}^4 \otimes \mbox{\R}^{r+1}$. This
scalar product is automatically $K$-invariant and $Q_{[e]}$-Hermitian.
On $W$ we choose any $K$-invariant Euclidean scalar product (which 
exists by compactness of $K$). It is automatically $Q_{[e]}$-Hermitian
because $Q_{[e]} \subset {\rm End} (W)$ is precisely the image of 
${\fr o}(3) \subset {\fr k} = {\rm Lie}\, K$ under the isotropy
representation of $\fr k$ on the $\fr k$-invariant subspace
$W \subset T_{[e]}M \cong T_{[e]}M(V) \oplus W$. It remains
to show that ${\rm Isom} (M,h)$ does not contain any transitive solvable
Lie subgroup if $r>0$. In fact, $M$ is homotopy equivalent to the
simply connected real Grassmannian ${\rm SO}(3+r)/{\rm SO}(3) \times
{\rm SO}(r)$ of oriented $3$-planes in $\mbox{\R}^{3+r}$ ($r>0$). 
On the other hand, if $(M,h)$ admits a transitive solvable group
of isometries then $M$ must be homotopy equivalent to a (possibly
trivial) torus, which  contradicts the fact that $M$ is simply connected
(and not contractible). $\Box$

\noindent  
The last argument proves, in fact, the following theorem.

\begin{thm}\label{NTThm} 
Let ${\fr p}(\Pi )$ be any extended Poincar\'e algebra of signature
$(p,q)$, $p>3$  and $M = M(\Pi )$ the manifold constructed in 
Thm.~\ref{mainThm}. Then $M$ does not admit any transitive
(topological) action by a solvable Lie group. 
\end{thm}
\section{Bundles associated to the quaternionic manifold $(M,Q)$}
\label{bundleSec} 
To any almost quaternionic manifold $(M,Q)$ one can canonically
associate the following bundles over $M$: the twistor bundle 
$Z(M)$, the canonical ${\rm SO}(3)$-principal bundle $S(M)$ and the  
Swann bundle $U(M)$. The {\bf twistor bundle} (or {\bf twistor space}) 
$Z(M) \rightarrow M$ is the subbundle of $Q$ whose fibre $Z(M)_m$ at $m\in M$
consists of all complex structures subordinate to the 
quaternionic structure $Q_m$, i.e.\ 
\[ Z(M)_m = \{ A\in Q_m| A^2 = -{\rm Id}\} \, .\] 
So $Z(M)$ is a bundle of $2$-spheres. The fibre $S(M)_m$ of the 
${\rm SO}(3)$-principal bundle $S(M)$ at $m\in M$
consists of all hypercomplex structures $(J_1,J_2,J_3)$ subordinate
to $Q_m$. Finally, 
\[ U(M) = S(M) \times_{{\rm SO}(3)} (\mbox{\Ha}^{\ast}/\{ \pm 1\}) \]
is associated to the action of ${\rm SO}(3) \cong {\rm Sp}(1)/\{ \pm 1\}$
on $\mbox{\Ha}^{\ast}/\{ \pm 1\}$ induced by left-multiplication of 
unit quaternions on $\mbox{\Ha}$. The total space $Z(M)$ carries a
canonical almost complex structure $\mbox{\J}$, which is integrable
if $Q$ is quaternionic, see \cite{A-M-P}. Similarly, one can define an
almost hypercomplex structure $(\mbox{\J}_1,\mbox{\J}_2,\mbox{\J}_3)$
on $U(M)$,  which is integrable if $Q$ is quaternionic, cf.\ \cite{P-P-S}. 
We recall the definition of the  complex structure $\mbox{\J}$
on the twistor space $Z = Z(M)$ of a quaternionic manifold $(M,Q)$.
Since $Q$ is $1$-integrable, there exists a quaternionic connection
$\nabla$ on $M$, see Def.\ \ref{connDef} and Def.\ \ref{integrDef}.  
The holonomy of $\nabla$ preserves not only $Q\subset {\rm End}(TM)$ 
but also its sphere subbundle $Z \subset Q$, simply because ${\rm Id}$ is a 
parallel section of ${\rm End}(TM)$.  Let
\begin{equation} \label{connEqu} 
TZ = TZ^{ver} \oplus TZ^{hor}
\end{equation}
be the corresponding decomposition into the vertical space $TZ^{ver}$
tangent to the fibres of the twistor bundle $Z\rightarrow M$ and its 
$\nabla$-horizontal complement $TZ^{hor}$.  The complex structure 
$\mbox{\J}$ preserves the decomposition
(\ref{connEqu}). Let $m\in M$ be a point in $M$ and 
$z:=J \in Z_m \subset Q_m$ a complex structure on $T_mM$ subordinate to $Q_m$. 
Then $\mbox{\J}_z \in {\rm End}(T_zZ)$ is defined by:
\[ \mbox{\J}A := JA \quad \mbox{and} \quad \mbox{\J}\tilde{X} = \widetilde{JX}
\]
for all $A \in T_zZ^{ver} = T_zZ_m = \{ A\in Q_m|AJ =-JA\}$ and all 
$X\in T_mM$, where  $\tilde{X} \in T_zZ^{hor}$ denotes the 
$\nabla$-horizontal lift of $X$. It was proven in \cite{A-M-P}
that $\mbox{\J}$ does not depend on the choice of quaternionic connection 
$\nabla$. 

If $(M,Q)$ admits a  quaternionic pseudo-K\"ahler metric $g$
(of nonzero scalar curvature) then it is known that $(Z(M),\mbox{\J})$
admits a complex contact structure and a pseudo-K\"ahler-Einstein
metric \cite{S1}, that $S(M)$ admits a pseudo-3-Sasakian structure
\cite{Ko} and that $(U(M),\mbox{\J}_1,\mbox{\J}_2,\mbox{\J}_3)$
admits a pseudo-hyper-K\"ahler metric \cite{Sw1}. Moreover, all these
special geometric structures are canonically associated to the data $(M,Q,g)$. 
We recall that a {\bf complex contact structure} on a complex 
manifold $Z$ is a holomorphic distribution 
$\cal D$ of codimension one whose Frobenius form 
$[\cdot , \cdot ]: \wedge^2 {\cal D} \rightarrow TZ/{\cal D}$ is 
(pointwise) nondegenerate; for the definition of
3-Sasakian structure see \cite{I-K} and \cite{T}. 
If a Lie group $G$ acts (smoothly) on an almost quaternionic manifold 
$(M,Q)$ preserving $Q$ then there is an induced $\mbox{\J}$-holomorphic 
action on $Z$. Similarly, if a Lie group $G$ acts on a 
quaternionic pseudo-K\"ahler manifold $(M,Q,g)$ preserving the data
$(Q,g)$ then it acts on any of the bundles $Z(M)$, $S(M)$ and 
$U(M)$ preserving all the special geometric structures mentioned above. 

\begin{thm}
Let $(M(\Pi ) = G(\Pi )/K,Q)$ be the homogeneous quaternionic
manifold associated to an extended Poincar\'e algebra of signature
$(p,q)$, $p\ge 3$, $Z(\Pi ) := Z(M(\Pi ))$ its twistor space,
$S(\Pi ) := S(M(\Pi ))$ its canonical ${\rm SO}(3)$-principal 
bundle and $U(\Pi ) := U(M(\Pi ))$ its Swann bundle. Then 
$G(\Pi )$  acts transitively on the manifolds $Z(\Pi )$ and
$S(\Pi )$ and acts on $U(\Pi )$ with an orbit of codimension one. 
 \end{thm}

\noindent
{\bf Proof:} Since $G = G(\Pi )$ acts transitively on the base $M = M(\Pi )$
of any of the bundles $Z(\Pi ) \rightarrow M$, $S(\Pi )  \rightarrow M$ 
and $U(\Pi )  \rightarrow M$, it is sufficient to consider the 
action of the stabilizer $K = {\rm Spin}(3)\cdot {\rm Spin}_0(r,q)$
on the fibres $Z(\Pi )_{[e]}$, $S(\Pi )_{[e]}$ and $U(\Pi )_{[e]}$, 
$[e] = eK \in M = G/K$. The subgroup ${\rm Spin}_0(r,q) \subset K$ acts
trivially on $S(\Pi )_{[e]}$ and hence also on $Z(\Pi )_{[e]}$ and
$U(\Pi )_{[e]}$, whereas ${\rm Spin}(3)$ acts transitively on the
set $S(\Pi )_{[e]}$ of hypercomplex structures subordinate to
$Q_{[e]}$ and hence also on the set $Z(\Pi )_{[e]}$ of complex structures
subordinate to $Q_{[e]}$. From this it follows that $G$ acts transitively
on $Z(\Pi )$ and $S(\Pi )$ and with an orbit of codimension one on $U(\Pi )$.
$\Box$ 

{}From now on we denote by $m_0 := [e] = eK \in M = G/K$ the canonical base 
point of $M$ and fix the complex structure $J_1 \in Q_{m_0}$ as base point 
$z_0 := J_1 \in Z_{m_0}$ in $Z = Z(\Pi )$. 
\begin{cor}
$(Z=Gz_0 \cong G/G_{z_0}, \mbox{\J})$ is a homogeneous complex manifold
of the group $G$. The stabilizer of the point $z_0 \in Z$ in $G$ is
the centralizer of $J_1$ in $K$: $G_{z_0} = {\rm Z}_{K}(J_1) 
= {\rm Z}_{{\rm Spin}(3)}(J_1) \cdot  {\rm Spin}_0(r,q)$,  
${\rm Z}_{{\rm Spin}(3)}(J_1) \cong {\rm U}(1)$. 
\end{cor}

Now we are going to construct a natural holomorphic immersion  
$Z \rightarrow \bar{Z} = \bar{Z}(\Pi ) = G^{\mbox{\C}}/H$ of
$Z$ into a homogeneous complex manifold of the complexified
linear group $G^{\mbox{\C}} \subset {\rm Aut} ({\fr r}^{\mbox{\C}})$,
where $(G_{z_0})^{\mbox{\C}} \subset H \subset G^{\mbox{\C}}$ are
closed complex Lie subgroups.  

First of all, we give an explicit description of the complex structure
$\mbox{\J}$ on the twistor space $Z$. The choice of base point $z_0 \in
Z$ determines a $G$-equivariant diffeomorphism $Z = Gz_0 
\stackrel{\sim}{\rightarrow} G/G_{z_0}$, which maps $z_0$
to the canonical base point $eG_{z_0} \in G/G_{z_0}$. From now on we
will identify $Z$ and $G/G_{z_0}$ via this map.
The complex structure $\mbox{\J}$ being $G$-invariant, it is completely
determined by the $G_{z_0}$-invariant complex structure $\mbox{\J}_{z_0}$
on $T_{z_0}Z$. In order to describe $\mbox{\J}_{z_0}$ we introduce the
following $G_{z_0}$-invariant complement ${\fr z} = {\fr z}(\Pi )$ to
${\fr g}_{z_0} = {\rm Lie} \, G_{z_0} = \mbox{\R}e_2\wedge e_3 \oplus 
{\fr o}(r,q)$ in $\fr g$: 
\[ {\fr z} = \mbox{\R}e_1\wedge e_2 + \mbox{\R}e_1\wedge e_3  + {\fr m} \, .\]
The $G_{z_0}$-invariant decomposition ${\fr g} = {\fr g}_{z_0} + {\fr z}$
determines a $G_{z_0}$-equivariant isomorphism $T_{z_0}Z \cong 
{\fr g}/{\fr g}_{z_0} \stackrel{\sim}{\rightarrow} {\fr z}$. Using it 
we can consider the $G_{z_0}$-invariant complex structure $\mbox{\J}_{z_0}$
as a $G_{z_0}$-invariant complex structure on $\fr z$. 

\begin{prop}
The $G$-invariant complex structure $\mbox{\J}$ on the twistor space
$Z = G/G_{z_0}$ is given on $T_{z_0}Z \cong {\fr z} = 
\mbox{\R}e_1\wedge e_2 + \mbox{\R}e_1\wedge e_3  + {\fr m}$ by:
\[ \mbox{\J}_{z_0}e_1\wedge e_2 = e_1\wedge e_3\, ,\quad 
\mbox{\J}_{z_0}|{\fr m} = J_1\, .\]
\end{prop}

\noindent
{\bf Proof:} Let $\nabla$ be the $G$-invariant quaternionic connection
on $(M,Q)$ constructed in Lemma \ref{NomizuLemma} and $L_x = 
\sum_{\alpha =1}^3 \omega_{\alpha}(x)J_{\alpha} + \bar{L}_x$, $x\in 
{\fr m}$, its Nomizu operators, where $\bar{L}_x \in {\rm z}(Q) \cong 
{\fr gl}(d,\mbox{\Ha})$ ($d =\dim_{\mbox{\Ha}} M$). The connection
$\nabla$ induces the decomposition $T_{z_0}Z  = T_{z_0}Z^{ver} 
\oplus T_{z_0}Z^{hor} \cong {\fr z} = {\fr z}^{ver}\oplus  {\fr z}^{hor}$
into vertical space and horizontal space. The vertical space is 
$T_{z_0}Z^{ver} = T_{z_0}Z_{z_0} = \mbox{\R}J_2 \oplus \mbox{\R}J_3
\subset Q_{z_0}$ and ${\fr z}^{ver} = \mbox{\R}e_1 \wedge e_2 \oplus
\mbox{\R}e_1\wedge e_3$ respectively, the identification being 
$J_2 \mapsto -e_1\wedge e_2$, $J_3 \mapsto -e_1\wedge e_3$. 
{}For any vector $x\in {\fr m}$ we consider the curve $c(t) = 
\exp tx K \in G/K =M$ ($t\in \mbox{\R}$) and define a lift 
$s(t) \in Z_{c(t)}$ by the differential
equation ${\cal L}_Xs = 0$ with initial condition 
$s(0) = z_0 = J_1$. Here $X = \alpha (x)$  is the fundamental vector field 
on $M$ associated to $x$ (as defined on p.\ \pageref{fundvfDef}) and 
${\cal L}_X$ is the Lie derivative with respect to $X$. Then
$s(t) = (\exp tx) z_0$ is precisely the orbit of $z_0\in Z = G/G_{z_0}$
under the 1-parameter subgroup of $G$ generated by $x$. 
The vector 
\begin{eqnarray*} \frac{d}{dt}|_{t=0} (s - t \nabla_Xs) &=& 
s'(0) -  \nabla_Xs|_{t=0} \\
&=& s'(0) + [L_x,J_1] =
s'(0) - 2\omega_2(x)J_3 + 2\omega_3(x)J_2 \in T_{z_0}Z
\end{eqnarray*} 
is horizontal. It is precisely the horizontal lift of 
$X(m_0) \in T_{m_0}M$ and corresponds to 
\begin{equation}
\tilde{x} := x + 2\omega_2(x)e_1\wedge e_3 -2\omega_3(x)e_1\wedge e_2
\in {\fr z} \label{xtildeEqu} \end{equation} 
under the identification $T_{z_0}Z \cong {\fr z}$. 
This shows that
\begin{eqnarray*}
 {\fr z}^{hor} &=& \{ \tilde{x}|x\in {\fr m}\} \\
&=& \mbox{\R} (e_2 +e_1\wedge e_3) + \mbox{\R} (e_3 -e_1\wedge e_2) +\\
& & \mbox{\R}e_0 + \mbox{\R}e_1 + {\rm span} \{ e_{\iota a}| 
\iota = 0,\ldots , 3, \: a = 1,\ldots ,r+q\} + W\, ,
\end{eqnarray*}
see Lemma \ref{NomizuLemma}. Now  the formulas for $\mbox{\J}_{z_0}$
follow easily. In fact, it is clear that $\mbox{\J}_{z_0}$ coincides
with $J_1$ on the $J_1$-invariant subspace ${\fr m} \cap {\fr z}^{hor}
= \mbox{\R}e_0 + \mbox{\R}e_1 + {\rm span} \{ e_{\iota a}| 
\iota = 0,\ldots , 3, \: a = 1,\ldots ,r+q\} + W = e_2^{\perp} \cap e_3^{\perp}
\subset {\fr m}$. 
The equation  $\mbox{\J}_{z_0}e_1\wedge e_2 = e_1\wedge e_3$
follows immediately from $J_1J_2 = J_3$, since $e_1\wedge e_2$, $e_1\wedge 
e_3 \in {\fr z}^{ver}$ are identified with the vertical vectors 
$-J_2$, $-J_3\in T_{z_0}Z^{ver}$. It is now sufficient to check that
$\mbox{\J}_{z_0}e_2 = e_3$. This is done in the next computation:
\[ \mbox{\J}_{z_0}e_2 = \mbox{\J}_{z_0}(e_2 + e_1 \wedge e_3) -
\mbox{\J}_{z_0}(e_1\wedge e_3) = \mbox{\J}_{z_0}\tilde{e}_2 + e_1 \wedge
e_2 = \widetilde{J_1e_2} + e_1 \wedge e_2 = \tilde{e}_3 + e_1 \wedge e_2
= e_3\, . \Box \]

We denote by ${\fr z}^{1,0}$ (respectively, ${\fr z}^{0,1}$) the eigenspace
of $\mbox{\J}_{z_0}\in {\rm End} ({\fr z})$ for the eigenvalue $i$ 
(respectively, $-i$). The integrability of the complex structure $\mbox{\J}$
implies that ${\fr h} := {\fr h}(\Pi ) := ({\fr g}_{z_0})^{\mbox{\C}} + 
{\fr z}^{0,1} \subset {\fr g}^{\mbox{\C}}$ is a (complex) Lie subalgebra.
Let $H = H(\Pi ) \subset G^{\mbox{\C}} \subset {\rm Aut}({\fr r}^{\mbox{\C}})$
\label{H} be the corresponding connected linear Lie group. Let us also consider
the Lie algebra ${\fr h}(V) := {\fr h} \cap {\fr g}(V)^{\mbox{\C}} =
({\fr g}_{z_0})^{\mbox{\C}} + {\fr z}(V)^{0,1} \subset {\fr g}(V)^{\mbox{\C}}$,
where ${\fr z}(V) = {\fr z}\cap {\fr g}(V)$ and ${\fr z}(V)^{0,1} := 
{\fr z}^{0,1} \cap {\fr g}(V)^{\mbox{\C}}$. 

\begin{prop}
\label{algIIProp} $H = H(\Pi ) \subset G^{\mbox{\C}} 
\subset {\rm Aut}({\fr r}^{\mbox{\C}})$ are complex algebraic subgroups. 
\end{prop}

\noindent
{\bf Proof:} It follows from Prop.\ \ref{algProp} that 
${\fr g}^{\mbox{\C}} \subset {\rm der}({\fr r}^{\mbox{\C}})$ is a complex 
algebraic subalgebra and hence $G^{\mbox{\C}} 
\subset {\rm Aut}({\fr r}^{\mbox{\C}})$ a complex 
algebraic subgroup. It only remains to show that $\fr h$ is 
a complex algebraic subalgebra. Let us consider the decomposition
${\fr h} = ({\fr g}_{z_0})^{\mbox{\C}} + {\fr z}^{0,1}$. The 
subalgebra $({\fr g}_{z_0})^{\mbox{\C}}$  is algebraic. In fact,
${\fr g}_{z_0} = {\rm z}_{\fr k}(e_2\wedge e_3)$ is a centralizer
in the real algebraic subalgebra ${\fr k} \subset {\fr g}$. 
If the subalgebra $\langle {\fr z}^{0,1}\rangle \subset {\fr h}$
generated by the subspace ${\fr z}^{0,1}$ is an algebraic subalgebra
of ${\rm der}({\fr r}^{\mbox{\C}})$, then $\fr h$ is generated by algebraic
linear Lie algebras and hence is itself algebraic, see \cite{O-V}. 
The algebraicity of $\langle {\fr z}^{0,1}\rangle$ is proven in the
next lemma. $\Box$ 

\begin{lemma}
${\fr z}^{0,1}$ generates the algebraic subalgebra 
$\langle {\fr z}^{0,1}\rangle = {\fr o}(E')^{\mbox{\C}} +  {\fr z}^{0,1}
\subset {\fr g}^{\mbox{\C}}$. 
\end{lemma}

\noindent
{\bf Proof:} First we compute the subalgebra $\langle {\fr z}^{0,1}\rangle$
of ${\fr g}^{\mbox{\C}}$ generated by ${\fr z}^{0,1} = 
{\fr z}^{0,1}(V) + W^{0,1}$. Note that 
\begin{eqnarray*} {\fr z}^{0,1}(V) &=& 
{\rm span} \{ e_1\wedge e_2 + i e_1\wedge e_3,
e_0+ie_1, e_2+ie_3\}\\ 
& & +\,   {\rm span} \{ e_{0a}+ie_{1a}, e_{2a}+ie_{3a}
|a=1,\ldots , r+q\}\, .\end{eqnarray*} 
It is easy to check that 
\begin{eqnarray*} [{\fr z}^{0,1}(V), {\fr z}^{0,1}(V)] &=&  
{\fr o}(E')^{\mbox{\C}} + {\rm span} \{ e_1\wedge e_2 + i e_1\wedge e_3,
 e_2+ie_3\}\\
 & & +\,    {\rm span} \{  e_{2a}+ie_{3a}
|a=1,\ldots , r+q\}\end{eqnarray*}  
and 
\[ [{\fr z}^{0,1}(V), W^{0,1}] = W^{0,1} \, .\] 
We show that $[W^{0,1}, W^{0,1}] \subset \mbox{\C}(e_2+ie_3)$. 
This shows that $\langle {\fr z}^{0,1}\rangle \supset 
{\fr o}(E')^{\mbox{\C}} +  {\fr z}^{0,1}$. Since the 
right-hand side is closed under Lie brackets, we conclude that
$\langle {\fr z}^{0,1}\rangle = {\fr o}(E')^{\mbox{\C}} +  {\fr z}^{0,1}$. 
It is easy to check that the subalgebra
\[ {\fr o}(E')^{\mbox{\C}} +  {\rm span} \{ e_1\wedge e_2 + i e_1\wedge e_3,
 e_2+ie_3\} +   {\rm span} \{  e_{0a}+ie_{1a}, e_{2a}+ie_{3a}
|a=1,\ldots , r+q\} + W^{0,1}\]
coincides with its derived Lie algebra. Therefore, it is an algebraic
subalgebra of ${\fr g}^{\mbox{\C}}$, see \cite{O-V}. Now, to prove
that $\langle {\fr z}^{0,1}\rangle$ is algebraic, it is sufficient
to check that $\mbox{\C}(e_0+ie_1)$ is algebraic. The element
$e_0+ie_1 \in {\fr g}^{\mbox{\C}}$ is conjugated to the algebraic
element $e_0 \in {\fr g}^{\mbox{\C}}$ via $\exp (-ie_1) \in G^{\mbox{\C}}$. 
Consequently, it is algebraic. $\Box$ 

\begin{thm}\label{EMBThm} For the twistor space $Z$ of any of the homogeneous 
quaternionic ma\-ni\-folds $(M=G/K,Q)$ 
constructed in Thm.\ \ref{mainThm} there is a natural  open  
$G$-equivariant holomorphic immersion
\[ Z \cong G/G_{z_0} \rightarrow 
\bar{Z} = \bar{Z}(\Pi ) := G^{\mbox{\C}}/H\]
into a homogeneous complex (Hausdorff) manifold of the 
complex algebraic group $G^{\mbox{\C}}$. This immersion
is a (universal) finite covering over its image, which is an  
open $G$-orbit. 
(The group $H$ was defined on p.\ \pageref{H}.) 
\end{thm}

\noindent
{\bf Proof:} It follows from Prop.\ \ref{algIIProp} that $H \subset
G^{\mbox{\C}}$ is closed. This shows that $G^{\mbox{\C}}/H$ is a homogeneous
complex Hausdorff manifold. The inclusions $G \subset G^{\mbox{\C}}$ 
and $G_{z_0} \subset H$ define a $G$-equivariant map
$G/G_{z_0} \rightarrow  G^{\mbox{\C}}/H$, which is an immersion since 
${\fr g}\cap {\fr h} = {\fr g}_{z_0}$.  The differential of this immersion at 
$eG_{z_0}$ is canonically 
identified with the restriction 
$\phi :{\fr z} \rightarrow {\fr g}^{\mbox{\C}}/{\fr h}$ 
of the canonical projection ${\fr g}^{\mbox{\C}} \rightarrow 
{\fr g}^{\mbox{\C}}/{\fr h}$ to ${\fr z} \subset {\fr g} = {\fr g}_{z_0}
+ {\fr z}$. Obviously, the complex linear extension $\phi^{\mbox{\C}}$
maps ${\fr z}^{1,0}$ isomorphically to ${\fr g}^{\mbox{\C}}/{\fr h}$
and ${\fr z}^{0,1}$ to zero. This shows that 
$Z\cong G/G_{z_0} \rightarrow  G^{\mbox{\C}}/H$ is open and holomorphic 
with respect to the $G$-invariant complex structure $\mbox{\J}$ on $Z$ and the 
canonical complex structure on $G^{\mbox{\C}}/H$. $\Box$ 

\begin{thm}
The homogeneous complex manifold $\bar{Z} = G^{\mbox{\C}}/H$ carries 
a $G^{\mbox{\C}}$-invariant holomorphic hyperplane distribution
${\cal D} \subset T\bar{Z}$. The hyperplane ${\cal D}_{z_0} =
T_{z_0}Z^{hor} \subset T_{z_0}Z = T_{z_0}{\bar Z}$ is the horizontal
space associated to the $G$-invariant quaternionic connection
$\nabla$ on $M$ constructed in Lemma \ref{NomizuLemma}. Moreover,
$\cal D$ defines a complex contact structure on $\bar{Z} = 
\bar{Z}(\Pi )$ if and only if $\Pi$ is nondegenerate. In this case
the restriction ${\cal D}|Z$ coincides with the canonical complex
contact structure on the twistor space $Z$ of the quaternionic 
pseudo-K\"ahler manifold $M$. 
\end{thm}

\noindent
{\bf Proof:} Recall that we identify $T_{z_0}Z$ with the 
${\fr g}_{z_0}$-invariant subspace ${\fr z} \subset {\fr g} = 
{\fr g}_{z_0} + {\fr z}$ complementary to ${\fr g}_{z_0}$. 
Hereby the subspace ${\cal D}_{z_0} = T_{z_0}Z^{hor}$ is identified
with ${\fr z}^{hor} = \{ \tilde{x}|x\in {\fr m}\} 
\subset {\fr z}$, where $\tilde{x} = x + 2\omega_2(x)e_1\wedge e_3 
-2\omega_3(x)e_1\wedge e_2$ is the $\nabla$-horizontal lift of $x$,
see equation (\ref{xtildeEqu}). The subspace ${\fr z}^{hor} \subset
{\fr z}$ is $\mbox{\J}_{z_0}$-invariant by the very definition
of the complex structure  $\mbox{\J}_{z_0}$ and 
$\mbox{\J}_{z_0}|{\fr z}^{hor}$ is given by $\mbox{\J}_{z_0}\tilde{x}
= \widetilde{J_1x}$ for all $x\in {\fr m}$. The subspace 
$({\fr z}^{hor})^{1,0} = ({\fr z}^{hor})^{\mbox{\C}} \cap 
{\fr z}^{1,0} \subset {\fr z}^{1,0}$ is identified with the 
$i$-eigenspace ${\cal D}_{z_0}^{1,0} 
\subset T^{1,0}_{z_0}Z = T^{1,0}_{z_0}\bar{Z}$ of $\mbox{\J}_{z_0}$ on
${\cal D}_{z_0}$. In order to prove that ${\cal D}_{z_0}^{1,0}$ 
extends to a $G^{\mbox{\C}}$-invariant holomorphic distribution
${\cal D}^{1,0} \subset T^{1,0}\bar{Z}$ it is sufficient to check
the following lemma.

\begin{lemma}
The complex Lie algebra ${\fr h} = {\fr g}_{z_0}^{\mbox{\C}} +
{\fr z}^{0,1}$ preserves the projection of $({\fr z}^{hor})^{1,0}
= \{ \tilde{x} -i\widetilde{J_1x}|x\in {\fr m}\} = 
{\rm span}_{\mbox{\C}}\{ 
e_2 + e_1\wedge e_3 -i(e_3 -e_1\wedge e_2), e_0 -ie_1,
e_{0a}-ie_{1a}, e_{2a} - ie_{3a}|a=1,\ldots , r+q\} + W^{1,0}$ 
into ${\fr g}^{\mbox{\C}}/{\fr h}$, i.e.\ 
\[ [{\fr h}, ({\fr z}^{hor})^{1,0}] \subset {\fr h} +   ({\fr z}^{hor})^{1,0}
\, .\] 
\end{lemma}

Now ${\cal D}^{1,0}$ defines a complex contact structure
if and only if the Frobenius form 
$\wedge^2{\cal D}^{1,0} \stackrel{[\cdot ,\cdot]}{\rightarrow} 
 T^{1,0}\bar{Z}/{\cal D}^{1,0}$ is nondegenerate, 
which is equivalent to the nondegeneracy of the skew symmetric
complex bilinear form
\[ \omega : \wedge^2 ({\fr z}^{hor})^{1,0} 
\stackrel{[\cdot ,\cdot]}{\rightarrow} {\fr g}^{\mbox{\C}} 
\rightarrow {\fr g}^{\mbox{\C}}/({\fr h}   +   ({\fr z}^{hor})^{1,0}) \cong 
\mbox{\C}\, .\] 

\begin{lemma}
Let $({\fr z}^{hor})^{1,0} =  ({\fr z}(V)^{hor})^{1,0} + W^{1,0}$ be the
decomposition of $({\fr z}^{hor})^{1,0}$ induced by the decomposition
${\fr z} = {\fr z}(V) + W$. Then 
$\omega(({\fr z}(V)^{hor})^{1,0}, W^{1,0}) = 0$,\linebreak[3] 
$\omega| \wedge^2({\fr z}(V)^{hor})^{1,0}$ is nondegenerate 
and $\omega|\wedge^2W^{1,0}$
is given by 
\begin{equation} \label{omegaEqu} 
\omega  (s -iJ_1s,t -i J_1t) = 2 (\langle e_2, [s,t]\rangle
+i\langle e_3, [s,t]\rangle)(e_2 -ie_3) \pmod{{\fr h} + 
({\fr z}^{hor})^{1,0}}\, .\end{equation} 
\end{lemma}

{}From the lemma it follows that $\omega$ is nondegenerate if and only if
$\omega|\wedge^2W^{1,0}$ is nondegenerate. The explicit formula for
$\omega|\wedge^2W^{1,0}$ given in equation (\ref{omegaEqu}) now shows that
$\omega$ is nondegenerate if and only if $b(\Pi ) = b_{\Pi , 
(e_1,e_2,e_3)}$ is nondegenerate, which is in turn equivalent to the
nondegeneracy of $\Pi$ by Cor.\ \ref{bpiCor}. This proves that
${\cal D}$ is a complex contact structure if and only if $\Pi$ 
is nondegenerate. If $\Pi$ is  nondegenerate then ${\cal D}_{z_0}
= T_{z_0}Z^{hor}$ is precisely the horizontal space 
associated to the Levi-Civita connection $\nabla^g$ of the
quaternionic  pseudo-K\"ahler manifold $(M,Q,g)$ and $\cal D$
is the canonical complex contact structure on its twistor space $Z$, 
as defined by Salamon \cite{S1}. $\Box$
\section{Homogeneous quaternionic supermanifolds associated to superextended 
Poincar\'e algebras}\label{superSec} 
In this section we will show that our main result, Thm.\ \ref{mainThm}, 
has a natural su\-per\-geo\-me\-tric analogue, Thm.\ \ref{supermainThm}. The
fundamental idea is to replace the map $\Pi : \wedge^2W \rightarrow V$ 
defined on the exterior square $\wedge^2W$ by an ${\fr o}(V)$-equivariant 
linear map $\Pi : \vee^2W \rightarrow V$ defined on the symmetric
square $\vee^2W = {\rm Sym}^2W$. We will freely use the language of 
supergeometry. The necessary background is outlined in the appendix. 

\subsection{Superextended Poincar\'e algebras} \label{sextSec} 
Let $(V,\langle \cdot , \cdot \rangle )$ be a pseudo-Euclidean vector space,
$W$ a $C\! \ell^0 (V)$-module and $\Pi : \vee^2W \rightarrow V$ 
an ${\fr o}(V)$-equivariant linear map. We recall that ${\fr o}(V)$
acts on $W$ via ${\rm ad}^{-1} : {\fr o} (V) 
\rightarrow {\fr spin}(V) \subset C\! \ell^0 (V)$, see equation 
(\ref{ad1Equ}). 

Given these data we extend the Lie bracket on ${\fr p}_0 := {\fr p}(V)$
to a super Lie bracket (see Def.\ \ref{slbDef}) 
$[\cdot ,\cdot ]$ on the $\mbox{\Z}_2$-graded vector space ${\fr p}_0
+{\fr p}_1$, ${\fr p}_1 = W$, by the following requirements:
\begin{enumerate}
\item[1)] The adjoint representation (see Def.~\ref{superadjDef}) 
of ${\fr o}(V)$ on ${\fr p}_1$
concides with the natural representation of ${\fr o}(V) \cong 
{\fr spin}(V)$ on $W = {\fr p}_1$ and $[V,{\fr p}_1] = 0$. 
\item[2)] $[s,t] = \Pi (s\vee t)$ for all $s,t \in W$. 
\end{enumerate}
The super Jacobi identity  follows from 1) and 2). 
The resulting super Lie algebra will be denoted by ${\fr p} (\Pi )$. 

\begin{dof} \label{SEPADef} 
Any super Lie algebra ${\fr p}(\Pi )$ as above is called a 
{\bf superextended Poincar\'e algebra} (of {\bf signature} $(p,q)$
if $V \cong \mbox{\R}^{p,q}$). ${\fr p} (\Pi )$ is called 
{\bf nondegenerate} if $\Pi$  is nondegenerate, 
i.e.\ if the map $W \ni s\mapsto \Pi (s\vee \cdot ) \in W^{\ast} \otimes
V$ is injective. 
\end{dof}
The structure of superextended Poincar\'e algebra on the vector space
${\fr p}(V) + W$ is completely determined by the map $\Pi : \vee^2W \rightarrow
V$. An explicit basis for the vector space 
$(\vee^2W^{\ast}\otimes V)^{{\fr o}(V)}$ of such ${\fr o}(V)$-equivariant
linear maps was constructed in \cite{A-C2} for all $V$ and $W$. 

\subsection{The canonical supersymmetric bilinear form $b$}
Let $V = \mbox{\R}^{p,q}$ be the standard
pseudo-Euclidean vector space with scalar product $\langle \cdot , \cdot
\rangle$ of signature $(p,q)$. 
{}From now on we fix a decomposition $p = p' + p''$ and assume
that $p' \equiv 3\pmod{4}$, see Remark 7 below. 
We denote by $(e_i) = (e_1,\ldots ,e_{p'})$ the first $p'$ basis vectors
of the standard basis of $V$ and by $(e_i') = (e_1',\ldots , e_{p''+q}')$
the remaining ones. The two complementary orthogonal subspaces  
of  $V$ spanned by these bases are denoted by $E =  \mbox{\R}^{p'} = 
\mbox{\R}^{p',0}$ and $E' = \mbox{\R}^{p'',q}$ respectively. The vector 
spaces $V$, $E$ and $E'$ are 
oriented by their standard orthonormal bases. E.g. 
the orientation of Euclidean $p'$-space $E$ defined by the basis 
$(e_i)$ is $e_1^{\ast}\wedge \cdots \wedge 
e_{p'}^{\ast} \in \wedge^{p'} E^{\ast}$. Here $(e_i^{\ast})$ denotes
the basis of $E^{\ast}$ dual to $(e_i)$. Now let 
${\fr p} (\Pi ) = {\fr p} (V) + W$ be a superextended Poincar\'e algebra
of signature $(p,q)$ and $(\tilde{e}_i)$ any orhonormal
basis of $E$. Then we define a $\mbox{\R}$-bilinear form 
$b_{\Pi , (\tilde{e}_i)}$ on the $C\! \ell^0(V)$-module $W$ by: 
\begin{equation} b_{\Pi , (\tilde{e}_i)} (s,t) = \langle \tilde{e}_1, 
[\tilde{e}_2\ldots \tilde{e}_{p'}s,t] \rangle = \langle \tilde{e}_1, \Pi 
(\tilde{e}_2 \ldots \tilde{e}_{p'}s\vee t)\rangle \, ,\quad s,t \in W \, .
\label{superbEqu} 
\end{equation} 
We put $b = b(\Pi ) := b_{\Pi , (e_i)}$ for the standard basis
$(e_i)$ of $E$. As in \ref{sextSec} we consider $W = {\fr p}_1$ as 
$\mbox{\Z}_2$-graded vector space of purely odd degree and recall that,
by Def.\ \ref{superbilDef}, an even supersymmetric (respectively, super skew 
symmetric) bilinear form on $W$ is simply an ordinary skew
symmetric (respectively, symmetric) bilinear form on $W$. 

\noindent 
{\bf Remark 7:} Equation (\ref{superbEqu}) defines an even super skew 
symmetric bilinear form on $W$ if $p'  \equiv 1\pmod{4}$. 
{}For even $p'$ the above formula does
not make sense, unless one assumes that $W$ is a   $C\! \ell (V)$-module
rather than a $C\! \ell^0(V)$-module. Here we are only interested in the case 
$p'  \equiv 3\pmod{4}$. Moreover, later on, for the construction
of homogeneous quaternionic supermanifolds we will put $p' =3$. 

\begin{thm} \label{superbThm} The bilinear form 
$b$ has the following properties:
\begin{enumerate}
\item[1)] $b_{\Pi , (\tilde{e}_i)} = \pm b$ if $\tilde{e}_1
\wedge \cdots \wedge \tilde{e}_{p'} = \pm e_1 \wedge \cdots \wedge 
e_{p'}$. In particular, our definition of $b$ does not depend on the choice
of positively oriented orthonormal basis  of $E$. 
\item[2)] $b$ is an even supersymmetric bilinear form.
\item[3)] $b$ is invariant under the  connected subgroup $K(p',p'') 
= {\rm  Spin}(p') \cdot {\rm Spin}_0(p'',q) \subset 
{\rm Spin}_0(p,q)$ (and is not ${\rm Spin}_0(p,q)$-invariant, 
unless $p''+q=0$). 
\item[4)] Under the identification 
${\fr o} (V) = \wedge^2 V = \wedge^2 E + \wedge^2 E' + E\wedge E'$, 
see equation (\ref{wedgeEqu}), 
the subspace $E\wedge E'$ acts on $W$ by $b$-symmetric endomorphisms
and the subalgebra $\wedge^2 E \oplus \wedge^2 E' \cong {\fr o}(p') 
\oplus {\fr o}(p'',q)$ acts on $W$ by $b$-skew symmetric endomorphisms. 
\end{enumerate}
\end{thm}

\noindent
{\bf Proof:} The proof is the same as for Thm.\ \ref{bThm}, up to
the modifications caused by the fact that $\Pi$ 
is now symmetric instead of  skew symmetric. Part 2) e.g.\ follows 
from the next computation, in which we use that $p' \equiv 3\pmod{4}$: 
\begin{eqnarray*} 
b(t,s) &=& \langle e_1,[e_2 \ldots e_{p'}t,s]\rangle\\
&=& -\langle e_1 ,[e_4\ldots e_{p'}t,e_2e_3s]\rangle + \langle e_1 , 
{\rm ad} (e_2e_3)
[e_4 \ldots e_{p'}t,s]\rangle\\
&=& -\langle e_1 ,[e_4 \ldots e_{p'}t,e_2e_3s]\rangle = \cdots =
-\langle e_1 ,[t,e_2 \ldots e_{p'}s] \rangle \\
&=& - \langle e_1 ,[e_2 \ldots e_{p'}s,t] \rangle  = -b(s,t)\, . \Box
\end{eqnarray*}

\begin{dof} The bilinear form $b = b(\Pi ) = b_{\Pi ,(e_1,\ldots , e_{p'})}$ 
defined above is called
the {\bf canonical supersymmetric bilinear form} on $W$ associated 
to the ${\fr o}(V)$-equivariant map $\Pi :\vee^2W 
\rightarrow V = \mbox{\R}^{p,q}$ and the decomposition $p = p' + p''$.  
\end{dof}

\begin{prop} \label{superkerProp} The kernels of the linear 
maps $\Pi : W \rightarrow 
W^{\ast}\otimes V$ and $b = b(\Pi ) : W \rightarrow W^{\ast}$ coincide:
${\rm ker}\, \Pi = {\rm ker}\, b$. 
\end{prop}

\noindent {\bf Proof:} See the proof of Prop.\ \ref{kerProp}. $\Box$

\begin{cor} \label{superbpiCor} ${\fr p} (\Pi )$ is nondegenerate 
(see Def.\ \ref{SEPADef}) if and only if 
$b (\Pi )$ is non\-de\-ge\-ne\-rate.
\end{cor} 

\subsection{The main theorem in the super case}
Any superextended Poincar\'e algebra ${\fr p} = {\fr p}(\Pi) =
{\fr p}(V) + W$ has an even derivation $D$ with eigenspace 
decomposition ${\fr p} = {\fr o}(V) + V + W$  and corresponding
eigenvalues $(0,1,1/2)$.  Therefore, the super Lie algebra 
${\fr p} = {\fr p}_0 + {\fr p}_1 = {\fr p}(V) + W$ is canonically
extended to a super Lie algebra ${\fr g} = {\fr g}(\Pi ) 
= \mbox{\R} D + {\fr p} = {\fr g}_0 + {\fr g}_1$, where 
${\fr g}_0 = \mbox{\R} D + {\fr p}_0 =
\mbox{\R} D + {\fr p}(V) = {\fr g}(V)$ and ${\fr g}_1 =
{\fr p}_1 = W$.

\begin{prop} \label{superfaithProp} The adjoint representation (see 
Def.\ \ref{superadjDef}) of 
${\fr g}$ is faithful and moreover it induces a faithful representation on
its ideal ${\fr r} = \mbox{\R}D + V +W \subset {\fr g} = {\fr o}(V) + {\fr r}$.
\end{prop}
By Prop.\ \ref{superfaithProp} we can consider $\fr g$ as subalgebra
of the super Lie algebra ${\fr gl}({\fr r})$ (defined in \ref{supergroupSec}), 
i.e.\ $\fr g$ is a linear super Lie algebra (see Def.\ \ref{superadjDef}).
We denote by $G = G(\Pi )$ the corresponding linear Lie supergroup,
see \ref{supergroupSec}. Its 
underlying Lie group is $G_0 = G(V)$ and has ${\fr g}_0 = {\fr g}(V)$ as 
Lie algebra. For the construction of homogeneous quaternionic supermanifolds
we will assume that $V = {\mbox{\R}}^{p,q}$ with $p\ge 3$. Then we fix the 
decomposition $p = p' + p''$, where now $p' = 3$ and $p'' = p-3 =: r$. As 
before, we have a corresponding orthogonal decomposition $V = E + E'$, the
subalgebra ${\fr k} = {\fr k}(p',p'') =  {\fr k}(3,r) = {\fr o}(3)
\oplus {\fr o}(r,q) \subset {\fr o}(p,q) = {\fr o}(V)$ preserving 
this decomposition and the corresponding linear Lie subgroup $K = K(3,r)
\subset  G_0$. We are interested in the homogeneous supermanifold
(see \ref{superhomogSec}):
\[ M = M(\Pi ) := G/K = G(\Pi ) /K \, .\]
Its underlying manifold is the homogeneous manifold
\[ M_0 := G_0/K = G(V)/K = M(V)\, .\]

\begin{thm}\label{supermainThm}
\begin{enumerate}
\item[1)] There exists a $G$-invariant quaternionic structure $Q$ on
$M = G/K$.
\item[2)] If $\Pi$ is nondegenerate (see Def.~\ref{SEPADef}) then
there exists a $G$-invariant pseudo-Riemannian metric $g$ on $M$
such that $(M,Q,g)$ is a quaternionic pseudo-K\"ahler supermanifold.
\end{enumerate}
\end{thm}

\noindent
{\bf Proof:} First let $(Q_0,g_0)$ denote the $G_0$-invariant 
quaternionic pseudo-K\"ahler structure on $M_0 = G_0/K$
which was introduced in the proof of Thm.\ \ref{mainThm}
(previously it was denoted simply by $(Q,g)$). 
As in that proof, see equation (\ref{extEqu}), 
the corresponding $K$-invariant quaternionic structure on 
$T_{eK}M_0 \cong {\fr g}_0/{\fr k}$ is extended to a $K$-invariant quaternionic
structure $Q_{eK}$ on $T_{eK}M \cong {\fr g}/{\fr k}$ defining a $G$-invariant
almost quaternionic structure $Q$ on $G/K$ (see \ref{supermetSec} and 
\ref{superhomogSec}). Similarly, the 
$G_0$-invariant pseudo-Riemannian metric $g_0$ on $M_0$ corresponds
to a $K$-invariant pseudo-Euclidean scalar product on ${\fr g}_0/{\fr k}$.
This scalar product is extended 
by $b = b_{\Pi , (e_1,e_2,e_3)}$, in the obvious way,   to
a $K$-invariant supersymmetric bilinear form $g_{eK}$ on ${\fr g}/{\fr k}$,
which is nondegenerate if $\Pi$ is nondegenerate.   Let us first treat 
the case where $\Pi$ is nondegenerate. In this case $g_{eK}$
defines a $G$-invariant $Q$-Hermitian pseudo-Riemannian metric
$g$ on $G/K$. The Levi-Civita connection of the homogeneous
pseudo-Riemannian supermanifold $(M,g)$ is computed in Lemma 
\ref{superLCLemma} below.
As on  p.\ \pageref{mPage},  ${\fr g}/{\fr k}$ is identified 
with the complement ${\fr m} = {\fr m}(V) + W$ to $\fr k$ in $\fr g$. 
We use the $g_0$-orthonormal basis $(e_{\iota i})$ of ${\fr m}(V)$
introduced on p.\ \pageref{eiotaiPage} and recall that $\epsilon_i=
g_0(e_{\iota i},e_{\iota i})$. 
Also we will continue to write $x\wedge_g  y$ for the $g_{eK}$-skew symmetric
endomorphism of $\fr m$ defined for $x,y\in {\fr m}(V)$ by:  
$x\wedge_g  y (z) := g_{eK}(y,z)x - g_{eK}(x,z)y$, $z\in {\fr m}$. 

\begin{lemma} \label{superLCLemma}
The Nomizu map $L^g = L(\nabla^g)$ associated to the Levi-Civita connection 
$\nabla^g$ of the homogeneous pseudo-Riemannian supermanifold $(M = G/K,g)$
is given by the following formulas:
\begin{eqnarray*}
L^g_{e_0} &=&  L^g_{e_{\alpha a}} = 0\, ,\\
L^g_{e_{\alpha}} &=& \frac{1}{2}J_{\alpha} + \bar{L}^g_{e_{\alpha}}\, ,
\end{eqnarray*}
where 
\[ \bar{L}^g_{e_{\alpha}} = \frac{1}{2}\sum_{i=0}^{r+q}
\epsilon_i e_{\alpha  i} \wedge_ge_i' - \frac{1}{2}\sum_{i=0}^{r+q}
\epsilon_i e_{\gamma i}  \wedge_ge_{\beta i} \in {\rm z}(Q_{eK})\, ,\]
$L^g_{e_a'} \in {\rm z}(Q_{eK})$ is given by:
\begin{eqnarray*}
L^g_{e_a'} |{\fr m}(V) &=& \sum_{\iota =0}^3 e_{\iota a} \wedge_g e_{\iota}
\, ,\\
L^g_{e_a'} |W &=& \frac{1}{2} e_1e_2e_3e_a'\, .
\end{eqnarray*}
{}For all $s\in W$ the Nomizu operator $L_s^g\in {\rm z} (Q_{eK})$ maps the
subspace ${\fr m}(V) \subset {\fr m} = {\fr m}(V) + W$ into $W$ and
$W$ into ${\fr m}(V)$. The restriction $L_s^g|W$ ($s\in W$) is completely 
determined by $L^g_s|{\fr m}(V)$ (and vice versa) according to the
relation
\[ g_{eK}(L_s^gt,x) = g_{eK}(t,L_s^gx)\, ,\quad s,t\in W\, ,\quad 
x\in {\fr m}(V)\, .
\]
Finally,  $L^g_s|{\fr m}(V)$ ($s\in W$) is completely 
determined by its values on the quaternionic basis $(e_i')$, 
$i = 0,\ldots , r+q$, which are as follows:
\begin{eqnarray*}
L_s^ge_0' &=& \frac{1}{2}s\, ,\\
L_s^ge_a' &=& \frac{1}{2}e_1e_2e_3e_a's\, .
\end{eqnarray*}
(It is understood that $L_x^g  = {\rm ad}_x|{\fr m}$ for all
$x\in {\fr k}$, cf.\ equation (\ref{superNomizuIEqu}).)
In the above formulas $a = 1, \ldots , r+q$, $\alpha = 1,2,3$ and 
$({\alpha},{\beta},{\gamma})$ is a 
cyclic permutation of $(1,2,3)$. 
\end{lemma}

\noindent
{\bf Proof:} This follows from equation (\ref{superKoszulEqu}) by a 
straightforward computation. $\Box$ 

\begin{cor}
\label{superLCCor} The Levi-Civita connection $\nabla^g$ of the homogeneous
almost quaternionic pseudo-Hermitian supermanifold $(M,Q,g)$ preserves
$Q$ and hence $(M,Q,g)$ is a quaternionic  pseudo-K\"ahler supermanifold.
\end{cor}

By Cor.\ \ref{superLCCor} we have already established part 2) of 
Thm.\ \ref{supermainThm}. Part 1) is a consequence of 2) provided that
$\Pi$ is nondegenerate. It remains to discuss the case of degenerate
$\Pi$. From the ${\fr o}(V)$-equivariance of $\Pi$ it follows that
$W_0 = {\rm ker}\, \Pi \subset W$ is an ${\fr o}(V)$-submodule.
Let $W'$ be a complementary ${\fr o}(V)$-submodule. Then 
$W_0$ and $W'$ are   $C\! \ell^0(V)$-submodules; we put  
$\Pi':= \Pi|\vee^2W'$. We denote by $(M':= M(\Pi'),
Q',g')$ the corresponding quaternionic  pseudo-K\"ahler supermanifold and  by 
$L' := L^{g'}: {\fr g}(\Pi') \rightarrow
{\rm End}({\fr m}(\Pi'))$ the Nomizu map associated to its Levi-Civita
connection. By the next lemma, 
we can extend the map $L'$ to a torsionfree Nomizu map 
$L:{\fr g}(\Pi) \rightarrow {\rm End}({\fr m}(\Pi))$, whose
image normalizes $Q_{eK}$. This proves the 1-integrability of $Q$
(by Cor.~\ref{superquatCor}), completing the proof of
Thm.~\ref{supermainThm}. $\Box$  

By Cor.\ \ref{superLCCor} we can decompose $L_x' = \sum_{\alpha=1}^3
\omega_{\alpha}'(x)J_{\alpha} + \bar{L}_x'$, where 
the $\omega_{\alpha}'$ are 1-forms on ${\fr m}(\Pi')$ and 
$\bar{L}_x' \in {\rm z}({Q'}_{eK})$ belongs to the centralizer of the 
quaternionic structure ${Q'}_{eK} = Q_{eK}|{\fr m}(\Pi')$ on ${\fr m}(\Pi')$. 

\begin{lemma} \label{superNomizuLemma} 
The Nomizu map $L':{\fr g}(\Pi') \rightarrow {\rm End}({\fr m}(\Pi'))$ 
associated to the Levi-Civita connection of $(M(\Pi') = G(\Pi')/K,g')$
can be extended to the Nomizu map $L:{\fr g}(\Pi) \rightarrow 
{\rm End}({\fr m}(\Pi))$ of a  $G(\Pi)$-invariant 
quaternionic connection $\nabla$ on the homogeneous almost 
quaternionic supermanifold
$(M(\Pi) = G(\Pi)/K,Q)$. The extension is defined as follows:
\[ L_x:= \sum_{\alpha=1}^3 \omega_{\alpha}(x)J_{\alpha} + \bar{L}_x\, ,
\quad x\in {\fr m}(\Pi)\, ,\]
where $\bar{L}_x \in{\rm z}(Q)$ (the centralizer is taken 
in ${\fr gl}({\fr m}(\Pi))$) is defined below and 
the 1-forms $\omega_{\alpha}$ on ${\fr m}(\Pi) := {\fr m}(\Pi') + W_0$
satisfy $\omega_{\alpha}|{\fr m}(\Pi') := \omega_{\alpha}'$ and 
$\omega_{\alpha}|W_0 := 0$. The operators $\bar{L}_x$ are given by:
\[ \bar{L}_x|{\fr m}(\Pi') := \bar{L}_x' \quad \mbox{if} \quad
x\in {\fr m}(\Pi')\, ,\]
\[ \bar{L}_{e_0}|W_0 := \bar{L}_{e_{\alpha a}}|W_0 
:= \bar{L}_{e_{\alpha}}|W_0 := 0\, ,
\]
\[ \bar{L}_{e_a'}|W_0 := \frac{1}{2} e_1e_2e_3e_a'\, ,\]
\[ \bar{L}_s|W_0 := 0 \quad \mbox{if} \quad s\in W\, ,\]
\[ \bar{L}_sx := (-1)^{\tilde{x}}L_xs -[s,x] \quad \mbox{if} \quad s\in W_0\, ,\quad
x\in {\fr m}(\Pi')\, .\]
(It is understood that $L_x = {\rm ad}_x|{\fr m}(\Pi)$ for all $x\in{\fr k}$.)
In the above formulas, as usual, $a = 1, \ldots , r+q$, 
$\alpha = 1,2,3$ and $\tilde{x} \in \mbox{\Z}_2 = \{ 0,1\}$ 
stands for the $\mbox{\Z}_2$-degree
of $x\in {\fr m}(\Pi ') = {\fr m}(\Pi ')_0 + {\fr m}(\Pi ')_1 =
{\fr m}(V) + W'$. 
\end{lemma}

\noindent
{\bf Proof:} The proof is similar to that of Lemma \ref{NomizuLemma}. $\Box$

\begin{appendix}
\section{Supergeometry}\label{A1Sec}
In this appendix we summarize the supergeometric material needed
in \ref{superSec}. Standard references on supergeometry are
\cite{M}, \cite{L} and \cite{K}, see also \cite{B-B-H}, \cite{Ba}, \cite{Be}, 
\cite{Bern}, \cite{B-O}, \cite{DW}, \cite{F}, \cite{O1}, \cite{O2}, 
\cite{Sch} and \cite{S-W}.  (D.A.\ Leites has informed us that he will
soon publish a monograph on supergeometry.) 
\subsection{Supermanifolds}
Let $V = V_0 + V_1$ be a $\mbox{\Z}_2$-graded vector space. 
We recall that  an element
$x\in V$ is called {\bf homogeneous} (or {\bf of pure degree})
if $x\in V_0 \cup V_1$. The {\bf degree} of a homogeneous
element $x\in V$ is the number $\tilde{x} \in \mbox{\Z}_2 = \{ 0, 1\}$
such that $x\in V_{\tilde{x}}$. The element $x\in V$ said to be even
if $\tilde{x} = 0$ and odd if $\tilde{x} = 1$. 
{}For $V_0$ and $V_1$ of finite
dimension, the {\bf dimension} of $V$ is defined
as $\dim V := \dim V_0|\dim V_1 = m|n$ and a {\bf basis} of $V$ is 
by definition a tuple $(x_1,\ldots , x_m,
\xi_1, \ldots , \xi_n)$ such that $(x_1,\ldots , x_m)$ is a basis
of $V_0$ and $(\xi_1, \ldots , \xi_n)$ is a basis of $V_1$. 

\begin{dof}\label{superalgDef}
Let $\bf A$ be a $\mbox{\Z}_2$-graded algebra. The {\bf supercommutator}
is the bilinear map $[\cdot , \cdot  ] : {\bf A}\times {\bf A} 
\rightarrow {\bf A}$ 
defined by:  
\[ [a,b] := ab - (-1)^{\tilde{a}\tilde{b}}ba\]
for all homogeneous elements $a,b \in {\bf A}$. The algebra 
$\bf A$ is called {\bf 
supercommutative} if $[a,b] = 0$ for all $a,b \in {\bf A}$. 
A  $\mbox{\Z}_2$-graded supercommutative 
associative (real) algebra $A = A_0 + A_1$  will simply be called
a {\bf superalgebra}. 
\end{dof}

\noindent 
{\bf Example 1:} The exterior algebra $\wedge E = \wedge^{even}E 
+\wedge^{odd}E$ over a finite dimensional vector space $E$ is a 
superalgebra. 

\begin{dof}
\label{slbDef} A {\bf super Lie bracket} on a 
$\mbox{\Z}_2$-graded vector space  
$V = V_0 + V_1$ is a bilinear map $[ \cdot , \cdot ] : V \times V 
\rightarrow V$ such that for all $x,y,z \in V_0 \cup V_1$ we 
have:  
\begin{enumerate} 
\item[i)] $\widetilde{[x,y]} = \tilde{x} +  \tilde{y}$,
\item[ii)] $[x,y] = - (-1)^{\tilde{x}\tilde{y}} [y,x]$ and 
\item[iii)] $[x,[y,z]] = [[x,y],z] + (-1)^{\tilde{x}\tilde{y}} [y,[x,z]]$ 
(``super Jacobi identity''). 
\end{enumerate} 
The $\mbox{\Z}_2$-graded algebra with underlying $\mbox{\Z}_2$-graded 
vector space $V$ and product defined by the super Lie bracket 
$[\cdot , \cdot ]$ is called a {\bf super Lie algebra}. 
\end{dof} 

\noindent 
{\bf Example 2:} 
The supercommutator of any associative $\mbox{\Z}_2$-graded algebra
$\bf A$  is a super Lie bracket and hence defines on it the 
structure of super Lie algebra. For example, we may take ${\bf A} =
{\rm End}(V)$ with the obvious structure of 
$\mbox{\Z}_2$-graded associative algebra (with unit). The  corresponding
super Lie algebra is denoted by ${\fr gl}(V)$ and is called the 
{\bf general linear super Lie algebra}.

Let $M_0$ be a (differentiable) manifold of dimension $m$. We denote
by ${\cal C}_{M_0}^{\infty}$ its sheaf of functions. Sections of
the sheaf ${\cal C}_{M_0}^{\infty}$ over an open set $U \subset M_0$
are simply smooth functions on $U$: ${\cal C}_{M_0}^{\infty}(U) = 
C^{\infty}(U)$. Now let ${\cal A} = {\cal A}_0 + {\cal A}_1$ be a sheaf
of superalgebras over $M_0$. 

\begin{dof}\label{supermfDef}
The pair $M = (M_0,{\cal A})$ is called a (differentiable) 
{\bf supermanifold} of {\bf dimension} $\dim M = m|n$ over $M_0$ 
if for all $p\in M_0$
there exists an open neighborhood $U\ni p$ and a rank 
$n$ free sheaf
${\cal E}_U$ of ${\cal C}^{\infty}_{U}$-modules over 
$U$ such that 
${{\cal A}|}_U \cong \wedge {\cal E}_U$ (as sheaves of 
superalgebras). A  {\bf function} 
on $M$ (over an open set $U \subset M_0$) is by definition 
a  section of $\cal A$ (over $U$). 
The sheaf ${\cal A} = {\cal A}_M$ 
is called the {\bf sheaf of functions} on $M$
and $M_0$ is called the  manifold {\bf underlying} the supermanifold $M$. 
Let  $M = (M_0, {\cal A}_M)$ and 
$N = (N_0, {\cal A}_N)$ be 
supermanifolds. A {\bf morphism} $\varphi : M\rightarrow N$ is a pair
$\varphi = (\varphi_0, \varphi^{\ast})$, where $\varphi_0 :
M_0 \rightarrow N_0$ is a smooth map and $\varphi^{\ast} :
{\cal A}_N \rightarrow (\varphi_0)_{\ast}{\cal A}_M$ is a 
morphism of sheaves of superalgebras. It is called
an {\bf isomorphism} if $\varphi_0$ is a diffeomorphism and
$\varphi^{\ast}$ an isomorphism. An isomorphism
$\varphi : M \rightarrow M$ is called an {\bf automorphism} of $M$.
The set of all morphisms $\varphi : M\rightarrow N$ (respectively,
automorphisms $\varphi : M \rightarrow M$) is denoted by 
${\rm Mor}(M,N)$ (respectively, ${\rm Aut}(M))$.
\end{dof}
{}From Def.\ \ref{supermfDef} it follows that there exists a canonical 
epimorphism of sheaves $\epsilon^{\ast} : 
{\cal A} \rightarrow {\cal C}^{\infty}_{M_0}$, which is called the 
{\bf evaluation map}. Its kernel is the ideal
generated by ${\cal A}_1$: $\ker \epsilon^{\ast} = 
\langle {\cal A}_1 \rangle = {\cal A}_1 + {\cal A}_1^2$. 

Given supermanifolds $L$, $M$, $N$ and morphisms $\psi \in {\rm Mor}(L,M)$ 
and $\varphi \in {\rm Mor}(M,N)$, there is a {\bf composition} 
$\varphi \circ \psi \in {\rm Mor}(L,N)$  defined by:
\[(\varphi \circ \psi )_0 =  \varphi_0 \circ \psi_0 \quad \mbox{and} \quad  
(\varphi \circ \psi )^{\ast} = \psi^{\ast}\circ \varphi^{\ast}\, . \]
Here we have used the same symbol $\psi^{\ast}$ for the map
$(\varphi_0)_{\ast} {\cal A}_M \rightarrow (\varphi_0)_{\ast}(\psi_0)_{\ast}
{\cal A}_L = (\varphi_0 \circ \psi_0)_{\ast} {\cal A}_L$ induced by 
$\psi^{\ast} :  {\cal A}_M \rightarrow (\psi_0)_{\ast} {\cal A}_L$. 
Similarly, if $\varphi : M \rightarrow N$ is an isomorphism, then we can define
it {\bf inverse} isomorphism by: 
\[ \varphi^{-1} := (\varphi_0^{-1}, (\varphi^{\ast})^{-1}) :
N \rightarrow M\, . \]
Here, again, we have used the same notation 
$(\varphi^{\ast})^{-1}$ for the map ${\cal A}_M \rightarrow
(\varphi_0^{-1})_{\ast}{\cal A}_N$ induced by 
$(\varphi^{\ast})^{-1} : (\varphi_0)_{\ast} {\cal A}_M \rightarrow 
{\cal A}_N$. Finally, for every supermanifold $M = (M_0, {\cal A})$, 
there is the {\bf identity} automorphism ${\rm Id}_M := ({\rm Id}_{M_0}, {\rm 
Id}_{{\cal A}})$. The above operations turn the set ${\rm Aut}(M)$ 
into a group. 

\noindent 
{\bf Example 3:} Let $E \rightarrow M_0$ be a (smooth) vector bundle
of rank $n$ over the $m$-dimensional manifold 
$M_0$ and ${\cal E}$ its sheaf of sections. It is a 
locally free sheaf of ${\cal C}^{\infty}_{M_0}$-modules and 
$SM(E):= (M_0,\wedge{\cal E})$ is a supermanifold of dimension $m|n$. 
Its evaluation map is the canonical projection onto $0$-forms
$\wedge {\cal E} =  {\cal C}^{\infty}_{M_0} + \sum_{j=1}^n\wedge^j{\cal E}
\rightarrow {\cal C}^{\infty}_{M_0}$. It is well known, see \cite{Ba},
that any supermanifold is isomorphic 
to a supermanifold of the form
$SM(E)$. However, the isomorphism is not canonical, unless $n=0$. 

\noindent 
{\bf Example 4:} Any manifold $(M_0, {\cal C}^{\infty}_{M_0})$ of 
dimension $m$ can be considered as a supermanifold of dimension $m|0$. 
In fact, it is associated to the vector bundle of rank $0$ over $M_0$
via the construction of Example 3. For any supermanifold $M =
(M_0,{\cal A})$ the pair 
$({\rm Id}_{M_0},\epsilon^{\ast})$ defines a canonical morphism 
$\epsilon : M_0 \rightarrow M$. The composition of $\epsilon$
with the canonical constant map $p: \{ p\} \rightarrow M_0$ ($p\in M_0$)
defines a morphism $\epsilon_p = (p,\epsilon^{\ast}_p): \{ p\} \rightarrow
M$. The epimorphism $\epsilon^{\ast}_p: {\cal A} \rightarrow \mbox{\R}$
onto the constant sheaf is called the {\bf evaluation at}  $p$:
\[ \epsilon^{\ast}_pf  = (\epsilon^{\ast}f)(p)  \, .\] 
$f(p) := \epsilon^{\ast}_pf \in \mbox{\R}$ is called the {\bf value} of 
$f$ {\bf at} the point $p$.  

\noindent 
{\bf Example 5:} Let $V = V_0 + V_1$ be a $\mbox{\Z}_2$-graded vector space
of dimension $m|n$ and $E_V = V_1 \times V_0 \rightarrow V_0$ the trivial
vector bundle over $V_0$ with fibre $V_1$.  Then to $V$ we can
canonically associate the supermanifold $SM(V) := SM(E_V)$, see Example 3.

Let $(x_0^i) = (x^1_0,\ldots , x^m_0)$ be local coordinates for $M_0$
defined on an open set $U \subset M_0$, 
${\cal E}_U$ a rank $n$ free sheaf of ${\cal C}^{\infty}_{U}$-modules
over $U$ and $\phi : \wedge{\cal E}_U \rightarrow {{\cal A}|}_U$ an 
isomorphism.
We can choose sections $(\xi^j_0) = (\xi_0^1,\ldots , \xi_0^n)$ of 
${\cal E}_U$ which generate  ${\cal E}_U$ freely over ${\cal C}^{\infty}_U$. 
Then any section of $\wedge {\cal E}_U$ is of the form:
\begin{equation} \label{fEqu} 
f = \sum_{\alpha \in \mbox{\Z}_2^n} f_{\alpha} 
(x^1_0,\ldots , x^m_0)\xi_0^{\alpha}\, , \quad 
f_{\alpha} (x^1_0,\ldots , x^m_0) \in C^{\infty}(U) 
\, ,
\end{equation}   
where $\xi_0^{\alpha} := (\xi_0^1)^{\alpha_1} \wedge \ldots \wedge 
(\xi_0^n)^{\alpha_n}$ for $\alpha = (\alpha_1,\ldots ,\alpha_n)$. 
The tuple $(\phi,x_0^i,\xi_0^j)$ is called a {\bf local coordinate
system} for $M$ over $U$. The open set $U\subset M_0$
is called a {\bf coordinate neighborhood} for $M$. Any function
on $M$ over $U$ is of the form $\phi (f)$ for some section $f$ as in
(\ref{fEqu}). The functions $x^i := \phi (x_0^i) \in {\cal A}(U)_0$,
$\xi^j := \phi(\xi_0^j) \in {\cal A}_1$ are called {\bf local
coordinates} for $M$ over $U$. The evaluation map 
$\epsilon^{\ast}: {\cal A} \rightarrow {\cal C}^{\infty}_{M_0}$ 
is expressed in a local coordinate system simply by putting
$\xi_0^1 = \cdots = \xi_0^n = 0$ in (\ref{fEqu}):
\[ \epsilon^{\ast}(\phi (f)) = f(x_0^1,\ldots , x_0^m, 0,\ldots ,0)
:= f_{(0,\ldots ,0)}(x_0^1,\ldots , x_0^m)\, .\]
In particular, we have $\epsilon^{\ast}x^i = x_0^i$ and 
$\epsilon^{\ast}\xi^j = 0$. 

\noindent
{\bf Example 6:} Let $V = V_0 + V_1$ be a $\mbox{\Z}_2$-graded vector space
and $(x^i,\xi^j) = (x^i,\ldots , x^m,$ $\xi^1, \ldots , \xi^n)$ a basis of the 
$\mbox{\Z}_2$-graded vector space $V^{\ast} = {\rm Hom}(V,\mbox{\R})$.
Then $(x^i,\xi^j)$ can be considered as global coordinates on the 
supermanifold $SM(V)$, i.e.\ as local coordinates for $SM(V)$ over
$V_0 = SM(V)_0$, see Example 5. 
 
Let $M$ and $N$ be supermanifolds of dimension $m|n$ and $p|q$ respectively.
In local coordinates $(x^i,\xi^j)$ for $M$ and 
$(y^k,\eta^l)$ for $N$ a morphism $\varphi$ is expressed by
$p$ even functions $y^k(x^1,\ldots , x^m, \xi^1, \ldots , \xi^n) := 
\varphi^{\ast}y^k$ and $q$ odd functions 
$\eta^l(x^1,\ldots , x^m, \xi^1, \ldots , \xi^n)$ $:= \varphi^{\ast}\eta^l$.

There exists a supermanifold $M\times N = (M_0 \times N_0, 
{\cal A}_{M\times N})$ called the {\bf product} of the
supermanifolds $M$ and $N$
and morphisms $\pi_M : M\times N \rightarrow M$, $\pi_N :
M\times N \rightarrow N$ such that 
$(\pi^{\ast}_Mx^i, \pi^{\ast}_Ny^k, \pi^{\ast}_M\xi^j, \pi^{\ast}_N\eta^l)$
are local coordinates for $M\times N$ over $U\times V$ if 
$(x^i,\xi^j)$ are local coordinates for $M$ over $U$ and 
$(y^k,\eta^l)$  are local coordinates for $N$ over $V$.
The morphism $\pi_1 = \pi_M$ (respectively, $\pi_2 = \pi_N$) is called the
{\bf projection} of $M\times N$ onto the first (respectively, second)
factor. Given morphisms $\varphi_i :M_i \rightarrow N_i$, $i=1,2$, 
there is a corresponding morphism $\varphi_1 \times \varphi_2 : 
M_1 \times M_2 \rightarrow N_1 \times N_2$ such that
$\pi_{N_1}  \circ (\varphi_1 \times \varphi_2) = \varphi_1 \circ \pi_{M_1}$
and 
$\pi_{N_2} \circ (\varphi_1 \times \varphi_2) = \varphi_2 \circ \pi_{M_2}$. 

As next, we will discuss the notion for tangency on supermanifolds. 
{}For this purpose, we recall the following definition.
\begin{dof}
An endomorphism $X = X_0 + X_1 \in {\rm End}({\bf A}) = 
{\rm End}_{\mbox{\R}}({\bf A})$ (here $\tilde{X}_{\alpha} = \alpha$, 
$\alpha = 0,1$) 
of a $\mbox{\Z}_2$-graded algebra ${\bf A}$ is called
a {\bf derivation} if it satisfies the Leibniz-rule
\[ X_{\alpha}(ab) = X_{\alpha}(a)b + (-1)^{\alpha \tilde{a}}aX_{\alpha}(b)\]
for all homogeneous $a,b \in {\bf A}$ and $\alpha \in \mbox{\Z}_2$.
The $\mbox{\Z}_2$-graded vector space of all derivations of $\bf A$
is denoted by ${\rm Der}\, {\bf A} = ({\rm Der}\, {\bf A})_0 + 
({\rm Der}\, {\bf A})_1$.
\end{dof}
Notice that the supercommutator on ${\rm End}({\bf A})$
restricts to a super Lie bracket on ${\rm Der}\, {\bf A}$. 

\begin{dof}
Let $M = (M_0,{\cal A})$ be a supermanifold. The {\bf tangent
sheaf} of $M$ is the sheaf of derivations of $\cal A$ and is 
denoted by ${\cal T}_M = ({\cal T}_M)_0 + ({\cal T}_M)_1$.
A {\bf vector field} on $M$ is a section of ${\cal T}_M$. 
The {\bf cotangent sheaf} is the sheaf ${\cal T}_M^{\ast} =
{\rm Hom}_{\cal A}({\cal T}_M,{\cal A})$. The {\bf full tensor superalgebra}
over ${\cal T}_M$ is the sheaf of superalgebras generated
by tensor products ($\mbox{\Z}_2$-graded over $\cal A$) of 
${\cal T}_M$ and ${\cal T}_M^{\ast}$. It is denoted by 
$\otimes_{\cal A} \langle {\cal T}_M, {\cal T}_M^{\ast}\rangle$. 
A {\bf tensor field} on $M$ is a section of 
$\otimes_{\cal A} \langle {\cal T}_M, {\cal T}_M^{\ast}\rangle$.
\end{dof}

Explicitly, a section $X \in {\cal T}_M(U)$ ($U \subset M_0$ open) associates
to any open subset $V \subset U$ a derivation ${X|}_V \in  {\rm Der}\, 
{\cal A}(V)$ such that
\[ {X|}_V ({f|}_V) = {{X|}_U(f)|}_V\quad \mbox{for all} \quad 
f \in {\cal A}(U) \, .\] 

\noindent 
Given local coordinates $(x^i,\xi^j)$ on $M$ over $U$ there exist
unique vector fields 
\[ \frac{\partial}{\partial x^i}\in {\cal T}_M(U)_0 \quad
\mbox{and} \quad \frac{\partial}{\partial \xi^j}\in {\cal T}_M(U)_1\]
such that
\[ \frac{\partial x^k}{\partial x^i} = \delta_i^k\, ,\quad 
\frac{\partial \xi^l}{\partial x^i} = 0\, ,\quad 
\frac{\partial x^k}{\partial \xi^j} = 0\, ,\quad 
\frac{\partial \xi^l}{\partial \xi^j} = \delta_j^l\, .\] 
(Here, of course, $\delta_i^k \in {\cal A}(U)$ is the unit element in 
the algebra ${\cal A}(U)$ if $i=k$ and is zero otherwise.) Moreover,
the vector fields $(\partial /\partial x^i, \partial /\partial \xi^j)$
freely generate ${\cal T}_M(U) \cong {\rm Der}\, {\cal A}(U)$ over
${\cal A}(U)$. This shows that ${\cal T}_M$ is a locally 
free sheaf of rank $m|n = \dim M$ over $\cal A$. It is 
also a sheaf of super Lie algebras; simply because 
${\rm Der}\, {\cal A}(U)$ is a subalgebra of the super Lie algebra
${\rm End}_{\mbox{\R}}\, {\cal A}(U)$ for all open $U \subset M_0$.  
As in the case of ordinary manifolds, given a vector field $X$
there exists a unique derivation ${\cal L}_X$ of the full tensor 
superalgebra  $\otimes_{\cal A} \langle {\cal T}_M, {\cal T}_M^{\ast}\rangle$
over ${\cal T}_M$ compatible with contractions 
such that
\[ {\cal L}_Xf = X(f) \quad \mbox{and} \quad {\cal L}_XY = [X,Y]\]
for all functions $f$ and vector fields $Y$ on $M$ ($[X,Y]$ is the 
{\em super}commutator of vector fields). 
\begin{dof}
Let $M = (M_0,{\cal A})$ be a supermanifold. 
A {\bf tangent vector} to $M$ at $p\in M_0$ is an $\mbox{\R}$-linear map 
$v = v_0 + v_1: {\cal A}_p \rightarrow \mbox{\R}$ such that
\[ v_{\alpha}(fg) = v_{\alpha}(f){\epsilon}^{\ast}_p(g) + 
{(-1)}^{\alpha \tilde{f}}
{\epsilon}^{\ast}_p(f) v_{\alpha}(g)\, , \quad \alpha = 0,1\, ,  \]
for all germs of functions $f,g\in {\cal A}_p$ of pure degree. 
(${\cal A}_p$ denotes the stalk of $\cal A$ at $p$.)
The $\mbox{\Z}_2$-graded vector space of all tangent vectors to $M$
at $p$ is denoted by $T_pM = (T_pM)_0 + (T_pM)_1$ and is called 
the {\bf tangent space} to $M$ at $p$. 
\end{dof}

Let $X$ be a vector field defined on some open set $U\subset M_0$ and 
$p\in U$. Then we can define the {\bf value} $X(p) \in T_pM$ of 
$X$ {\bf at} $p$: 
\[ X(p) (f) := \epsilon^{\ast}_p(X(f))\, ,\quad f\in {\cal A}_p\, .\] 
However, unless $\dim M = m|n = m|0$, a vector field is not determined by
its values at all $p\in M_0$. The above definition of value at a point $p$
is naturally extended to arbitrary tensor fields $S$; the value of $S$ 
at $p$ is denoted by $S(p)$. Again, unless $\dim M = m|n = m|0$, a 
tensor field on $M$ is not determined by
its values at all $p\in M_0$.

Given a morphism $\varphi : M\rightarrow N$, to any local
vector field $X\in {\cal T}_M(U)$ on $M$ we can associate
a vector field $d\varphi X\in (\varphi^{\ast}{\cal T}_N)(U)$
on $N$ with values in ${\cal A}_M$ which is defined by:
\[ (d\varphi X)(f) = X(\varphi^{\ast}f)\, ,\quad 
f\in {\cal A}_N(V)\, ,\]
where $V\subset N_0$ is an open set such that
$\varphi^{-1}_0(V) \supset U$. We recall that 
$\varphi^{\ast}{\cal T}_N$ is the sheaf of 
${\cal A}_M$-modules over $M_0$ defined by:
\[ \varphi^{\ast}{\cal T}_N := {\cal A}_M \otimes_{\varphi_0^{-1}{\cal A}_N}
\varphi_0^{-1}{\cal T}_N \, .\] 
Here the action of $\varphi_0^{-1}{\cal A}_N$ on ${\cal A}_M$ is defined by the
map \[ \varphi_0^{-1}{\cal A}_N \rightarrow \varphi_0^{-1}\varphi_{\ast}
{\cal A}_M \rightarrow {\cal A}_M\]
induced by $\varphi^{\ast}: {\cal A}_N \rightarrow  \varphi_{\ast}
{\cal A}_M$. By the above construction, we obtain a section
$d\varphi$ of the sheaf ${\rm Hom}_{\cal A}({\cal T}_M, \varphi^{\ast}
{\cal T}_N)$, which is expressed with respect to local coordinates
$(u^1,\ldots ,u^{m+n}) = (x^1,\ldots , x^m, \xi^1,\ldots , \xi^n)$    on $M$
and $(v^1,\ldots ,v^{p+q}) = (y^1,\ldots , y^p, \eta^1,\ldots ,$ $\eta^q)$ 
on $N$ by the Jacobian matrix 
$(\frac{\partial \varphi^{\ast}v^i}{\partial u^j})$. 
The value $d\varphi (p) \in {\rm Hom}(T_pM,T_{\varphi_0(p)}N)$ of the 
differential at a point $p\in M_0$ is defined by: 
\[ d\varphi (p)X(p) = (d\varphi X)(p) \]
for all vector fields $X$ on $M$. The differential $d\epsilon$  of 
the canonical morphism $\epsilon :M_0 \rightarrow M$ provides the  
canonical isomorphism $d\epsilon (p): T_pM_0 \stackrel{\sim}{\rightarrow}  
(T_pM)_0$ for all $p\in M$. 
\begin{dof}
A morphism $\varphi : M\rightarrow N$ is called an {\bf immersion}
(respectively, a {\bf submersion}) if $d\varphi$ has constant rank 
$m|n = \dim M$ (respectively, $p|q = \dim N$), i.e.\ 
if for all local coordinates as above the matrix 
$(\epsilon^{\ast}(\partial \varphi^{\ast}y^i/\partial x^j))$ has
constant rank $m$ (respectively, $p$) and 
$(\epsilon^{\ast}(\partial \varphi^{\ast}\eta^i/\partial \xi^j))$ has
constant rank $n$ (respectively, $q$). An immersion $\varphi :M \rightarrow N$ 
is called
{\bf injective} (respectively, an {\bf embedding}, a {\bf closed embedding})
and is denoted by $\varphi : M\hookrightarrow N$ if $\varphi_0 : M_0 
\rightarrow N_0$ is 
injective (respectively, an embedding,   
a closed embedding).  
Two immersions $\varphi : M\rightarrow N$ and ${\varphi}' : M'\rightarrow N$
are called {\bf equivalent} if there exists an isomorphism 
$\psi : M\rightarrow M'$ such that $\varphi = {\varphi}'\circ \psi$. A {\bf 
submanifold} (respectively, an {\bf embedded submanifold}, a {\bf closed
submanifold}) is an equivalence class of injective immersions
(respectively,  embeddings, closed embeddings).  
\end{dof}
Notice that for any supermanifold the canonical morphism $\epsilon 
= ({\rm Id}_{M_0}, \epsilon^{\ast}): M_0 \hookrightarrow M$ is a closed
embedding. 

As for ordinary manifolds, immersions and submersions admit adapted
coordinates:
\begin{prop}\label{L-KProp}
(see \cite{L}, \cite{K}) Let $M$ and $N$ be supermanifolds of dimension 
$\dim M = m|n$ and $\dim N = p|q$. 
A morphism $\varphi :M \rightarrow N$
is an immersion (respectively, submersion)
if and only if for all $(x,y = \varphi_0(x)) \in M_0\times \varphi_0(M_0)
\subset M_0 \times N_0$ there exists local coordinates 
$(x^1, \ldots, x^m, \xi^1, \ldots , \xi^n)$ 
for $M$ defined near $x$ and $(y^1,\ldots , y^p, \eta^1,\ldots , \eta^q)$ 
for $N$ defined near 
$y$ such that 
\[  \varphi^{\ast}y^i =\left\{ \begin{array}{r@{\quad \mbox{for} \quad}l}
x^i & i = 1,\ldots , m\le p\\
0 & m < i\le p
\end{array}\right. 
\quad \mbox{and} \quad
\varphi^{\ast}\eta^j = \left\{ \begin{array}{r@{\quad \mbox{for} \quad}l}
\xi^j & j = 1,\ldots , n\le q\\
0 & n < j\le q
\end{array}\right. \] 
(respectively, 
\[  \varphi^{\ast}y^i = x^i \quad \mbox{for} \quad i = 1,\ldots , 
p\le m\quad \mbox{and} \quad \varphi^{\ast}\eta^j 
= \xi^j\quad \mbox{for} \quad  j = 1,\ldots , q\le n)\, . \]
\end{prop}
{}For any submanifold $\varphi : M\hookrightarrow N$ we define its
{\bf vanishing ideal}  ${\cal J}_y 
\subset ({\cal A}_N)_y$ {\bf at} $y\in N_0$ as follows:  
${\cal J}_y := ({\cal A}_N)_y$ if $y\not\in \varphi_0(M_0)$ and ${\cal J}_y := 
{\rm ker}(r\circ \varphi^{\ast})$ if $y = \varphi_0(x)$, $x\in M_0$, 
where $r : ((\varphi_0)_{\ast}{\cal A}_M)_y \rightarrow ({\cal A}_M)_x$ is the 
natural restriction homomorphism.  The union ${\cal J} := \cup_{y\in N_0}
 {\cal J}_y \subset {\cal A}_N$ is called the {\bf vanishing ideal} of the
submanifold $\varphi : M\hookrightarrow N$. 
By Prop.\ \ref{L-KProp}, it  has the following
property (P): 
{\it If $y \in \varphi_0(M_0)$ then there exists 
$p-m$  even  functions  $(y^i)$ and $q-n$ odd functions $\eta^j$ on $N$ 
vanishing at $y$ whose germs at $y$ generate ${\cal J}_y$  and which 
can be complemented to local coordinates for $N$ 
defined near $y$. (Here $m|n = \dim M$ and $p|q =\dim N$.)
We denote by $Z_U \subset N_0$ the submanifold defined by  
the equations $\epsilon^{\ast}y^i =  0$ on some sufficiently small 
open neighborhood $U\subset N_0$  of $y$. Then $U$ can be chosen 
such that the germs of the functions $(y^i,\eta^j)$ at 
$z$ generate ${\cal J}_z$ for all $z\in Z_U$.}    

Conversely, let $\varphi_0 : M_0\hookrightarrow N_0$ be any  
submanifold and ${\cal J}_y \subset ({\cal A}_N)_y$, $y\in N$,  
a collection of ideals with 
the above property (P) and such that 
$J_y = ({\cal A}_N)_y$ for $y \not\in \varphi_0(M_0)$, then there exists a 
supermanifold $M = (M_0, {\cal A}_M)$ and an injective immersion  
$\varphi = (\varphi_0, \varphi^{\ast}): M \hookrightarrow N$ with 
${\cal J} = \cup _{y\in N}
 {\cal J}_y$ as its vanishing ideal. 

Notice that if $M\hookrightarrow N$ 
is a {\em closed} submanifold then its vanishing ideal can be defined
directly as the sheaf of ideals ${\cal J} := {\rm ker}\, \varphi^{\ast}$.

\subsection{Pseudo-Riemannian metrics, connections and quaternionic
structures on supermanifolds} \label{supermetSec}
Let $\bf A$ be a superalgebra and $T$ a free $\bf A$-module of 
rank $m|n$. 
\begin{dof}
\label{superbilDef} An even (respectively, odd) {\bf bilinear form} 
on $T$ is a biadditive map $g:T\times T \rightarrow A$
such that
\[ g(aX,bY) = (-1)^{\tilde{b}\tilde{X}}abg(X,Y) \]
\[ (\mbox{respectively,} \quad 
g(aX,bY) = (-1)^{\tilde{b}\tilde{X}+\tilde{a} + \tilde{b}}abg(X,Y))\, , \]
for all homogeneous $a,b\in {\bf A}$ and $X,Y \in T$. A bilinear form
$g$ on $T$ is called {\bf supersymmetric} 
(respectively, {\bf super skew symmetric}) if 
\[ g(X,Y) = (-1)^{\tilde{X}\tilde{Y}}g(Y,X)\]
\[ (\mbox{respectively,} \quad 
g(X,Y) = -(-1)^{\tilde{X}\tilde{Y}}g(Y,X))\, ,\]
for all homogeneous $X,Y \in T$. It is called {\bf nondegenerate} if
\[ T\ni X \mapsto g(X,\cdot ) \in T^{\ast} = {\rm Hom}_A(T,{\bf A})
\]
is an isomorphism of $\bf A$-modules. The $\bf A$-module of bilinear forms
on $T$ is denoted by ${\rm Bil}_{\bf A}(T) = {\rm Bil}_{\bf A}(T)_0 + 
{\rm Bil}_{\bf A}(T)_1$. 
\end{dof}
Notice that ${\rm Bil}_{\bf A}(T) \cong {\rm Hom}_{\bf A}(T,T^{\ast}) 
\cong T^{\ast}\otimes_{\bf A} T^{\ast}$, where $T^{\ast} = 
{\rm Hom}_{\bf A}(T,{\bf A})$. 

{}For a supermanifold $M = (M_0,{\cal A})$ the sheaf 
${\rm Bil}_{\cal A}\, {\cal T}_M$ of bilinear forms on ${\cal T}_M$ is defined 
in the obvious way such that $({\rm Bil}_{\cal A}\, {\cal T}_M)(U) = 
{\rm Bil}_{{\cal A}(U)}({\cal T}_M(U))$ 
for every open subset $U\subset M_0$ with the
property that ${\cal T}_M(U)$ is a free ${\cal A}(U)$-module. 
We have obvious isomorphisms of sheaves of $\cal A$-modules:  
${\rm Bil}_{\cal A} \, {\cal T}_M 
\cong  {\rm Hom}_{\cal A}({\cal T}_M,{\cal T}_M^{\ast}) \cong
{\cal T}^{\ast}_M \otimes_{\cal A}  {\cal T}^{\ast}_M$. 
Since any section $g$ of ${\rm Bil}_{\cal A} \, {\cal T}_M$ 
can be considered as a tensor field on $M$, it has a well defined value 
$g(p) \in {\rm Bil}_{\mbox{\R}}(T_pM)$ for all $p\in M_0$. The 
restriction $g(p)|(T_pM)_0\times (T_pM)_0$  defines 
a section $g_0$ of ${\rm Bil}_{{\cal C}_{M_0}^{\infty}}\, {\cal T}_{M_0}$
via the canonical identification  $d\epsilon (p): T_pM_0 
\stackrel{\sim}{\rightarrow}(T_pM)_0$.  
\begin{dof}\label{superRiemDef} 
A {\bf pseudo-Riemannian metric} on a supermanifold $M = (M_0, {\cal A})$,
$M_0$ connected, is an even nondegenerate supersymmetric section $g$ of
${\rm Bil}_{\cal A}\, {\cal T}_M$. The {\bf signature} $(k,l)$ of $g$ is the
signature of the pseudo-Riemannian metric $g_0$ on $M_0$. The 
pseudo-Riemannian metric $g$ is said to be a 
{\bf Riemannian metric}  if $g_0$ is
Riemannian. Let $M = (M_,{\cal A})$ be a supermanifold and $\cal E$ 
a locally free sheaf of ${\cal A}$-modules. A {\bf  connection} on
$\cal E$ is an even section $\nabla$ of the sheaf
${\rm Hom}_{\cal A}({\cal T}_M, {\rm End}_{\mbox{\R}}\, {\cal E})$, which 
to any vector field $X$ on $M$  associates a section 
$\nabla_X$ of ${\rm End}_{\mbox{\R}}\, {\cal E}$ such that
\[ \nabla_Xfs = X(f) s + (-1)^{\tilde{X}\tilde{f}}f\nabla_Xs\]
for all vector fields  $X$ on $M$, functions $f$ on $M$ and 
sections $s$ of $\cal E$ of pure degree. 
The {\bf curvature} $R$ of $\nabla$ is the even super skew symmetric
section of ${\rm End}_{\cal A}\, {\cal E} \otimes_{\cal A} 
{\rm Bil}_{\cal A}\, {\cal T}_M$ defined by:  
\[ R(X,Y) := \nabla_X\nabla_Y - (-1)^{\tilde{X}\tilde{Y}}\nabla_Y\nabla_X
-\nabla_{[X,Y]}\] 
for all vector fields $X$, $Y$ on $M$ of pure degree. 
A {\bf connection on a supermanifold} $M$ is by definition a 
connection on its tangent sheaf ${\cal T}_M$. Its {\bf torsion} is
the even super skew symmetric section of 
${\cal T}_M \otimes_{\cal A} {\rm Bil}_{\cal A}\, {\cal T}_M$ 
defined by: 
\[  T(X,Y) := \nabla_XY - (-1)^{\tilde{X}\tilde{Y}}\nabla_YX - [X,Y]\]
for all vector fields $X$, $Y$ on $M$ of pure degree. 
\end{dof}
As for  ordinary manifolds, a connection $\nabla$ on a supermanifold $M$
induces a connection on the full tensor superalgebra $\otimes_{\cal A} \langle 
{\cal T}_M, {\cal T}_M^{\ast}\rangle$ . 
In particular, if $g$ is e.g.\ an even section of ${\rm Bil}_{\cal A}\,
{\cal T}_M \cong {\cal T}_M^{\ast} \otimes {\cal T}_M^{\ast}$ then we have
\[ (\nabla_Xg)(Y,Z) = Xg(Y,Z) - g(\nabla_XY,Z) -
(-1)^{\tilde{X}\tilde{Y}}g(Y,\nabla_XZ) \]
for all vector fields $X$, $Y$ and $Z$ on $M$ of pure degree. 
As in the case of ordinary manifolds, see e.g.\  
\cite{O'N}, one can prove that a pseudo-Riemannian supermanifold $(M,g)$ has a 
unique torsionfree connection $\nabla = \nabla^g$ such that $\nabla g = 0$. 
We will call this connection the {\bf Levi-Civita connection} of $(M,g)$.
It is computable from the following superversion of the Koszul-formula:
\begin{eqnarray}  2g(\nabla_XY,Z) 
&=& Xg(Y,Z) + (-1)^{\tilde{X}(\tilde{Y} +\tilde{Z})}Yg(Z,X)\nonumber \\ 
& & - \, (-1)^{\tilde{Z}(\tilde{X} + \tilde{Y})}Zg(X,Y)  
- g(X,[Y,Z]) \nonumber \\
& &+\, (-1)^{\tilde{X}(\tilde{Y} + \tilde{Z})}g(Y,[Z,X])  
+ (-1)^{\tilde{Z}(\tilde{X} + \tilde{Y})}g(Z,[X,Y])
\label{superKoszulIEqu}
\end{eqnarray}
for all vector fields $X$, $Y$ and $Z$ on $M$ of pure degree. 

Next we are going to define quaternionic supermanifolds and quaternionic 
K\"ahler supermanifolds. First we define the notion of almost quaternionic
structure.  

\begin{dof}\label{superalmostqstrDef} 
Let $M= (M_0,{\cal A})$ be a supermanifold. An {\bf almost complex structure}
on $M$ is an even global section $J \in ({\rm End}_{\cal A}\, {\cal T}_M)
(M_0)$ such that $J^2 = - {\rm Id}$. An {\bf almost hypercomplex} structure
on $M$ is a triple $(J_{\alpha}) = (J_1,J_2,J_3)$ of almost complex 
structures on $M$ satisfying $J_1J_2 = J_3$. An {\bf almost quaternionic
structure} on $M$ is a subsheaf $Q \subset {\rm End}_{\cal A}\, {\cal T}_M$
with the following property: for every $p\in M_0$ there exists an open 
neighborhood $U \subset M_0$ and an almost hypercomplex structure
$(J_{\alpha})$ on ${M|}_U$ such that $Q(U)$ is a free ${\cal A}(U)$-module
of rank $3|0$ with basis $(J_1,J_2,J_3)$. A pair $(M,J)$ 
(respectively, $(M,(J_{\alpha}))$, $(M,Q))$ as above is called an
{\bf almost complex supermanifold} (respectively, {\bf almost
hypercomplex supermanifold}, {\bf almost quaternionic supermanifold}). 
\end{dof}

Second we introduce the basic compatibility conditions between 
almost quaternionic structures, connections and  pseudo-Riemannian
metrics. 
\begin{dof}
Let $(M,Q)$ be an almost quaternionic supermanifold of dimension $\dim M =
m|n$. A connection $\nabla$ on $(M,Q)$ is called an {\bf almost quaternionic
connection} if $\nabla$ preserves $Q$, i.e.\ if $\nabla_XS$ is a section of 
$Q$ for any vector field $X$ on $M$ and any section $S$ of $Q$. 
A {\bf quaternionic connection} on $(M,Q)$ is a torsionfree almost 
quaternionic connection. If the almost quaternionic structure
$Q$ on $M$ admits a quaternionic connection, then it is called
{\bf 1-integrable} or {\bf quaternionic structure}. In this case the pair
$(M,Q)$ is called a {\bf quaternionic supermanifold}, provided that
$m = \dim M_0 >4$. 
\end{dof}

\begin{dof}
A (pseudo-) Riemannian metric $g$ on an almost quaternionic supermanifold
$(M,Q)$ is called {\bf Hermitian} if sections of $Q$ are $g$-skew 
symmetric, i.e.\ if $g(SX,Y) = -g(X,SY)$ for all 
sections $S$ of $Q$ and vector fields $X$, $Y$ on $M$. In this case
the triple $(M,Q,g)$ is called an {\bf almost quaternionic} (pseudo-)
{\bf Hermitian supermanifold}. If, moreover, the Levi-Civita connection
$\nabla^g$ of the $Q$-Hermitian metric $g$ is quaternionic, then
$(M,Q,g)$ is  called a {\bf  quaternionic} (pseudo-) {\bf K\"ahler
supermanifold}, provided that $m = \dim M_0 >4$. 
\end{dof}
We recall that to any (pseudo-) Riemannian
metric $g$ on a supermanifold $M$ we have associated the (pseudo-) Riemannian
metric $g_0$ on the manifold $M_0$, see Def.\ \ref{superRiemDef}.  
Now we will associate 
an almost quaternionic structure $Q_0$ on $M_0$ to any almost 
quaternionic structure $Q$ on $M$. To any even section $S$ of 
${\rm End}_{\cal A}\, {\cal T}_M$ we associate a section $S_0$
of ${\rm End}_{{\cal C}_{M_0}^{\infty}}\, {\cal T}_{M_0}$
defined by $S_0(p) = S(p)|T_pM_0$, where $S(p) \in 
({\rm End}_{\mbox{\R}}(T_pM))_0$ is the value of the tensor field
$S$ at $p\in M_0$ and, as usual, $T_pM_0$ is canonically identified with
$(T_pM)_0$. Let $Q_0\subset 
{\rm End}_{{\cal C}_{M_0}^{\infty}}\, {\cal T}_{M_0}$ be the subsheaf of 
${\cal C}_{M_0}^{\infty}$-modules spanned by the local sections
of the form $S_0$, where $S$ is a local section of $Q$. Then $Q_0$ is a 
locally free sheaf of rank $3$ and hence it defines a rank $3$ subbundle
of the vectorbundle ${\rm End}\, TM_0$, more precisely, an almost
quaternionic structure on $M_0$ in the usual sense, see Def.~\ref{acsDef}.  

Finally we give the definition of quaternionic supermanifolds and of 
quaternionic K\"ahler supermanifolds of dimension $4|n$. 
\begin{dof}
An almost quaternionic supermanifold $(M,Q)$ of dimension 
$\dim M = 4|n$ is called a {\bf quaternionic supermanifold}
if 
\begin{enumerate}
\item[i)] there exists a quaternionic connection $\nabla$ on
$(M,Q)$, 
\item[ii)] $(M_0,Q_0)$ is a quaternionic manifold in the sense
of Def.\ \ref{quatmfDef}.
\end{enumerate}
An almost quaternionic Hermitian supermanifold $(M,Q,g)$ of 
dimension $\dim M =4|n$ is called a {\bf quaternionic K\"ahler
supermanifold} if 
\begin{enumerate}
\item[i)] the Levi-Civita connection $\nabla^g$ preserves $g$ and 
\item[ii)] $(M_0,Q_0)$ is a quaternionic K\"ahler manifold in the sense
of Def.\ \ref{KaehlerDef}.
\end{enumerate}
\end{dof}

\subsection{Supergroups} \label{supergroupSec}
Let $V = V_0 + V_1$ be a  $\mbox{\Z}_2$-graded vector space. 
We recall that ${\fr gl}(V)$ denotes the general linear super Lie algebra.
Its super Lie bracket is the supercommutator on ${\rm End}_{\mbox{\R}}\, V$, 
see Example 2. 

\begin{dof}
\label{superadjDef} A {\bf representation} of a super Lie algebra
$\fr g$ on a $\mbox{\Z}_2$-graded vector space $V$ is a homomorphism 
of super Lie algebras ${\fr g} \rightarrow {\fr gl}(V)$. 
The {\bf adjoint representation} of $\fr g$ is the 
representation ${\rm ad}: {\fr g} \ni x \mapsto {\rm ad}_x \in 
{\fr gl}({\fr g})$ defined by:
\[ {\rm ad}_xy = [x,y]\, ,\quad x,y \in {\fr g}\, .\] 
Here $[\cdot ,\cdot ]$ denotes the super Lie bracket in $\fr g$. 
A {\bf linear super Lie algebra} is a subalgebra ${\fr g} 
\subset {\fr gl}(V)$. 
\end{dof}
Notice that, by the super Jacobi identity, see Def.\ \ref{slbDef}, 
${\rm ad}_x$ is a derivation of 
$\fr g$ for all $x\in {\fr g}$.  

Now let $E$ be a module over a superalgebra $\bf A$; say 
$E = V \otimes {\bf A}$, where $V$ is a $\mbox{\Z}_2$-graded vector space
(and $\otimes$ stands for the $\mbox{\Z}_2$-graded tensor product over
$\mbox{\R}$).    The invertible elements of $({\rm End}_{\bf A}\, E)_0$
form a group, which is denoted by ${\rm GL}_{\bf A}(E)$. Given a 
$\mbox{\Z}_2$-graded vector space $V$, the correspondence
\[ {\bf A} \mapsto {\rm GL}_{\bf A}(V \otimes {\bf A})\]
defines a covariant functor from the category of 
superalgebras into the category of groups. We can compose this functor with the
contravariant functor
\[ M = (M_0,{\cal A}) \mapsto {\cal A}(M_0)\]
from the category of supermanifolds into that of 
superalgebras. The resulting contravariant functor 
\[ M \mapsto {\rm GL}(V)[M] := {\rm GL}_{{\cal A}(M_0)}(V \otimes 
{\cal A}(M_0))\] 
from the category of supermanifolds into that of groups
is denoted by ${\rm GL}(V)[\cdot ]$. 

\begin{dof}\label{supergroupDef} 
A {\bf supergroup} $G[\cdot ]$ is a contravariant functor 
from the category of supermanifolds into that of groups. The supergroup
${\rm GL}(V)[\cdot ]$, defined above, is called the {\bf general
linear supergroup} over the $\mbox{\Z}_2$-graded vector space $V$. 
A {\bf subgroup} of a supergroup $G[\cdot ]$ is a supergroup
$H[\cdot ]$ such that $H[M]$ is  a subgroup of $G[M]$ 
for all supermanifolds $M$ and $H[\varphi ]  = G[\varphi ]|H[N]$ for all
morphisms $\varphi : M\rightarrow N$. 
We will write $H[\cdot ] \subset G[\cdot ]$ if  
$H[\cdot ]$ is a subgroup of the supergroup $G[\cdot ]$. 
A {\bf linear supergroup}  is a subgroup $H[\cdot ] \subset 
{\rm GL}(V)[\cdot ]$. 
\end{dof}

\noindent 
{\bf Example 7:} To any linear super Lie algebra ${\fr g} \subset
{\fr gl}(V)$ we can associate the linear supergroup 
$G[\cdot ] \subset {\rm GL}(V)[\cdot ]$ defined by:
\[ G[M] := \langle \exp ({\fr g} \otimes {\cal A}(M_0))_0\rangle\]
for any supermanifold $M = (M_0,{\cal A})$. The right-hand
side is the subgroup of \linebreak[3] ${\rm GL}_{{\cal A}(M_0)}$ $(V \otimes 
{\cal A}(M_0))$ generated by the exponential image of the
Lie algebra $({\fr g} \otimes {\cal A}(M_0))_0 
\subset ({\fr gl}(V) \otimes {\cal A}(M_0))_0 
\cong ({\rm End}_{{\cal A}(M_0)}(V \otimes {\cal A}(M_0)))_0$. 
The convergence of the exponential series, which is locally uniform in all 
derivatives of arbitrary order, follows from the analyticity of the exponential
map by a (finite) Taylor expansion with respect to the odd coordinates. 
$G[\cdot ]$ is called the {\bf linear supergroup associated to the 
linear super Lie algebra} $\fr g$. 

{}For any supermanifold $M$ we denote by $\Delta_M = (\Delta_{M_0},
\Delta_M^{\ast})$ the {\bf diagonal embedding} defined by: 
$\Delta_{M_0}(x) = (x,x)$ ($x\in M_0$) and 
$\Delta_M^{\ast} \pi_i^{\ast}f = f$ for all functions $f$ on $M$. 
Here $\pi_i : M\times M \rightarrow M$ denotes the projection
onto the $i$-th factor of $M\times M$ ($i = 1,2$). 
\begin{dof}\label{superLiegroupDef}
A {\bf Lie supergroup} is a supermanifold $G = (G_0,{\cal A}_G)$ whose 
underlying manifold is a Lie group $G_0$ (which, for convenience, 
we will always assume to be connected) with 
neutral element $e \in G_0$, 
multiplication $\mu_0 : G_0\times G_0 \rightarrow G_0$ and 
inversion $\iota_0 : G_0 \rightarrow G_0$, together with morphisms
$\mu = (\mu_0 ,\mu^{\ast}): G\times G \rightarrow G$ and 
$\iota = (\iota_0 , \iota^{\ast}): G \rightarrow G$ such that
\begin{enumerate}
\item[i)] $\mu \circ ({\rm Id}_G \times \mu ) = \mu \circ  
(\mu \times {\rm Id}_G) \in {\rm Mor}(G\times G \times G, G)$,
\item[ii)] $\mu \circ ({\rm Id}_G \times \epsilon_e) = 
\pi_1: G \times \{ e\} \stackrel{\sim}{\rightarrow} G$, 
$\mu \circ (\epsilon_e \times {\rm Id}_G)
= \pi_2 : \{ e\} \times G \stackrel{\sim}{\rightarrow} G$ and  
\item[iii)] $\mu \circ ({\rm Id}_G \times \iota ) \circ \Delta_G =
\mu \circ (\iota \times {\rm Id}_G) \circ \Delta_G = {\rm Id}_G$. 
\end{enumerate}
Here $\epsilon_e = (e, \epsilon_e^{\ast}): \{ e\} \hookrightarrow  G_0
\stackrel{\epsilon}{\hookrightarrow} G$ is the canonical 
embedding. The morphisms $\mu$ and $\iota$ are called, respectively,  
{\bf multiplication} and {\bf inversion} in the Lie supergroup $G$. 
A {\bf Lie subgroup} of a Lie supergroup $G = (G_0, {\cal A}_G)$  
is a submanifold $\varphi : H = (H_0, {\cal A}_H) 
\hookrightarrow G$ such that 
the immersion $\varphi_0 : H_0 \hookrightarrow G_0$ induces on $H_0$ 
the structure of  Lie subgroup of $G_0$ and $\varphi$ induces on $H$
the structure of Lie supergroup with underlying Lie group $H_0$.
We will write $H \subset G$ if $H$ is a Lie subgroup of $G$. 
An {\bf action} of a Lie supergroup $G$ on a supermanifold $M$
is a morphism $\alpha : G \times M \rightarrow M$ such that
\begin{enumerate}
\item[i)] $\alpha \circ ({\rm Id}_G \times \alpha ) = 
\alpha \circ (\mu \times {\rm Id}_M)$ and  
\item[ii)] $\alpha \circ (\epsilon_e \times {\rm Id}_M) 
= \pi_2 : \{ e\} \times M \stackrel{\sim}{\rightarrow} M$.  
\end{enumerate}
\end{dof}
Given an action $\alpha : G \times M \rightarrow M$ we have the notion
of {\bf fundamental vector field} $X$ associated to $x\in T_eG$. It is
defined by $X(f) := x(\alpha^{\ast}f)$. The correspondence $x\mapsto X$
defines an even $\mbox{\R}$-linear map $T_eG \rightarrow {\cal T}_M(M_0)$ 
from the $\mbox{\Z}_2$-graded vector space $T_eG$ to the free 
${\cal A}_M(M_0)$-module of global vector fields on $M$. 

\begin{dof} \label{superinvDef} 
Let an action $\alpha$ of a Lie supergroup $G$ on a 
supermanifold $M$ be given. 
A tensor field $S$ on $M$ is called {\bf $G$-invariant} if the
Lie derivative ${\cal L}_XS = 0$ for all fundamental vector fields
$X$. (Notice that, in particular, this defines the notion of $G$-invariant
pseudo-Riemannian metric on $M$.) An almost quaternionic structure
$Q$ on $M$ is called {\bf $G$-invariant} if  ${\cal L}_XS$ is a section
of $Q$ for all sections $S$ of $Q$ and fundamental vector fields $X$ on $M$.  
\end{dof}
We can specialize the above definition to the (left-) action
$\mu :G\times G \rightarrow G$ given by the multiplication in the Lie 
supergroup $G$. The $G$-invariant tensor fields with respect
to that action are called {\bf left-invariant}. A tensor field 
$S$ is called {\bf right-invariant} if $L_XS = 0$ for all 
left-invariant vector fields $X$ on $G$. The right-invariant
vector fields are precisely the fundamental vector fields
for the action $\mu$. 

{}For any Lie supergroup $G$, we can define a group homomorphism
$\mu_l : G_0 \rightarrow {\rm Aut}(G)$ by:
\[ \mu_l(g) := \mu \circ (\epsilon_g \times {\rm Id}_G) : 
\{ g\} \times G \cong G \rightarrow G\, ,\quad g\in G_0\, .\]
Here $\epsilon_g: \{ g\} \hookrightarrow G$ is the canonical embedding
and $\{ g\} \times G$ is canonically identified with $G$ via the projection 
$\pi_2 : \{ g\} \times G \stackrel{\sim}{\rightarrow} G$ onto 
the second factor. Similarly, we can define a group antihomomorphism
$\mu_r : G_0 \rightarrow {\rm Aut}(G)$ by:
\[ \mu_r(g) := \mu \circ ({\rm Id}_G \times \epsilon_g): 
G \times \{ g\} \cong G \rightarrow G\, ,\quad g\in G_0\, .\]
Notice that the canonical embedding $\epsilon :G_0 \hookrightarrow G$
is $G_0$-equivariant with respect to the usual left- (respectively, right-) 
action on $G_0$ and the action on $G$ defined by $\mu_l$ (respectively, 
$\mu_r$). 

Given a Lie supergroup $G$ and a supermanifold $M$ there is a natural
group structure on ${\rm Mor}(M,G)$ with multiplication defined by:
\[ \varphi \cdot \psi := \mu \circ (\varphi \times \psi ) \circ \Delta_M
\, ,\quad \varphi , \psi \in {\rm Mor}(M,G)\, .\]
The correspondence $M\mapsto {\rm Mor}(M,G)$ defines a supergroup.

\begin{dof}\label{subordinateDef} 
The supergroup $G[\cdot ] := {\rm Mor}(\cdot , G)$ is called the 
supergroup {\bf subordinate} to the Lie supergroup $G$. 
\end{dof}

\noindent 
{\bf Example 8:} Let $V$ be a  $\mbox{\Z}_2$-graded vector space. Then
${\rm GL}_{\mbox{\R}}(V)$ is, by definition, the group of invertible
elements of $({\rm End}_{\mbox{\R}}\, V)_0$. The Lie group 
${\rm GL}_{\mbox{\R}}(V) \cong {\rm GL}_{\mbox{\R}}(V_0) \times 
{\rm GL}_{\mbox{\R}}(V_1)$ is an open submanifold of the 
vector space $({\rm End}_{\mbox{\R}}\, V)_0 = 
SM({\rm End}_{\mbox{\R}}\, V)_0$, see Example 5. We define the submanifold
\[ {\rm GL}(V) :=  {SM({\rm End}_{\mbox{\R}}\, V)|}_{{\rm GL}_{\mbox{\R}}(V)}
\hookrightarrow SM({\rm End}_{\mbox{\R}}\, V)\, .\]
The manifold underlying 
the supermanifold  ${\rm GL}(V) = ({\rm GL}(V)_0, {\cal A}_{{\rm GL}(V)})$ 
is the above Lie group: ${\rm GL}(V)_0 = {\rm GL}_{\mbox{\R}}(V)$. 
From the definition of the supermanifold ${\rm GL}(V)$ it is clear that 
${\rm Mor}(M,{\rm GL}(V))$ is canonically identified with the set
${\rm GL}_{{\cal A}(M_0)}(V \otimes {\cal A}(M_0)) = {\rm GL}(V)[M]$
(cf.\ Def.\ \ref{supergroupDef}) for any supermanifold 
$M = (M_0,{\cal A})$. Moreover, ${\rm GL}(V)$
has a unique structure of Lie supergroup inducing the canonical
group structure on ${\rm GL}(V)[M]$ for any supermanifold $M$. 
In other words, the general linear supergroup ${\rm GL}(V)[\cdot ]$ is the 
supergroup subordinate to the Lie supergroup ${\rm GL}(V)$, see 
Def.\ \ref{subordinateDef}. 

\begin{dof}
The Lie supergroup ${\rm GL}(V)$
is called the {\bf general linear Lie supergroup}. 
A {\bf linear Lie supergroup} is a Lie subgroup of
${\rm GL}(V)$. 
\end{dof}

\noindent 
There exists a unique morphism
\[ {\rm Exp} : SM({\fr gl}(V)) \rightarrow {\rm GL}(V)\]
such that 
\[ {\rm Exp} \circ \varphi = \exp (\varphi ) \] 
for all supermanifolds $M = (M_0,{\cal A})$ and $\varphi \in 
{\rm Mor}(M,SM({\fr gl}(V))) = ({\fr gl}(V) \otimes {\cal A}(M_0))_0
= ({\rm End}_{{\cal A}(M_0)}\, V\otimes {\cal A}(M_0))_0$, 
where,  on the right-hand side, 
$\exp : ({\rm End}_{{\cal A}(M_0)}\, V\otimes {\cal A}(M_0))_0
\rightarrow {\rm GL}_{{\cal A}(M_0)} (V\otimes {\cal A}(M_0)) =
{\rm Mor}(M,{\rm GL}(V))$ is the exponential map
for even endomorphisms of $V\otimes {\cal A}(M_0)$ 
(for the definition of the supermanifold $SM({\fr gl}(V))$ see Example 5). 
The underlying map ${\rm Exp}_0 : 
{\fr gl}(V)_0 = SM({\fr gl}(V))_0 \rightarrow
{\rm GL}(V)_0$ is the ordinary exponential map for the Lie group
${\rm GL}(V)_0$: ${\rm Exp}_0 = \exp$. The morphism $\rm Exp$ is called the 
{\bf exponential morphism} of ${\fr gl}(V)$. 

\noindent 
{\bf Example 9:} Let ${\fr g} \subset {\fr gl}(V)$ be a linear 
super Lie algebra and $G[\cdot ] \subset {\rm GL}(V)[\cdot ]$
the correponding linear supergroup, see Example 7. We denote by $G_0 
\subset GL(V)_0$ the connected linear Lie group with Lie algebra ${\fr g}_0$. 
It is the immersed Lie subgroup generated by the exponential image of 
${\fr g}_0 = {\fr gl}(V)_0$. We define an  ideal ${\cal J}_g  
\subset ({\cal A}_{{\rm GL}(V)})_g$ as follows: A germ of function 
$f \in ({\cal A}_{{\rm GL}(V)})_g$  ($g\in G_0$) belongs to 
${\cal J}_g$ if and only if $r(\varphi^{\ast}f) = 0$
for all injective immersions 
$\varphi \in G[M] \subset Mor(M,{\rm GL}(V))$, where 
$r : ((\varphi_0)_{\ast}{\cal A}_M)_g \rightarrow 
({\cal A}_M)_{\varphi_0^{-1}(g)}$ is the natural restriction map. 
We claim that ${\cal J}_G := \cup_{g\in G_0} {\cal J}_g$ is the vanishing 
ideal of a submanifold
$G \hookrightarrow {\rm GL}(V)$. Due to the invariance of 
${\cal J}_G$ under the group $\mu_l(G_0) \subset {\rm Aut}(G)$
it is sufficient to prove the claim locally over an open set of the form 
$\exp U \subset G_0$, $U \subset {\fr g}_0$ an open neighborhood of 
$0\in {\fr g}_0$. The local statement follows from the fact that
the  morphism $SM({\fr g})\hookrightarrow SM({\fr gl}(V)) 
\stackrel{{\rm Exp}}{\rightarrow} {\rm GL}(V)$ has maximal rank
at $0\in {\fr g}_0 = SM({\fr g})_0$ and hence defines a submanifold 
${SM({\fr g})|}_U\hookrightarrow  {\rm GL}(V)$
for some open neighborhood of $U$ of $0\in {\fr g}_0$. The vanishing
ideal of this submanifold coincides with ${\cal J}_G$ over
$\exp U$ by the definition of the exponential morphism. 
Next we claim that multiplication $\mu : {\rm GL}(V) \times {\rm GL}(V)
\rightarrow {\rm GL}(V)$ and inversion $\iota :
{\rm GL}(V) \rightarrow {\rm GL}(V)$ have the property that
$\mu^{\ast} {\cal J}_G \subset {\cal J}_{G\times G}$ and 
$\iota^{\ast}{\cal J}_G  = {\cal J}_G$. Here ${\cal J}_{G\times G}$
is the vanishing ideal of the submanifold 
$G\times G \hookrightarrow {\rm GL}(V)
\times {\rm GL}(V)$. The claim follows from the fact that
$G[\cdot ]$ is a subgroup of ${\rm GL}(V)[\cdot ]$ and 
implies that the morphisms 
$G\times G\hookrightarrow {\rm GL}(V)
\times {\rm GL}(V) \stackrel{\mu}{\rightarrow} {\rm GL}(V)$ 
and $G \hookrightarrow {\rm GL}(V) \stackrel{\iota}{\rightarrow} 
 {\rm GL}(V)$ induce morphisms $G \times G \rightarrow G$ and
$G \rightarrow G$, which induce on $G$ the structure of Lie 
supergroup.  In other words, $G$ is a Lie subgroup of the 
general linear Lie supergroup ${\rm GL}(V)$. It is called
the {\bf linear Lie supergroup associated to the linear
Lie superalgebra} ${\fr g} \subset {\fr gl}(V)$. 
Notice that ${\rm Exp}^{\ast}$ maps ${\cal J}_G$ into the 
vanishing ideal of $SM({\fr g})  \hookrightarrow SM({\fr gl}(V))$
and hence the restriction $SM({\fr g}) \hookrightarrow 
SM({\fr gl}(V)) \stackrel{{\rm Exp}}{\rightarrow}
{\rm GL}(V)$ of the exponential morphism induces a morphism
$SM({\fr g})\rightarrow G$, which is called the {\bf exponential morphism}
of $\fr g$ and is again denoted by $\rm Exp$. Its differential at $0
\in {\fr g}_0$ yields an isomorphism ${\fr  g} \cong T_0SM({\fr g}) 
\stackrel{\sim}{\rightarrow} T_eG$.

\subsection{Homogeneous supermanifolds} \label{superhomogSec}
Let $G = (G_0, {\cal A}_G)$ be a Lie supergroup, $K \subset G_0$ a closed
subgroup and $\pi_0 : G_0 \rightarrow G_0/K$ the canonical projection.
Then the subsheaf ${\cal A}_G^K := {\cal A}_G^{\mu_r(K)} 
\subset {\cal A}_G$ of functions on $G$ invariant under the 
subgroup $\mu_r(K) \subset {\rm Aut}(G)$ is again a sheaf of 
superalgebras on $G_0$. Explicitly, a function $f \in {\cal A}_G(U)$
($U\subset G_0$ open) belongs to ${\cal A}_G^K(U)$ if it can be 
extended to a $\mu_r(K)$-invariant function over $UK \subset G_0$. 
Its pushed forward sheaf ${\cal A}_{G/K} := 
(\pi_0)_{\ast}{\cal A}_G^K$ is a sheaf of superalgebras
on the homogeneous manifold $G_0/K$. 

\begin{thm} \label{superhomogThm}
(cf.\ \cite{K}) Let $G \hookrightarrow {\rm GL}(V)$ be a closed linear 
Lie supergroup and $K \subset G_0$ a closed subgroup. Then
$G/K = (G_0/K, {\cal A}_{G/K})$ is a supermanifold with a canonical
submersion $\pi : G \rightarrow G/K$ and a canonical action 
$\alpha : G \times G/K \rightarrow G/K$. 
\end{thm}

\noindent
{\bf Proof:} Since $\mu_l(g)\in {\rm Aut}(G)$  induces isomorphisms
\[ {\cal A}_{G/K}(gU) \stackrel{\sim}{\rightarrow} {\cal A}_{G/K}(U)\] 
for all $g \in G_0$ and $U \subset G_0/K$ open, it is sufficient to check 
that  ${G/K|}_U = (U, {{\cal A}_{G/K}|}_U)$ is a supermanifold
for some open neighborhood $U$ of $eK \in G_0/K$.

\begin{lemma}\label{GL/KLemma}
Under the assumptions of Thm.\ \ref{superhomogThm}, there exists
local coordinates $(x,y,\xi )$ for ${\rm GL}(V)$ over some 
neighborhood $U = UK \subset {\rm GL}(V)_0$ of $e\in {\rm GL}(V)_0$
such that
\begin{enumerate}
\item[1)] $x = (x^i)$ and $y = (y^j)$ consist of even functions and 
$\xi = (\xi^k)$ of odd functions,
\item[2)] ${\cal A}_{{\rm GL}(V)}^K(U)$ is the subalgebra of 
${\cal A}_{{\rm GL}(V)}(U)$ which consists of functions $f(x,\xi )$
only of $(x,\xi )$, more precisely,
\[ f(x,\xi ) = \sum_{\alpha } f_{\alpha}(x)\xi^{\alpha} \, ,\]
where the $f_{\alpha}(x)\in {\cal A}_{{\rm GL}(V)}(U)$ are functions  
only of $x$, i.e.\ if $f_{\alpha}(x) \neq 0$ then 
$f_{\alpha}(x)$ does not belong to the ideal generated by $\xi$ and 
$\epsilon^{\ast}f_{\alpha}(x) \in C^{\infty}(U)$ are functions only of
$\epsilon^{\ast}x$ (independent of $\epsilon^{\ast}y$). 
Here we used the multiindex notation $\alpha = 
(\alpha_1, \ldots , \alpha_{2mn}) \in   \mbox{\Z}_2^{2mn}$ ($m|n = \dim V$)
and $\xi^{\alpha} = \prod_{j=1}^{2mn} (\xi^j)^{{\alpha}_j}$. 
\end{enumerate}
\end{lemma}

\noindent
{\bf Proof:} The natural (global) coordinates on the supermanifold 
${\rm GL}(V)$ are the matrix coefficients with respect to some basis of 
$V$. We will denote them simply by $(z^i,\zeta^j)$ instead of using
matrix notation. For these coordinates it is clear that 
$\mu_r(g)^{\ast}z^k$ is a linear combination 
(over the real numbers) of the even 
coordinates $z := (z^i)$ and  $\mu_r(g)^{\ast}\zeta^l$ is a 
linear combination (over the real numbers) 
of the odd coordinates $\zeta := (\zeta^j)$
for all $g\in {\rm GL}(V)_0$.  
In particular, we obtain a representation $\rho$ of $K \subset {\rm GL}(V)_0$
on the vector space spanned by the odd coordinates. Let $E_{\rho} 
\rightarrow G_0/K$ denote the vector bundle associated
to this representation. Any $\mu_r(K)$-invariant function 
on ${\rm GL}(V)$ linear in $\zeta$ defines a section of the
dual vector bundle $E_{\rho}^{\ast}$ and vice versa. Now, since
$E_{\rho}^{\ast}$ is locally trivial (like any vector bundle), we
can find $\mu_r(K)$-invariant local functions $\xi = (\xi^j)$ linear
in $\zeta$ such that $(z,\xi)$ are local coordinates for ${\rm GL}(V)$
over some open neighborhood $U = UK \subset {\rm GL}(V)_0$ of $e$. 
Next, by way of a local diffeomorphism in the even coordinates $z$,  
we can arrange that $z = (x,y)$, where the $x$ are $\mu_r(K)$-invariant
functions on ${\rm GL}(V)$ such that $\epsilon^{\ast}x \in C^{\infty}(U)$
induce local coordinates on $G_0/K$. Now any function $f \in
{\cal A}_{{\rm GL}(V)}(U)$ has a unique expression of the form 
\begin{equation}
\label{expansionEqu} \sum_{\alpha} f_{\alpha}(x,y)\xi^{\alpha} \, ,\
\end{equation}
where the $f_{\alpha}(x,y)$ are functions only of $(x,y)$. From the
$\mu_r(K)$-invariance of the $\xi^{\alpha}$ it follows that $f$ is 
$\mu_r(K)$-invariant if and only if the $f_{\alpha}(x,y)$ are
$\mu_r(K)$-invariant. The function $f_{\alpha}(x,y)$ is 
$\mu_r(K)$-invariant if and only if $\epsilon^{\ast}f_{\alpha}(x,y)
\in C^{\infty}(U)$ is invariant under the right-action of 
$K$ on ${\rm GL}_0$, i.e.\ if and only if $\epsilon^{\ast}f_{\alpha}(x,y)$
is a function only of $\epsilon^{\ast}x$. This shows that
$f \in {\cal A}_{{\rm GL}(V)}(U)^K$ if and only if the
coefficients $f_{\alpha}(x,y)$ in the expansion (\ref{expansionEqu})
are functions only of $x$ $\Box$ 

We continue the proof of Thm.\ \ref{superhomogThm}. From Lemma
\ref{GL/KLemma} it follows that ${\rm GL}(V)/K$ is a supermanifold. 
In fact, the $K$-invariant local functions $(x,\xi )$ on ${\rm GL}(V)$
constructed in that lemma induce local coordinates on 
${\rm GL}(V)/K$. Next we restrict the coordinates $(x,y,\xi )$ 
to the submanifold $G\hookrightarrow {\rm GL}(V)$. We can decompose
$x = (x',x'')$, $\xi = (\xi',\xi'')$ such that
$(x',y,\xi')$ restrict to local  coordinates on $G$ over some
open neighborhood (again denoted by) 
$U$ of $e \in G_0$ (notice that, by construction, $y$ restrict
to local coordinates on $K$).  Now,  as in the proof of the 
corresponding statement for ${\rm GL}(V)$ (see Lemma \ref{GL/KLemma}),
it follows that ${\cal A}_G^K(U)$ consists precisely of all functions
of the form $f(x',\xi') = \sum f_{\alpha}(x') (\xi')^{\alpha}$. 
This proves that $G/K$ is a supermanifold with local
coordinates $(x',\xi')$ over $U$. 

The inclusion ${\cal A}_G^K \subset {\cal A}_G$ defines the canonical
submersion $\pi : G \rightarrow G/K$. Finally, if $f$ is a $\mu_r(K)$-invariant
(local) function on $G$ then $\mu^{\ast}f$ is a 
(local) function on $G\times G$ invariant under the group
${\rm Id}_G \times \mu_r(K) \subset {\rm Aut}(G\times G)$.
This shows that the composition $G\times G \stackrel{\mu}{\rightarrow}
G \stackrel{\pi}{\rightarrow} G/K$ factorizes to a morphism $\alpha :
G \times G/K \rightarrow G/K$, which defines an action of $G$ on $G/K$.
$\Box$   

\begin{dof}
The supermanifold $M=G/K$ is called the {\bf homogeneous supermanifold}
associated to the pair $(G,K)$. 
\end{dof}

{}For the rest of the paper let ${\fr g} \subset {\fr gl}(V)$ be a linear
super Lie algebra, ${\fr k} \subset {\fr g}_0$ a subalgebra and 
${\fr g} = {\fr k} + {\fr m}$ a $\fr k$-invariant direct decomposition
compatible with the $\mbox{\Z}_2$-grading. We denote by 
$K \subset G_0 \subset G \subset {\rm GL}(V)$ the corresponding 
linear Lie supergroups (see Example 9) and assume that the (connected)  
Lie subgroups $K \subset G_0 \subset {\rm GL}(V)_0$ are closed. 
Then, by  Thm.\ \ref{superhomogThm}, $M = G/K$ is a supermanifold
with a canonical action of $G$. We have   the canonical identification 
${\fr m} \cong {\fr g}/{\fr k} \cong T_{eK}M$ given by $x\mapsto X(eK)$, 
where $X(eK)$ is the value of the fundamental vector field
$X$ on $M$ associated to $x\in {\fr m}$ at the base point 
$eK \in G_0/K = M_0$. 
We claim that any ${\rm ad}_{{\fr k}}$-invariant tensor $S_{eK}$  
over $\fr m$ defines a correponding $G$-invariant (see Def.\ 
\ref{superinvDef}) tensor field on $M$ 
such that $S(eK) = S_{eK}$. Here by a tensor  
over $\fr m$ we mean an element  of the full tensor superalgebra
$\otimes\langle {\fr m}, {\fr m}^{\ast}\rangle$ over $\fr m$. 
In fact, for any tensor $S_{eK}$ over $\fr m$ there exists
a corresponding left-invariant tensor field $S$ on $G$ such that
$S(eK) = S_{eK}$. In order for $S$ to define a tensor field on 
$G/K$ it is necessary and sufficient that $S$ is $\mu_r(K)$-invariant,
or, equivalently, that $S_{eK}$ is ${\rm ad}_{\fr k}$-invariant.
In particular, we have the following proposition: 
\begin{prop}
Let $g_{eK}$ be an ${\rm ad}_{{\fr k}}$-invariant nondegenerate supersymmetric
bilinear form on $\fr m$. Then there exists a unique $G$-invariant 
pseudo-Riemannian metric $g$ on $M =G/K$ (see Def.\ \ref{superRiemDef} 
and Def.\ \ref{superinvDef}) such that $g(eK) = g_{eK}$. 
Let $Q_{eK}$ be an ${\rm ad}_{{\fr k}}$-invariant quaternionic 
structure on $\fr m$  (i.e.\ ${\rm ad}: {\fr k} \rightarrow {\fr gl}({\fr m})$
normalizes $Q_{eK}$). Then there exists a unique $G$-invariant 
almost quaternionic structure $Q$ on $M$ (see Def.\ \ref{superalmostqstrDef}
and Def.\ \ref{superinvDef}) such that $Q(eK) = Q_{eK}$. 
\end{prop}

Finally, we need to discuss $G$-invariant connections on $M =G/K$. 

\begin{dof}
A connection $\nabla$ on a homogeneous supermanifold $M =G/K$ 
is called {\bf $G$-invariant} if 
\[ {\cal L}_X(\nabla_YS) = \nabla_{[X,Y]}S + (-1)^{\tilde{X}\tilde{Y}}
\nabla_Y ({\cal L}_XS)\] 
for all vector fields $X$ and $Y$on $M$. 
\end{dof}
Let $\nabla$ be a connection on a supermanifold $M$. For any vector
field $X$ on $M$ one defines the ${\cal A}_M(M_0)$-linear operator 
\[ L_X := {\cal L}_X -\nabla_X \, . \] 
We denote by $L_X(p) \in {\rm End}_{\mbox{\R}}\, T_pM$ its value at 
$p\in M_0$; it is defined by $L_X(p)Y(p) = (L_XY)(p)$ for all
vector fields $Y$ on $M$. 

{}For a $G$-invariant connection $\nabla$ on a homogeneous supermanifold
$M =G/K$ as above we define the {\bf Nomizu map} $L = L(\nabla ): 
{\fr g} \rightarrow {\rm End}(T_{eK}M)$, $x\mapsto L_x$, by the equation
\[ L_x := L_X(eK)\, ,\] 
where $X$ is the fundamental vector field on $M$ associated to
$x\in {\fr g}$. 
The operators $L_x \in {\rm End}(T_{eK}M)$ will be called 
{\bf Nomizu operators}. They have the following properties:
\begin{equation} \label{superNomizuIEqu} 
L_x = d\rho (x) \quad \mbox{for all} \quad x \in {\fr k}  
\end{equation}
and 
\begin{equation} \label{superNomizuIIEqu} 
L_{{\rm Ad}_kx} = \rho (k) L_x \rho (k)^{-1} \quad \mbox{for all} 
\quad x \in {\fr g}\, , \quad k\in K\, , 
\end{equation}
where $\rho :K \rightarrow {\rm GL}(T_{eK}M)$ is the isotropy 
representation (under the identification
$T_{eK}M \cong {\fr m}$ the representation $\rho$
is identified with adjoint representation of $K$ on $\fr m$). 

Conversely, any even linear map $L : {\fr g} \rightarrow {\rm End}(T_{[e]}M)$ 
satisfying (\ref{superNomizuIEqu}) and (\ref{superNomizuIIEqu})
is the Nomizu map of  a uniquely defined $G$-invariant connection
$\nabla = \nabla (L)$ on $M$. Its torsion tensor $T$ and curvature tensor
$R$  are expressed at $eK$ by: 
\[ T(\pi x, \pi y) = -(L_x \pi y - (-1)^{\tilde{x}\tilde{y}}
L_y \pi x + \pi [x,y]) \]
and 
\[ R(\pi x, \pi y) = [L_x,L_y]  + L_{[x,y]} \, ,\quad x,y\in {\fr g}\]
where $\pi  : {\fr g} \rightarrow T_{eK}M$ is the  canonical projection
$x \mapsto \pi x =  X(eK) = \frac{d}{dt}|_{t=0} \exp (tx)K$. 

Suppose now that we are given a $G$-invariant geometric structure
$S$ on $M$ (e.g.\ a $G$-invariant almost quaternionic structure $Q$)
defined by a corresponding $K$-invariant geometric structure
$S_{eK}$ on $T_{eK}M$. Then a $G$-invariant connection $\nabla$
preserves $S$ if and only if the corresponding Nomizu operators
$L_x$, $x\in {\fr g}$, preserve $S_{eK}$. So to construct a 
$G$-invariant connection preserving $S$ it is sufficient to
find a Nomizu map $L: {\fr g} \rightarrow {\rm End} (T_{eK}M)$
such that $L_x$ preserves $S_{e}$ for all $x\in {\fr g}$. 
We observe that, due to the $K$-invariance of $S_{e}$, the Nomizu
operators $L_x$ preserve $S_{e}$ already for $x\in {\fr k}$. 
The above considerations can be specialized as follows: 
\begin{prop}
Let $Q$ be a $G$-invariant almost quaternionic structure on a
homogeneous supermanifold $M = G/K$. There is a natural one-to-one
correspondence between $G$-invariant almost quaternionic connections on 
$(M,Q)$ and Nomizu maps $L: {\fr g} \rightarrow {\rm End} (T_{eK}M)$,
whose image normalizes $Q(eK)$, i.e.\ whose Nomizu operators $L_x$,
$x\in {\fr g}$, belong to the normalizer ${\rm n}(Q) \cong {\fr sp}(1)
\oplus {\fr gl}(d,\mbox{\Ha})$ ($d = (m+n)/4$) of the 
quaternionic structure $Q(eK)$ in the super Lie algebra ${\fr gl}(T_{eK}M)$. 
\end{prop}
\begin{cor}
\label{superquatCor} Let $(M = G/K,Q)$ be a homogeneous almost quaternionic
su\-per\-ma\-ni\-fold and $L: {\fr g} \rightarrow {\rm End} (T_{eK}M)$ a Nomizu
map such that 
\begin{enumerate}
\item[(1)] $L_x\pi y -(-1)^{\tilde{x}\tilde{y}}L_y\pi x = 
-\pi [x,y]$ for all $x,y \in {\fr g}$
(i.e.\ $T=0$) and 
\item[(2)] $L_x$ normalizes $Q(eK)\subset {\rm End} (T_{eK}M)$.
\end{enumerate}
Then $\nabla (L)$ is a $G$-invariant quaternionic connection on
$(M,Q)$ and hence $Q$ is\linebreak[3] 1-in\-te\-gra\-ble. 
\end{cor}
For use in \ref{superSec}, 
we give the formula for the Nomizu map
$L^g$ associated to the Levi-Civita connection $\nabla^g$ of a
$G$-invariant pseudo-Riemannian metric $g$ on a homogeneous
supermanifold $M = G/K$. Let $\langle \cdot ,\cdot \rangle = g(eK)$
be the $K$-invariant nondegenerate supersymmetric bilinear form on 
$T_{eK}M$ induced by $g$ (the value of $g$ at $eK$). 
Then $L_x^g \in {\rm End} (T_{eK}M)$, $x\in {\fr g}$, is given by the 
following Koszul type formula:
\begin{equation}\label{superKoszulEqu}
-2\langle L_x^g\pi y,\pi z\rangle = 
\langle \pi [x,y],\pi z\rangle - \langle \pi x, \pi [y,z]\rangle 
- (-1)^{\tilde{x}\tilde{y}}\langle \pi y, \pi [x,z]\rangle 
\, ,\quad x,y,z \in {\fr g}\, .
\end{equation}
\begin{cor}
\label{superqKCor} Let $(M = G/K,Q,g)$ be a homogeneous almost quaternionic
(pseudo-) Hermitian supermanifold and assume that $L_x^g$ 
normalizes
$Q(eK)$ for all $x\in {\fr g}$. Then the Levi-Civita connection
$\nabla^g = \nabla (L^g)$ is a $G$ invariant quaternionic connection
on $(M,Q,g)$ and hence $(M,Q,g)$ is a  quaternionic (pseudo-)
K\"ahler supermanifold if $\dim M_0 >4$. 
\end{cor} 
\end{appendix}
 
\end{document}